%% file: hid.tex
\documentclass{siamart250211}

\usepackage{amsmath}
\usepackage{amssymb}
\usepackage{hyperref}
\usepackage{graphicx}
\usepackage{stmaryrd} 
\usepackage{bbold} 
\usepackage{xcolor}
\usepackage{bm} 
\usepackage{comment}

\usepackage{tikz}

\newcommand{\real}{\mathbb{R}}
\newcommand{\znat}{\mathbb{Z}}
\newcommand{\nat}{\mathbb{N}}

\newcommand{\vA}{\mathsf{A}}
\newcommand{\vB}{\mathsf{B}}
\newcommand{\vC}{\mathsf{C}}
\newcommand{\vE}{\mathsf{E}}
\newcommand{\vQ}{\mathsf{Q}}
\newcommand{\vL}{\mathsf{L}}
\newcommand{\vM}{\mathsf{M}}
\newcommand{\vD}{\mathsf{D}}
\newcommand{\vF}{\mathsf{F}}
\newcommand{\vG}{\mathsf{G}}
\newcommand{\vI}{\mathsf{I}}
\newcommand{\vR}{\mathsf{R}}
\newcommand{\vu}{\mathsf{u}}
\newcommand{\vb}{\mathsf{b}}
\newcommand{\vv}{\mathsf{v}}
\newcommand{\vz}{\mathsf{z}}
\newcommand{\vy}{\mathsf{y}}

\newcommand{\cF}{\mathcal{F}}
\newcommand{\cO}{\mathcal{O}}
\newcommand{\hprod}{{.*}}

\newsiamremark{remark}{Remark}

\DeclareMathSymbol{\dv}{\mathord}{operators}{"3A}
\DeclareMathOperator{\sign}{sign}
\DeclareMathOperator{\supp}{supp}
\DeclareMathOperator*{\argmin}{argmin}

\renewcommand{\mod}[2]{#1\, \text{\scriptsize{\%}}\, #2}

\usetikzlibrary{3d,perspective}

\title{Hyperspectral imaging from diffractograms\thanks{Submitted to the editors
DATE.}
\funding{This work was partially funded by Lumos Imaging Inc. A. M.
and R. M. acknowledge funding from the Office of Naval Research
(N00014-22-1-2014, N000142512122). F. G. V. acknowledges
funding from the National Science Foundation DMS-2136198.}}
\author{Fernando González del Cueto\thanks{Lumos Imaging Inc., Salt Lake
City, UT, USA.}
\and
Fernando Guevara Vasquez\footnotemark[2]\ \thanks{Corresponding author. Department of Mathematics,
University of Utah, Salt Lake City, UT, USA. (\email{fguevara@math.utah.edu})}
\and Apratim~Majumder\footnotemark[2]\ \thanks{Department of Electrical and
Computer Engineering, University of Utah, Salt Lake City, UT, USA.
(\email{rmenon@eng.utah.edu}, \email{apratim.majumder@utah.edu})}
\and Rajesh Menon\footnotemark[2]\ \footnotemark[4]}

\headers{Hyperspectral imaging from diffractograms}{F. Gonz\'alez del Cueto, et al.}

\newlength{\imagesize} 
\begin{document}
\maketitle
\begin{abstract}
We consider an optical system with a periodic point spread function,
such as a camera with a periodic diffractive filter array placed near
the sensor. Such systems produce a monochrome image that we call
``diffractogram'', and that multiplexes the spatial and spectral
information about the scene.  We present a method to efficiently invert
this forward model and recover the hyperspectral scene.  The inverse
problem is formulated as a large scale $\ell_1$-penalized least squares
problem that is solved iteratively with the Alternating Direction Method
of Multipliers (ADMM). A Krylov subspace solver, such as LSQR, may
require tens to hundreds of iterations to approximatively solve the
augmented least-squares problem required in each ADMM iteration. Here we
present an exact direct solver for the augmented least-squares problems
that exploits the periodicity of the point spread function. Each direct
solve costs roughly as much as three LSQR iterations,
eliminating inner iterations entirely. We validate the method on both
synthetic and experimental datasets.
\end{abstract}


\section{Introduction}
\label{sec:intro}

\begin{figure}
\begin{center}
\begin{tikzpicture}

\filldraw[cyan!30,draw=black]
  (0,-1) arc[start angle=-30, end angle=30, radius=2]
  arc[start angle=180-30, end angle=180+30, radius=2];

\def\dx{1/8}
\def\ymin{-1}
\def\ymax{1}
\def\xref{3}
\def\q{6} 
\def\period{1,3,2,4} 
\pgfmathsetmacro{\N}{4*\q}
\pgfmathsetmacro{\dy}{(\ymax-\ymin)/\N}

\def\dfa{\period}
\foreach \i in {2,...,\q} {
\xdef\dfa{\dfa,\period}%
}

\begin{scope}[shift={(1,0)}]
\filldraw[fill=cyan!30,draw=black]
(0,\ymin)
\foreach \x [count=\i from 0] in \dfa  {
    -- ({\dx*\x},{\ymin+\i*\dy}) -- ({\dx*\x},{\ymin+(\i+1)*\dy})
}
-- (0,\ymax) -- cycle;
\end{scope}

\filldraw[fill=lightgray,draw=black] (2,-1) rectangle (2.2,1);

\node[inner sep=0] (obj) at (-4.5,0)
    {\includegraphics[width=3cm]{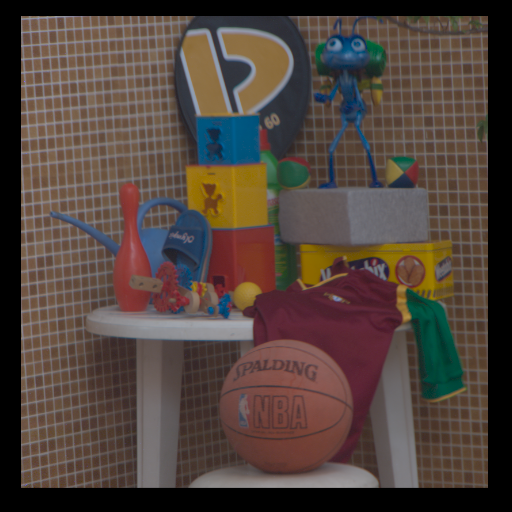}};

\node[inner sep=0] (obj) at (4.5,0)
    {\includegraphics[width=3cm]{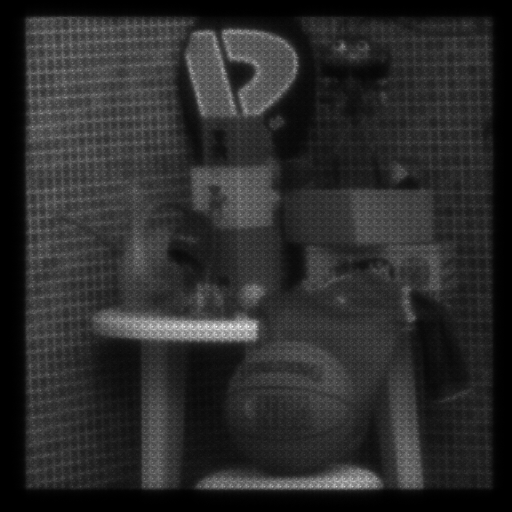}};

\draw (-3,0) -- (0,1) -- (2,0);
\draw (-3,0) -- (0,-1) -- (2,0);
\draw[dashed] (-3,0) -- (2,0);

\node at (-4.5,-2) {\shortstack{object $u(y,\lambda)$\\\small{(RGB
visualization)}}};
\node at (0,-2) {lens};
\node at (1+\dx/2,-2) {DFA};
\node at (2+0.1,-2) {sensor};
\node at (4.5,-2) {diffractogram $v(x)$};
\draw[->,thick] (-3,1.5) .. controls (0,2.4) .. node[midway,above] {forward problem}  (3,1.5);
\draw[->,thick] (3,-2.1) .. controls (0,-3) .. node[midway,below] {inverse problem} (-3,-2.1);
\end{tikzpicture}

\end{center}
 \caption{The forward problem consists in going from a hyperspectral
 object or scene $u(y,\lambda)$ to a diffractogram $v(x)$, which is the
 image obtained at a monochromatic sensor resulting from an optical
 system with a point spread function (PSF) that is periodic in space.
 The inverse problem consists in recovering the hyperspectral data
 $u(y,\lambda)$ from the diffractogram $v(x)$, and assumes knowledge of
 the point spread function of the system, for different wavelengths.
 The images are generated from hyperspectral data in
 \cite{Foster:2006:FMN} and a PSF from \cite{Majumder:2025:HDS}.}
 \label{fig:setup}
\end{figure}
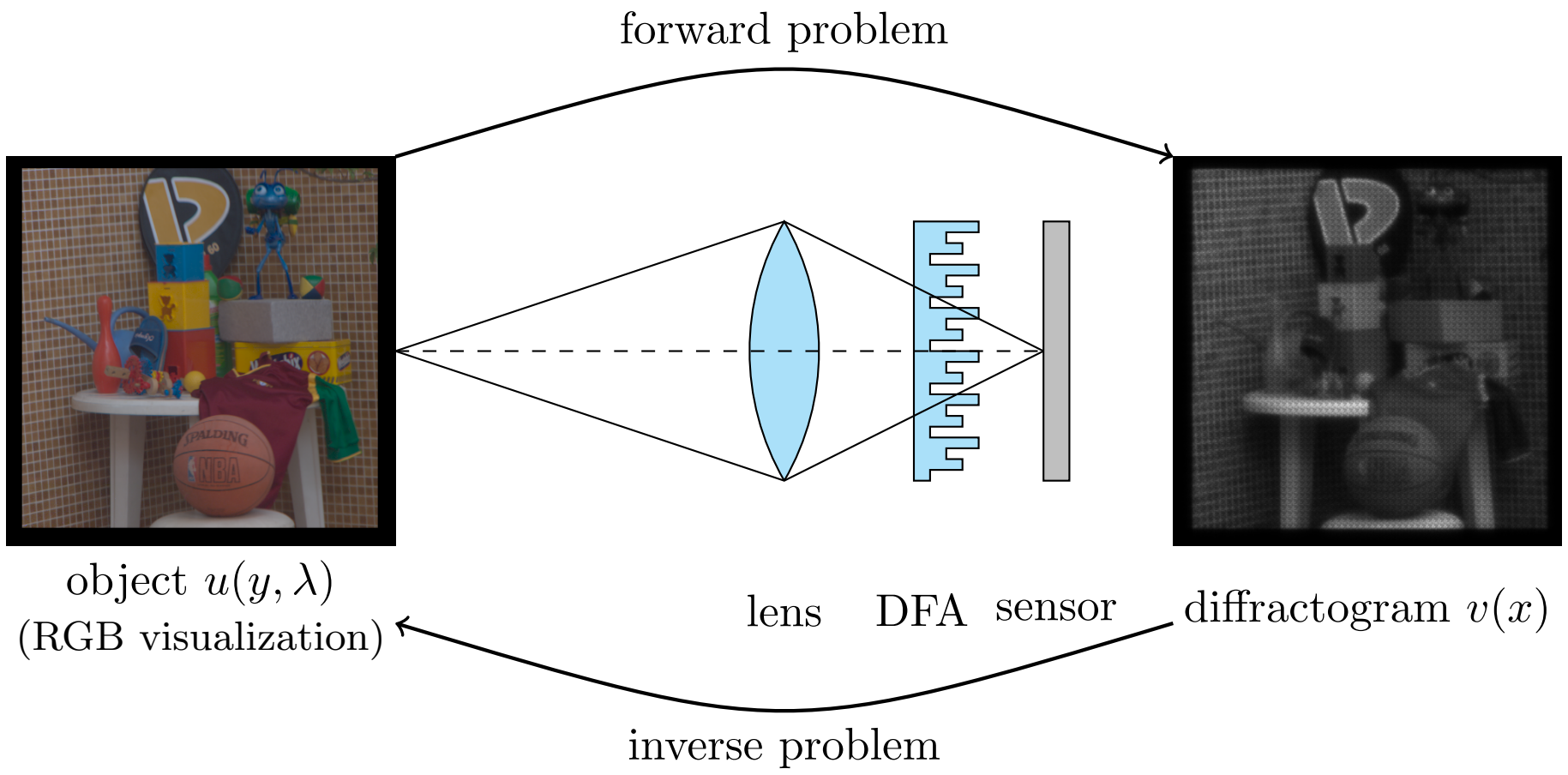
We consider the inverse problem of reconstructing the spatial and
spectral information of an object from its image through an optical
system that has a wavelength dependent 2D periodic mask. This could be
for example a diffractive filter array (DFA) \cite{Majumder:2025:HDS} or
a Bayer-like filter
\cite{Bayer:1976:CIA,Tisserand:2019:CBF,Bian:2024:BHI}. To fix ideas, we
consider the optical system in \cref{fig:setup}. We represent the object
or scene by a function $u(x,\lambda)$ at the \emph{object plane} and the
signal acquired by a monochromatic sensor at the \emph{image plane} by
$v(x)$, where $x = (x_1,x_2) \in \real^2$ and $\lambda>0$ is the wavelength.  The
optical system in \cref{fig:setup} is linear so it may be modeled as
\begin{equation}
 v(x) = \int d\lambda\int dy\,  A(x,y,\lambda) u(y,\lambda),
 \label{eq:3dconv}
\end{equation}
where $A(x,y,\lambda)$ is the point spread function (PSF), i.e. the
image, as a function of $x$, of a pointlike object located at $y$ and for
wavelength $\lambda$. The inverse problem here is to recover
$u(y,\lambda)$ from $v(x)$, assuming $A(x,y,\lambda)$ is known. We
assume the PSF has been appropriately magnified and reflected about the
origin so that the coordinates $x$ in the image plane and $y$ in the
object plane are identical. The periodicity of the imaging system
is assumed to be such that it endows the
PSF with an invariance to translations that are integer multiples of the period
$L = (L_1,L_2)$, where $L_1,L_2>0$. In other words we assume that
for all $x,y\in\real^2$ and wavelength $\lambda$, we have
\begin{equation}
 A(x + d ,y+d,\lambda) = A(x,y,\lambda),~
 \text{for all}~d = (n_1L_1,n_2L_2)~\text{with}~n_1,n_2 \in \znat.
 \label{eq:periodic}
\end{equation}
Such a periodic PSF is illustrated in \cref{fig:cartoonpsf} (a). Space
invariance is the special case where \eqref{eq:periodic} holds for any period
$L$. Indeed, for a particular $y$, we may choose $L=(|y_1|,|y_2|)$ and
$n=-(\sign(y_1),\sign(y_2))$, to see that $A(x,y,\lambda) = A(x-y,0,\lambda)$.
In real optical systems the PSF is space variant to different degrees. Space
invariance is a fundamental assumption in classic Fourier optics, see e.g.
\cite{Goodman:2005:IFO}, as the $y$ integral in \eqref{eq:3dconv} becomes a
convolution. We emphasize that the optical systems we
consider here cannot be modeled with a single convolution. 

The discrete model we use comes from approximating the $\lambda$
integral \eqref{eq:3dconv} by a quadrature formula with sample wavelengths
$\lambda_1,\ldots,\lambda_N$ and weights
$w_1,\ldots,w_N$:
\begin{equation}
 \vv \approx \sum_{i=1}^N w_i \vA_i \vu_i,
 \label{eq:discrete}
\end{equation}
where $\vv$ is a vector of appropriate size representing the signal
$v(x)$ acquired at the image plane, the $\vu_1,\ldots,\vu_N$ are vectors
representing $u(x,\lambda_1),\ldots,u(x,\lambda_N)$ (the unknowns) at
the object plane and the $\vA_i$ are matrices approximating the $y$
integral with $\lambda=\lambda_i$ fixed. With $\vA = [ w_1 \vA_1,
\ldots, w_N \vA_N]$, we formulate the inverse problem as the $\ell_1$
regularized linear least squares problem
\begin{equation}
 \min_{\vu_1,\ldots,\vu_N} \frac{1}{2} \| \vA \vu - \vv \|_2^2 +
 \sum_{j=1}^M \alpha_j \| \vR_j \vu \|_1,
 \label{eq:optim}
\end{equation}
where $\vR_1,\ldots,\vR_M$ are regularization matrices and
$\alpha_1,\ldots,\alpha_M\geq 0$ are the respective regularization
parameters. The vector of unknowns is $\vu = [\vu_1;\ldots;\vu_N]$, with
semicolon denoting vertical concatenation. We focus here on
regularization matrices $\vR_j$ that penalize the total variation
\cite{Rudin:1992:NTV,Dobson:1994:IET} in space, in wavelength and the
$\ell_1$ norm of the reconstructions.

Without exploiting the periodicity of the PSF, even the forward problem
is overwhelming to store and compute. For example if $\vu_i$ and $\vv$
are represented by $256 \times 256$ pixels, each matrix $\vA_i$ is
square of dimension $2^{16} = 65536$, potentially needing to store about
$4.3\times 10^9$ entries. With $N=25$ sample wavelengths, the inverse
problem is a large scale optimization problem involving about $1.6\times
10^6$ unknowns. Fortunately, the PSF is periodic so the optical system
can be approximated using the \emph{periodic circular convolution} or
PCC \cite{Gonzalez:2026:FCC} (\cref{sec:model}). This drastically
reduces the complexity of solving \eqref{eq:optim}, which we do with the
Alternating Direction Method of Multipliers (ADMM)
\cite{Parikh:2014:PA,Beck:2017:FOM}. The computational bottleneck in the
ADMM applied to our problem is the solution of large scale least squares
problems involving the $\vA_i$ (\cref{sec:admm}). Rather than using
Krylov subspace methods such as LSQR \cite{Paige:1982:LSQR}, we solve
the least squares problems with a direct method (\cref{sec:direct}),
which relies on the unitary equivalence between a PCC and a
block diagonal matrix \cite{Gonzalez:2026:FCC}. Although the direct
method imposes limits on the structure of the regularization matrix
$\vR$, it allows for $\ell_1$ and total variation regularization on both
space and wavelength. Before presenting numerical results we briefly
review in \cref{sec:rgb} a method to approximate how a hyperspectral
image would be perceived by the human eye. Then in \cref{sec:synth} we
present reconstructions using a reference hyperspectral data sets and an
experimentally acquired PSF.  Reconstructions with experimentally
acquired diffractograms and PSFs are shown in \cref{sec:realdata}. We
conclude with a summary and perspectives in \cref{sec:summary}.

\subsection{Related work}
We point to \cite{Hagen:2013:RSS} for a review of snapshot hyperspectral
imaging techniques. The method we present is tailored to the particular
technique introduced in \cite{Majumder:2025:HDS}, which uses a periodic
DFA close to the sensor. Although periodic DFAs are also used in
\cite{Wang:2014:CSS,Wang:2014:CSB,Wang:2015:UHS,Wang:2018:CMV}, the DFA
is placed relatively far from the sensor, so the PSF is no longer
periodic in the sense of \cref{eq:periodic}. Another example of a
hyperspectral imaging system with non-periodic PSF is the Scotch tape
camera \cite{Malone:2023:DSS}.  Without PSF periodicity, scalability
(in terms of image size or computational cost) is difficult to
achieve. Our method could be applied to other hyperspectral imaging
techniques with a periodic PSF resulting from a periodic array of
optical elements.  One example is the Bayer filter, which consists of a
periodic array of pigment filters deposited on the sensor
\cite{Bayer:1976:CIA,Tisserand:2019:CBF,Bian:2024:BHI}. Another example
is Fabry-Pérot interferometers built on the sensor \cite{Tack:2012:CHS}.
It may also be possible to obtain periodic PSFs by using a lenslet
array, fibers or mirrors as in Integral Field Spectroscopy, see e.g.
\cite{Bowen:1938:ISD,Arribas:1991:OFS,Bacon:1995:SHS}, however the
period is not necessarily rectangular as we assume here.

From an optimization perspective, we point to another efficient method
for minimizing the total variation of a multichannel image in
\cite{Wang:2008:NAM}.  The ADMM based reconstruction we present here is
more flexible as it allows for other $\ell_1$ based regularization
terms. Finally blind compressed sensing can be used for hyperspectral
imaging and does not require the PSF as we do in our case
\cite{Rajwade:2013:CHI}.

\section{Discretization of periodic kernels as periodic circular
convolutions} 
\label{sec:model} 
Since we discretize the integral in wavelength \eqref{eq:3dconv} with
a quadrature (\cref{sec:intro}), we can focus on the spatial integral
and momentarily drop the wavelength dependency. Hence our
goal here is to discretize integrals of the kind
\begin{equation}
 v(x) = \int dy\, A(x,y) u(y),
 \label{eq:onlyx}
\end{equation}
where $A(x,y)$ is a periodic kernel with period $L=(L_1,L_2)$. In
\cref{sec:padding} we make assumptions on the supports of the kernel and
of the signal so that integrals can be carried out on a finite region
$R$. We introduce a modified kernel $\widetilde A$ that in addition to
being an $L-$periodic kernel, both $\widetilde A(x,\cdot)$ and
$\widetilde A(\cdot,y)$ are periodic functions for any fixed $x,y\in
\real^2$. This allows to calculate the integrals \eqref{eq:onlyx} on
$R$, with an approach reminiscent to zero-padding for regular
convolution operations, see e.g. \cite{Press:2007:NR}. In
\cref{sec:sampling}, we choose the region $R$  and a
uniform discretization commensurate with the kernel period $L$ to sample
the kernel $\widetilde A$, preserving its periodicity. Then
\eqref{eq:onlyx} can be approximated by a periodic circular convolution
\cite{Gonzalez:2026:FCC}, a linear operation that is unitarily
equivalent to a block diagonal matrix and thus can be operated with
efficiently (\cref{sec:discretization}).

\begin{figure}
 \begin{center}
  \begin{tabular}{cc}
  \begin{tikzpicture}[x=0.8cm,y=0.8cm]
   \draw[->,thick] (-2,0) -- node[at end,above]{$x_1$} (4,0);
   \draw[->,thick] (0,-2) -- node[at end,left]{$x_2$} (0,4);
   \foreach \i in {-1,0,1,2,3} {
    \draw[red] (-1.3,\i) -- (3.3,\i);
    \draw[red] (\i,-1.3) -- (\i,3.3);
   }
   \node[below right] at (1,0) {$1$};
   \node[above left] at (0,1) {$1$};
   \node[below right] at (0,0) {$0$};
   \node[above left] at (0,0) {$0$};

   \node[minimum size= 0.3,inner sep=0,draw=blue,text=black] at (0.7,0.1) {\includegraphics[width=5mm]{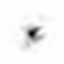}};
   \node[minimum size= 0.3,inner sep=0,draw=blue,text=black] at (1.7,1.1) {\includegraphics[width=5mm]{cartoonpsf1}};
   \node[minimum size= 0.3,inner sep=0,draw=blue,text=black] at (2.7,1.1) {\includegraphics[width=5mm]{cartoonpsf1}};
   \node[minimum size= 0.3,inner sep=0,draw=blue,text=black] at (0.5,1.5) {\includegraphics[width=5mm]{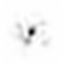}};
   \node[minimum size= 0.3,inner sep=0,draw=blue,text=black] at (1.5,2.5) {\includegraphics[width=5mm]{cartoonpsf2}};
  \end{tikzpicture}
   & 
  \begin{tikzpicture}[x=0.8cm,y=0.8cm]
  \draw[->,thick] (-1,0) -- node[at end,above]{$x_1$} (5,0);
  \draw[->,thick] (0,-1) -- node[at end,left]{$x_2$} (0,5);
  \draw[thick] (0,0) rectangle (4,4);
  \fill[black!20] (1,1) rectangle (3,3);
  \draw[dashed] (1,-1) -- (1,5);
  \draw[dashed] (3,-1) -- (3,5);
  \draw[dashed] (4,-1) -- (4,5);
  \draw[dashed] (-1,1) -- (5,1);
  \draw[dashed] (-1,3) -- (5,3);
  \draw[dashed] (-1,4) -- (5,4);
  \draw[<->] (0,-0.7) -- node[midway,below]{$a_1$} (1,-0.7);
  \draw[<->] (3,-0.7) -- node[midway,below]{$a_1$} (4,-0.7);
  \draw[<->] (-0.7,0) -- node[midway,left]{$a_2$} (-0.7,1);
  \draw[<->] (-0.7,3) -- node[midway,left]{$a_2$} (-0.7,4);
  \node[below right] at (4,0) {$R_1$};
  \node[above left]  at (0,4) {$R_2$};
   \node[below right] at (0,0) {$0$};
   \node[above left] at (0,0) {$0$};
  \end{tikzpicture}
   \\
   (a) & (b)
  \end{tabular}
 \end{center}
 \caption{(a) An illustration of the PSF periodicity with $L_1=L_2=1$. The PSF
 is invariant to translations by an integer multiple of the period (red
 lines).
 The PSF with center $y'=(0.7.0.1)$ is up to a translation
 the same as those in $y'+(1,1)$ and $y'+(2,1)$, but different from the PSF
 with center $y''=(0.5,1.5)$ which itself is identical up to translation to
 the one at $y''+(1,1)$. Assumption \cref{eq:cs:A} means that the PSF
 $A(\cdot,y)$ is supported inside a $2a_1 \times 2a_2$ box centered at $y$, shown
 in blue. (b) Applying the kernel $A$ to a signal $u$ supported inside the gray
 box (assumption \cref{eq:cs:u}), gives a signal $v$ supported in $[0,R_1]
 \times [0,R_2]$.}
 \label{fig:cartoonpsf}
\end{figure}

\subsection{Zero padding}
\label{sec:padding}
We assume that the kernel $A$ is supported near the diagonal, i.e.
 there are $a_1,a_2>0$ such that 
\begin{equation}
 \supp A \subset \{x,y \in \real^2 ~|~ |x_i-y_i| < a_i,~i=1,2 \}.
 \label{eq:cs:A}
\end{equation}
Physically, the support assumption means that the point spread function
$A(\cdot,y)$ is spatially localized to the box of  size $2a_1 \times
2a_2$ centered at $y$. Note that $2a_1$ and $2a_2$ are unrelated to
the period dimensions $L_1,L_2$, i.e. they could be larger or smaller. Similarly,
we assume that the signal $u$ is smooth and compactly supported in the
region
\begin{equation}
 \supp u \subset (a_1,R_1-a_1) \times (a_2,R_2-a_2),
 \label{eq:cs:u}
\end{equation}
for some $R_1,R_2>0$. The supports of the kernel and the signal are
chosen so that the result $v$ of applying the kernel $A$ to $u$ is supported
in $R \equiv (0,R_1) \times (0,R_2)$, see \cref{fig:cartoonpsf} (b).
With these assumptions on the support of the kernel and the signal, we
can replace the kernel $A$ by a periodized kernel $\widetilde A$, in a
manner similar to zero-padding for regular convolution problems, see
e.g. \cite{Press:2007:NR}.

\begin{lemma}
\label{lem:periodized}
Assume that $R_j > 2a_j$, $j=1,2$. Then for $x \in R$ we have
\begin{equation}
 v(x) = \int_R dy\, \widetilde A (x,y) u(y),
 \label{eq:atilde:int}
\end{equation}
where $\widetilde A$ is the periodization of the kernel $A$ defined by
\begin{equation}
 \widetilde A(x,y) = \sum_{i_1,i_2 \in \znat} A(x,y+(i_1 R_1,i_2 R_2)).
 \label{eq:atilde}
\end{equation}
\end{lemma}
\begin{proof}
First notice that the assumption $R_j>2a_j$, $j=1,2$, guarantees that
the supports of the terms in the series \eqref{eq:atilde} do not intersect. Indeed
for a fixed $y$, we have
\begin{equation}
 \supp A(\cdot,y+(i_1R_1,i_2R_2)) \subset y + (i_1R_1,i_2R_2) + (-a_1,a_1)
 \times (-a_2,a_2),
\end{equation}
which are disjoint subsets of $\real^2$ for different $(i_1,i_2)$. This
can be seen by the inequalities $a_l + i_l R_l < -a_l + (i_l+1)R_l$,
$l=1,2$. We conclude that for any $x,y \in \real^2$, the series
\eqref{eq:atilde} has at most one nonzero term.

For some $x \in R$, we see that
\begin{equation}
 \int_R dy\, \widetilde A (x,y) u(y) = \sum_{i_1,i_2 \in \znat} \int_R
 dy\, A(x,y+(i_1 R_1,i_2 R_2)) u(y).
 \label{eq:atilde:2}
\end{equation}
Because of \eqref{eq:cs:u}, we only need to consider $y \in (a_1,R_1-a_1)
\times (a_2,R_2-a_2)$. We claim that the integral on the right hand side
of \eqref{eq:atilde:2} vanishes for all $(i_1,i_2) \neq (0,0)$. Indeed
if either one of $i_1$ or $i_2$ is nonzero, we must have that $|x_l -
(y_l+i_lR_l)|>a_l$ for the corresponding $l=1,2$. Since only the term
$(i_1,i_2) = (0,0)$ matters in the sum in \eqref{eq:atilde:2}, we get
the desired result
\begin{equation}
 \int_R dy\, \widetilde A (x,y) u(y)  = \int_R dy\, A(x,y) u(y) = \int
 dy\, A(x,y) u(y).
\end{equation}
\end{proof}

\begin{lemma}
 \label{lem:atilde}
 If we further assume that $R_j$ is an integer multiple of $L_j$ for
 $j=1,2$, then the periodized kernel $\widetilde A(x,y)$ as defined in
 \cref{lem:periodized} satisfies the following.
 \begin{itemize}
  \item[(i)] The kernel $\widetilde A(x,y)$ is a $(L_1,L_2)-$periodic
  kernel, see
  \cref{eq:periodic} for the definition.
  \item[(ii)] For any $x,y \in \real^2$, $\widetilde A(\cdot,y)$ and
  $\widetilde A(x,\cdot)$ are periodic functions with unit cell $R$.
 \end{itemize}
\end{lemma}
\begin{proof}
First note that if $A(x,y)$ is an $(L_1,L_2)-$periodic kernel, then so
is the kernel $A(x,y+d)$ for any shift $d$ we choose.
Since $\widetilde A(x,y)$ is defined as a sum of shifted
$(L_1,L_2)-$periodic kernels, it must also be $(L_1,L_2)-$periodic,
proving (i).

By its definition \eqref{eq:atilde},  $\widetilde A(x,y)$ is periodic
with period $(R_1,R_2)$ as
a function of $y$ for any fixed $x$. Indeed for any $j_1,j_2 \in \znat$ we have
$\widetilde A (x, y + (j_1 R_1,j_2 R_2)) = \widetilde A (x,y)$. To show
periodicity in $x$, we use (i) and that $R_l$ is
an integer multiple of $L_l$, $l=1,2$, to obtain
\begin{equation}
\widetilde A(x+(j_1 R_1,j_2 R_2),y) = \widetilde A(x,y+(j_1 R_1,j_2
R_2)) = \widetilde A(x,y),
\end{equation}
which proves (ii).
\end{proof}

\subsection{Sampling the kernel}
\label{sec:sampling}
We sample the region $R$ with the points $(i_1h_1,i_2h_2)$,
$i_1=0,\ldots,n_1-1$ and $i_2=0,\ldots,n_2-1$. We choose the region $R$ and the
discretization steps $h_1,h_2$ commensurate with the period $L$ of the
kernel, i.e. there are $p_l,q_l \in \nat$ such that
\begin{equation}
 R_l = q_l L_l ~\text{and}~ h_l = \frac{R_l}{n_l} =
 \frac{L_l}{p_l},~\text{for}~l=1,2.
 \label{eq:commensurate}
\end{equation}
Since the kernel $\widetilde A$ is defined on $R\times R$, we may
organize the samples of $\widetilde A$ as a $4-$th order tensor
$\widetilde \vA$ of size $n_1 \times n_2 \times n_1 \times n_2$ with entries
\begin{equation}
 \widetilde \vA[i_1,i_2,j_1,j_2] = \widetilde
 A(i_1h_1,i_2h_2,j_1h_1,j_2h_2), 
 \label{eq:4tensordef}
\end{equation}
where  $i_l = 0,\ldots,n_l-1$ and $j_l=0,\ldots,n_l-1$ for $l=1,2$.
Hereinafter tensors are denoted with sans serif letters, with 
indexes starting at $0$ and understood modulo the
respective dimension. For instance
\begin{equation}
 \widetilde \vA [ i_1,i_2,j_1,j_2] = \widetilde \vA [ \mod{i_1}{n_1},
 \mod{i_2}{n_2}, \mod{i_1}{n_1}, \mod{i_2}{n_2} ],
 \label{eq:modulo}
\end{equation}
for any $i_1,i_2,j_1,j_2\in \znat$, where $\mod{i}{n} \in \{0,\ldots,n-1\}$ is the remainder in the Euclidean division of $i$ by $n$. The tensor notation and the
periodicity of $\widetilde A(\cdot,y)$ and $\widetilde A(x,\cdot)$ for
any fixed $x,y$, play well together making the tensor definition
\eqref{eq:4tensordef} valid for all $i_1,i_2,j_1,j_2 \in
\znat$. Indeed, by Euclidean division we can write $i_l = i'_l + i''_l
n_l$ and $j_l = j'_l + j''_l n_l$, where $i'_l = \mod{i_l}{n_l}$, $j'_l =
\mod{j_l}{n_l}$, and $i''_l,j''_l \in \znat$, for $l=1,2$. Now the period
of $\widetilde A(x,\cdot)$ and $\widetilde A(\cdot,y)$ is $R_l = h_l n_l$, for
$l=1,2$. Therefore we see that
\begin{equation}
\begin{aligned}
 \widetilde A  (i_1h_1,i_2h_2,j_1h_1,j_2h_2) 
&=\widetilde A(i_1'h_1,i_2'h_2,j_1'h_1,j_2'h_2)\\
&=\widetilde \vA[i_1',i_2',j_1',j_2']
=\widetilde \vA[i_1,i_2,j_1,j_2],
\end{aligned}
\end{equation}
where the last equality comes from \eqref{eq:modulo}. 

The crucial observation is that the $4-$th order tensor $\widetilde \vA$
is a $(p_1,p_2)-$periodic (discrete) kernel, a property a property
similar to \eqref{eq:periodic}, that is for any
$i_1,i_2,j_1,j_2\in \znat$ and $m_1,m_2,m_3,m_4 \in \znat$ 
\begin{equation}
 \widetilde \vA[i_1 + m_1 p_1,i_2 + m_2 p_2,j_1 + m_3 p_1,j_2 + m_4 p_2]
 = \widetilde \vA[i_1,i_2,j_1,j_2].
 \label{eq:periodic:discrete}
\end{equation}
This readily follows by writing the left hand side as an evaluation of
the kernel $\widetilde A$, and then using that
$p_l h_l = L_l$, $l=1,2$ and the $(L_1,L_2)-$periodicity of the kernel
$\widetilde A$ (\cref{lem:atilde}).

\subsection{An experimentally acquired kernel}
\label{sec:calibration}
The discretization we present here is motivated by the ``calibration'' or PSF
acquisition procedure in \cite{Majumder:2025:HDS}. In that case, a diffractive
filter array that is placed near the sensor has a period that is commensurate
with a $8\times 8$ sensor pixel region, so $h_1=h_2$ is fixed to be the pixel
pitch ($2.9\mu$m). To acquire $\widetilde \vA[:,:,j_1,j_2]$, a pinhole
restricting a narrowband diffraction limited focal spot acts as a point-like
object that is moved to scan $64$ different values for $j_1, j_2 \in
\{0,\ldots,7\}$, covering a small $8\times 8$ pixel patch or ``calibration
region''. The resulting $4-$th order tensor has all the information needed to
apply $\widetilde \vA$, near the calibration region.

Two examples of an experimentally acquired $4-$th order tensor
$\widetilde \vA$ are presented in \cref{fig:psf}, corresponding to $500$
nm and $620$ nm narrowband beams. Each $4-$th order tensor was originally of
dimensions $48 \times 48 \times 8 \times 8$, but is cropped to a $20
\times 20 \times 8 \times 8$ tensor for visualization purposes. The
visualization in \cref{fig:psf} consists of $8\times 8$ tiles, each
displaying $\widetilde \vA [:,:,j_1,j_2]$ for $j_1,j_2 \in
\{0,\ldots,7\}$. In other words, each tile is a cropping of a sensor
image of a point like object at different positions and at a particular
wavelength. The images was gamma corrected and color was added to
approximate the human perception of the wavelengths. This is
done via the same procedure to visualize hyperspectral images from
\cite{Foster:2019:HIC} and that we detail later in \cref{sec:rgb}.
Notice how the center of mass of each individual PSF moves to the bottom
(resp. right) when $j_1$ (resp. $j_2$) is increased. To
illustrate the drastic wavelength dependence of this PSF, we also
include zoomed in versions of $\widetilde \vA[:,:,4,3]$ for the same
wavelengths, in gamma corrected grayscale.

\newcommand{\plotpsf}[1]{%
 \setlength{\imagesize}{0.3\textwidth}
 \begin{tikzpicture}[x=\dimexpr\imagesize/8\relax,
                     y=-\dimexpr\imagesize/8\relax]
  \node[anchor=north west, inner sep=0] at (0,0)
    {\includegraphics[width=\imagesize]{#1}};
  \node at (-0.25,0.5) {$0$};
  \node at (-0.25,4.5) {$4$};
  \node at (-0.25,7.5) {$7$};
  \node at (0.5,-0.3) {$0$};
  \node at (4.5,-0.3) {$4$};
  \node at (7.5,-0.3) {$7$};
  \foreach \i in {1,2,3,4,5,6,7} {
   \draw[thick,blue] (0,\i)--(8,\i);
   \draw[thick,blue] (\i,0)--(\i,8);
  }
  \draw[->,thick] (0,0) -- node[at end,above]{$j_2$} (8.5,0);
  \draw[->,thick] (0,0) -- node[at end,left]{$j_1$} (0,8.5);
  \end{tikzpicture}%
}

\begin{figure}
 \begin{center}
  \begin{tabular}{cc}
   500 nm & 620 nm\\
   \plotpsf{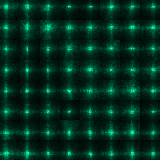} & \plotpsf{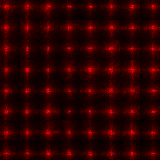}\\
   \includegraphics[width=0.175\textwidth]{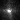} &
   \includegraphics[width=0.175\textwidth]{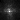} \\
   $j_1=4$, $j_2=3$ & 
   $j_1=4$, $j_2=3$ \\
  \end{tabular}

 \end{center}

 \caption{An example of experimentally acquired PSF for wavelengths
 $500$ nm and $620$ nm. The first row is a visualization of the
 corresponding fourth order tensors $\widetilde \vA_{500}$ and
 $\widetilde \vA_{620}$. We display 64 slices of $\widetilde
 \vA_\lambda [:,:,j_1,j_2]$  for $j_1,j_2 \in \{0,\ldots,7\}$ as $20\times 20$ pixel tiles separated
 by blue lines, for $\lambda \in \{500,620\}$ nm. The particular
 $j_1=4,j_2=3$ tiles are presented in grayscale in the second row, to
 highlight the wavelength dependence of the PSF.} 
 \label{fig:psf}
\end{figure}

\subsection{The discretization}
\label{sec:discretization}
Let us now discretize the integral in \eqref{eq:atilde:int} by the Riemann
sums
\begin{equation}
 v(i_1 h_1, i_2 h_2) \approx \sum_{j_1=0}^{n_1-1}\sum_{j_2=0}^{n_2-1} \widetilde A(i_1 h_1,
 i_2 h_2, j_1 h_1, j_2 h_2) u(j_1 h_1, j_2 h_2),
 \label{eq:riemann}
\end{equation}
for $i_l=0,\ldots,n_l-1$, $l=1,2$.
If we define the $2-$nd order $n_1 \times n_2$ tensors $\vu[i_1,i_2]=u(i_1 h_1,i_2h_2)$ and
$\vv[i_1,i_2] \approx v(i_1 h_1,i_2h_2)$ for $i_l \in
\{0,\ldots,n_l-1\}$, calculating the Riemann sums \eqref{eq:riemann}
is the tensor contraction operation
\begin{equation}
\vv[i_1,i_2] = \sum_{j_1=0}^{n_1-1}\sum_{j_2=0}^{n_2-1}
\widetilde \vA[i_1,i_2,j_1,j_2] \vu[j_1,j_2], 
\end{equation}
for $i_l=0,\ldots,n_l-1$, $l=1,2$.  Since the (discrete) kernel
$\widetilde \vA$ is $(p_1,p_2)-$periodic \cref{eq:periodic:discrete},
this (linear) tensor contraction is unitarily equivalent to a block
diagonal matrix of size $n_1n_2 \times n_1 n_2$ with $q_1q_2$ blocks,
each of dimension $p_1p_2 \times p_1 p_2$, where $n_l=p_lq_l$, $l=1,2$.
Crucially, the unitary equivalence is
given by appropriate Discrete Fourier Transforms that only depend on
$p_l, q_l$, $l=1,2$, which enables efficient operations with such
matrices \cite{Gonzalez:2026:FCC}. For conciseness we hereinafter use
the notation $|p| = p_1p_2$, $|q|=|q_1q_2|$ and $|n|=n_1n_2$. We note
that only $|q| |p|^2 = |p| |n|$ scalars need to be stored in memory to
represent $\widetilde \vA$, rather than the $|n|^2$ scalars needed with
a naïve approach.  Similarly, solving least squares problems with
$\widetilde \vA$ is relatively efficient, as we see later in
\cref{sec:direct}. 

\begin{remark}
\label{rem:trapz}
A detailed analysis of the discretization error in \eqref{eq:riemann}
is out of the scope of this paper. However the integrand is compactly
supported in $R$, and therefore the quadrature \eqref{eq:riemann} is a
composite trapezoidal rule. At least for one dimensional problems the
trapezoidal rule with smooth compactly supported integrands is
exponentially convergent, see e.g. \cite{Trefethen:2014:ECT}.  
\end{remark}

\section{Alternating Direction Method of Multipliers}
\label{sec:admm}
We solve \eqref{eq:optim} using the Alternating Direction Method of
Multipliers (ADMM), see \cite{Parikh:2014:PA,Beck:2017:FOM}. We briefly
review the ADMM method, with the main purpose being to expose the
computational bottlenecks and set the stage for a more detailed
explanation of how we deal with them in \cref{sec:direct}. By
introducing the auxiliary variables $\vz_1,\ldots,\vz_M$, the
optimization problem \eqref{eq:optim} is equivalent to the linearly
constrained problem 
\begin{equation}
\min_{\vu,\vz_1,\ldots,\vz_M} \frac{1}{2} \|\vA \vu - \vv\|_2^2 + \sum_{j=1}^M \alpha_j
\|\vz_j\|_1, \quad
\text{s.t.}~ \vz_j = \vR_j \vu~\text{for}~j=1,\ldots,M.
 \label{eq:optim3}
\end{equation}
For some penalty parameter $\rho>0$, consider the optimization problem
\begin{equation}
\min_{\vu,\vz_1,\ldots,\vz_M} \frac{1}{2} \|\vA \vu - \vv\|_2^2 + \sum_{j=1}^M \alpha_j
\|\vz_j\|_1 + 
\frac{\rho}{2} \|\vz_j - \vR_j \vu\|_2^2, \quad
\text{s.t.}~ \vz_j = \vR_j \vu~\text{for}~j=1,\ldots,M.
 \label{eq:optim4}
\end{equation}
Problem \eqref{eq:optim4} is still equivalent to \eqref{eq:optim3}
because the additional terms in the objective function vanish for
feasible points. The augmented Lagrangian is the Lagrangian of
\eqref{eq:optim4}, namely
\begin{equation}
 L_\rho(\vu,\vz,\vy) = \frac{1}{2} \|\vA \vu - \vv\|_2^2 + \sum_{j=1}^M
 \alpha_j \|\vz_j\|_1 + 
\frac{\rho}{2} \|\vz_j - \vR_j \vu\|_2^2 + \vy_j^T(\vR_j \vu - \vz_j),
\label{eq:auglag}
\end{equation}
where $\vy_j$ is the Lagrange multiplier corresponding to the linear
constraint $\vz_j = \vR_j \vu$ and we used the shorthands $\vy =
[\vy_1;\ldots;\vy_M]$ and $\vz = [\vz_1;\ldots;\vz_M]$. Notice that
$\vy_j$ and $\vz_j$ have the same dimensions. The ADMM method consists
in alternating between minimizing $L_\rho$ with respect to $\vu$, with
the other variables fixed and minimizing $L_\rho$ with respect to
$\vz_1,\ldots,\vz_M$, with the other variables fixed. The algorithm is summarized in
\cref{alg:admm}. Note from \eqref{eq:auglag} that the larger the
parameter $\rho$ is, the more we penalize auxiliary variables $\vz_j$
far from $\vR_j \vu$. 

\begin{algorithm}[h]
\caption{ADMM algorithm.}
\label{alg:admm}
\begin{tabbing}
xxx \= xxx \= xxx\kill
\textbf{Inputs:} Initial guesses $\vu^{(0)},\vz^{(0)},\vy^{(0)}$\\
for \> $k=0,1,\ldots$\\
\> (a) $\vu^{(k+1)} = \argmin_{\vu} L_\rho(\vu,\vz^{(k)},\vy^{(k)})$\\
\> (b) $\vz^{(k+1)} = \argmin_{\vz} L_\rho(\vu^{(k+1)},\vz,\vy^{(k)})$\\
\> (c) $\vy^{(k+1)}_j = \vy^{(k)}_j + \rho (\vR_j \vu^{(k+1)} -
\vz^{(k+1)}_j)$, $j=1,\ldots,M$.
\end{tabbing}
\end{algorithm}

Step (a) in \cref{alg:admm} amounts to solving the linear least squares
problem
\begin{equation}
 \min_{\vu} \frac{1}{2} \| \vA \vu - \vv \|_2^2 + \frac{\rho}{2} \sum_{j=1}^M
 \|\vR_j \vu - (\vz_j^{(k)} - \rho^{-1} \vy_j^{(k)}) \|_2^2,
 \label{eq:admm:a}
\end{equation}
and is by far the computational bottleneck in the ADMM for solving
\eqref{eq:optim}. Even for relatively small images, the matrix $\vA$ is
dense and too large to store in memory, so conventional direct methods (e.g. the
QR factorization) are impractical. Iterative methods such LSQR
\cite{Paige:1982:LSQR} are a possibility, assuming we can efficiently
compute matrix vector products with $\vA$, $\vR_j$ and their transposes.
For our particular problem, $\vA$ is a collection of PCC operations and
the $\vR_j$ are difference operators so efficient matrix vector products
are available. One advantage of iterative methods is that a precise
solution of the least squares problems may not be necessary because the
least squares problems are solved as part of the ADMM iteration.  The
computational time may thus be reduced by adjusting the number of
iterations/precision for the least squares solver. This could also be drawback,
since it introduces parameters to tune, e.g. a maximum number of
iterations and/or a desired
precision.  Rather than solving \eqref{eq:admm:a} iteratively, we
present a direct method in
\cref{sec:direct}, which relies on the block diagonal representation of
periodic circular convolutions \cite{Gonzalez:2026:FCC}. With this
approach a least squares solve has a cost comparable to that of three
iterations of LSQR.

By adding appropriate quadratic terms in $\vy^{(k)}_j$ to the augmented
Lagrangian \cref{eq:auglag}, we see that Step (b) in \cref{alg:admm} is
equivalent to
\begin{equation}
 \vz^{(k+1)} = \argmin_{\vz}  \sum_{j=1}^M \alpha_j \|\vz_j\|_1 +
 \frac{\rho}{2} \| \vR_j \vu^{(k+1)} - \vz_j
 + \rho^{-1} \vy^{(k)}_j \|_2^2,
 \label{eq:zproblem}
\end{equation}
which amounts to a proximal mapping involving the $1-$norm and
can be found explicitly in terms of inexpensive soft
thresholding operations, see e.g. \cite[\S 6]{Beck:2017:FOM}:
\begin{equation}
 \vz^{(k+1)}_j = \text{shrink}( \vR_j \vu^{(k+1)} + \rho^{-1}
 \vy^{(k)}_j, \alpha_j / \rho ),~j=1,\ldots,M,
 \label{eq:soft}
\end{equation}
where the operation $\text{shrink}(x,t) \equiv \sign(x) \max \{ |x| - t,0\}$ for
$x\in\real$ and $t>0$, is applied componentwise to the entries of the
vectors involved.

The idea behind Step (c) of \cref{alg:admm} is that the updated Lagrange
multipliers $\vy^{(k+1)}$ are chosen so that $0$ is in the subgradient
of $L_\rho$ with respect to $\vz$ and evaluated at $\vu^{(k+1)}$,
$\vz^{(k+1)}$. The parameter $\rho$ acts as a step size for Step (c).
For more details we refer to \cite[\S 15]{Beck:2017:FOM}.  For our
discussion, the important observation is that the computational cost of
this step is negligible compared to the cost of Step (a).

\section{Direct solution of the augmented least squares problem}
\label{sec:direct}

We show how to solve the least squares problem
\eqref{eq:admm:a}, which is the main computational bottleneck of the
ADMM algorithm applied to solving \eqref{eq:optim}. The direct method
we present here exploits the structure of $\vA$ and the regularization
matrices $\vR_j$, which are defined in \cref{sec:R}. The normal
equations for \eqref{eq:admm:a} are formulated in \cref{sec:normal}, but
even exploiting periodicity the cost of solving them could remain prohibitive. In
\cref{sec:smw} we show that the Sherman-Morrison-Woodbury formula can
bring down the cost significantly. The computational
cost of the direct method is compared to that of one LSQR iteration in
\cref{sec:estimate} by estimating the number of floating point
operations (flops) involved and finally, by wall-clock timing comparisons.

\subsection{The regularization matrices}
\label{sec:R}
We use $M=3$ regularization matrices so that $\|\vR_1 \vu\|_1$ computes
the spatial total variation, $\|\vR_2 \vu\|_1$  the total variation in
wavelength and $\vR_3 = \vI_{N|n|}$ so that $\|\vR_3 \vu\|_1 =
\|\vu\|_1$. We only consider total variation from a purely discrete
perspective. Following the discretization scheme in \cref{sec:discretization}, we
approximate the unknown signal $u(x,\lambda_i)$ at the $i-$th
wavelength, with a $2-$nd order tensor $\vu_i$ of size $n_1 \times n_2$,
for $i=1,\ldots,N$. In a slight abuse of notation, $\vu_i$ may also
represent the length $|n|$ vector representation of $\vu_i$, assuming column
major ordering. 

To compute the spatial total variation we first need to compute the
differences in the first and second spatial dimensions with the
operators $\vD_1$ and $\vD_2$, where $\vD_l$ is the \emph{periodic}
finite difference operator in the $l-$th dimension for a $n_1 \times
n_2$ tensor, $l=1,2$. In other words, $\vD_1$ and $\vD_2$ are the linear
operators such that
\begin{equation}
 \vD_1 \vv = \vv[1 \dv n_1,\dv] - \vv 
 ~\text{and}~
 \vD_2 \vv = \vv[\dv,1\dv n_2] - \vv,
 \label{eq:diffop}
\end{equation}
for some $n_1 \times n_2 \times N$ tensor $\vv$. Since we use  periodic
tensor indexing, the differences in \eqref{eq:diffop}
wrap around.  For instance, the difference $\vv[n_1,\dv] -
\vv[n_1-1,\dv]$ is the same as $\vv[0,\dv] - \vv[n_1-1,\dv]$. The choice
of periodic difference operators is crucial because it is compatible
with PCC operations, and is inspired by a similar choice in
\cite{Wang:2008:NAM}. Applying $\vD_l$ to each $\vu_j$ in
$[\vu_1;\vu_2;\ldots;\vu_N]$ can be easily expressed with the Kronecker
product $\vI_N \otimes \vD_l$ where $l=1,2$ and $\vI_N$ is the $N\times
N$ identity matrix. This gives the first regularization matrix $\vR_1$
of size $2N|n| \times N|n|$
\begin{equation} 
\vR_1 = 
\begin{bmatrix} 
\vI_N \otimes \vD_1\\ 
\vI_N \otimes \vD_2 
\end{bmatrix}.
\end{equation}

The total variation in wavelength is computed from the differences
between consecutive wavelengths, which can be expressed in terms of
$D^{(N)}$,  the $(N-1) \times N$ difference operator given by
\begin{equation}
 \vD^{(N)} = \begin{bmatrix}
 -1 & 1\\
 & -1 & 1\\
 && \ddots & \ddots\\
 &&&-1 & 1
 \end{bmatrix}.
\end{equation}
Hence the second regularization matrix is of size $(N-1)|n| \times N|n|$
and is given by $\vR_2 = \vD^{(N)}  \otimes I_{|n|}$.

Putting it all together, we can explicitly write 
\eqref{eq:optim} as
\begin{equation}
 \min_\vu \frac{1}{2} \|\vA \vu -\vv\|_2^2 + \alpha_1 \|\vI_N \otimes \vD_1  \vu\|_1 +
 \alpha_1 \|\vI_N \otimes \vD_2 \vu\|_1 + \alpha_2 \|\vD^{(N)} \otimes
 \vI_{|n|} \vu \|_1 + \alpha_3 \|\vu\|_1.
\end{equation}
Note that the regularization used in \cite{Majumder:2025:HDS} also
included spatial total variation and $\ell_1$ terms, but does not include
the total variation on the wavelength.

\subsection{The normal equations}
\label{sec:normal}
We now focus attention on the linear least squares problem
\eqref{eq:admm:a}, which we need to solve for each ADMM step with the $M=3$
regularization matrices defined above. To simplify the exposition we
rewrite it as
\begin{equation}
 \min_\vu \frac{1}{2} \| \vB \vu - \vb \|_2^2,
 \label{eq:lsq}
\end{equation}
where the system
matrix and right hand side are given by
\begin{equation}
\vB = \begin{bmatrix} 
\vA\\
\rho^{1/2} \vR_1\\
\rho^{1/2} \vR_2\\
\rho^{1/2} \vR_3
\end{bmatrix}
~\text{and}~
\vb = \begin{bmatrix}
\vv\\ 
\rho^{1/2}(\vz^{(k)}_1 - \rho^{-1} \vy^{(k)}_1)\\
\rho^{1/2}(\vz^{(k)}_2 - \rho^{-1} \vy^{(k)}_2)\\
\rho^{1/2}(\vz^{(k)}_3 - \rho^{-1} \vy^{(k)}_3)
\end{bmatrix},
\label{eq:Bandb}
\end{equation}
and $k$ is the current ADMM iteration number.  Clearly $\vB$ has more
rows than columns (overdetermined problem). As we shall see, $\vB$ is
guaranteed to be full rank provided $\rho>0$, so we can solve
\eqref{eq:lsq}  with the normal equations
\begin{equation}
 \vu = (\vB^* \vB)^{-1} \vB^* \vb.
 \label{eq:normal}
\end{equation}
A quick calculation reveals that
\begin{equation}
 \vB^* \vB = \vA^* \vA
 + \rho \left( \vI_N \otimes (\vD_1^* \vD_1 + \vD_2^* \vD_2)
 + \vL_N \otimes \vI_{|n|}
 + \vI_{N|n|}\right),
 \label{eq:bbs}
\end{equation}
where for $N\geq 2$, $\vL_{N}$ is the $N\times N$ discrete Laplacian with Neumann
boundary conditions, i.e.
\begin{equation}
 \vL_{N} = (\vD^{(N)})^* \vD^{(N)} =
 \begin{bmatrix}
  1 & -1\\
  -1 & 2 & -1 \\
   && \ddots\\
   &&&-1 & 2 & -1\\
   &&&& -1 & 1
 \end{bmatrix}.
 \label{eq:disclap}
\end{equation}
Notice from \eqref{eq:bbs} that $\vB^* \vB$ is a positive
semidefinite matrix plus a positive multiple of the identity and
therefore must be Hermitian positive definite and invertible. Hence
$\vB$ is full rank.

Moreover, the matrix $\vB^* \vB$ is a block matrix with $N \times N$
blocks each of size $|n| \times |n|$. Each block includes an operator of
the form $\vA_i^* \vA_j$. Since the $\vA_i$,  are $(p_1,p_2)-$periodic
convolutions and $(p_1,p_2)-$periodic convolutions are closed to products and adjoints,
$\vA_i^* \vA_j$ is also a $(p_1,p_2)-$periodic convolutions.  The blocks of
$\vB^* \vB$ also involve $\vD_1$, $\vD_2$ and $\vI_{n_1n_2}$ which are
all regular convolutions or $(1,1)-$periodic convolutions acting on signals of
dimension $n_1 \times n_2$.  As explained in \cite{Gonzalez:2026:FCC},
$(1,1)-$periodic convolutions are also $(p_1,p_2)-$periodic convolutions for any
integers $p_1,p_2>0$. We conclude that $\vB\vB^*$ can be evaluated using
$(p_1,p_2)-$periodic convolutions.  The following result from
\cite{Gonzalez:2026:FCC} comes in handy. 
\begin{lemma}
 \label{lem:blkdiag}
 Let $\vC$ be the matrix representation of a $(p_1,p_2)-$periodic
 circular convolution acting on signals of size $n_1 \times n_2$, where
 $n_j = p_j q_j$ for some integers $q_j$, $j=1,2$. Then  we have
 \begin{equation}
  \vC = \vQ \check \vC \vQ^*,
 \end{equation} 
 where $\check \vC$ is a block diagonal matrix with $|q|$ blocks of
 size $|p| \times |p|$ and $\vQ$ is a unitary matrix
 involving discrete Fourier transforms and depending only on $n_1,n_2$
 and $p_1,p_2$. Moreover applying $\vQ$ or its adjoint costs about the
 same as one DFT on data of dimensions $n_1 \times n_2$, i.e.
 $\cO(|n|\log|n|)$ flops.
\end{lemma}
By \cref{lem:blkdiag} we see that we only need to store $|p||n|$
scalars per block of $\vB^*\vB$, giving a total storage requirement of
$N^2 |p| |n|$ scalars.  If for instance $p_1=p_2=8$,
$n_1=n_2=256$ and $N=25$, storing $\vB^* \vB$ in double precision
arithmetic requires about 21GB, which could be prohibitive.

\subsection{The Sherman-Morrison-Woodbury formula}
\label{sec:smw}
To efficiently solve the normal equations, we observe that
$\vA^* \vA$ has at most rank $|n|$ and therefore $\vB^* \vB$ is a
$|n|$ rank update of a $|n| N \times |n| N$ invertible matrix,
i.e.  $\vB^*\vB = \vM + \vA^* \vA$, where 
\begin{equation}
 \vM = \rho \left(
   \vI_N \otimes (\vD_1^* \vD_1 + \vD_2^* \vD_2)
 + \vL_N \otimes \vI_{|n|}
 + \vI_{N|n|} \right),
\end{equation}
which is indeed Hermitian positive definite and thus invertible provided
$\rho>0$. Hence by the 
Sherman-Morrison-Woodbury (SMW) formula \cite{Woodbury:1950:IMM} we have
\begin{equation}
 (\vB^* \vB)^{-1} = \vM^{-1}  - \vM^{-1} \vA^* (\vI_{|n|} + \vA
 \vM^{-1} \vA^*)^{-1} \vA \vM^{-1}.
 \label{eq:smw}
\end{equation} 
Using the SMW formula to reduce the cost of solving the linear least
squares problems in ADMM is already mentioned in \cite{Parikh:2014:PA}
for Lasso.  As we shall see, all the matrices involved in \eqref{eq:smw}
are defined in terms of regular or periodic convolutions, which
significantly reduces the cost of solving a linear system with $\vB^*
\vB$. We explain in \cref{sec:M} how we solve systems with $\vM$. Then
in \cref{sec:inner}, we go over solving systems   with the inner matrix
in \eqref{eq:smw}, namely the $|n| \times |n|$ matrix $\vE = \vI_{|n|} +
\vA \vM^{-1} \vA^*$.

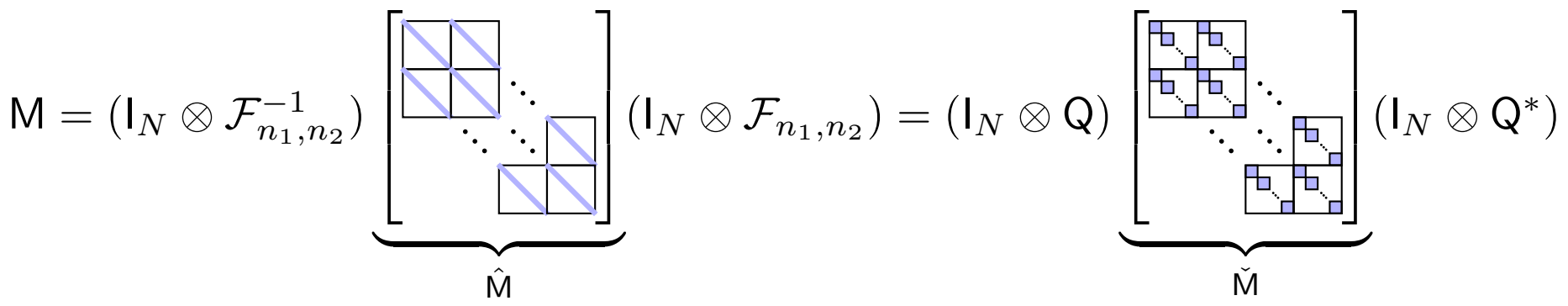
\begin{figure}
\newcommand{\Dblock}[1]{%
  \tikz[baseline=(current bounding box.center)]{
    \node[thin,draw,minimum width=#1, minimum height=#1,
          inner sep=0pt] (block) {};
    \draw[very thick,blue!30] (block.north west) -- (block.south east);
  }%
}%
\newcommand{\BDblock}[1]{%
\begin{tikzpicture}[x=#1/4,y=-#1/4,
      block/.style={draw,thin,fill=blue!30,minimum width=#1/4, minimum
      height=#1/4,inner sep=0pt},
      baseline=(current bounding box.center)
      ] 
\useasboundingbox (-0.5,-0.5) rectangle (3.5,3.5);
\node[block] at (0,0) {};
\node[block] at (1,1) {};
\fill[black] (1.8,1.8) circle (0.1);
\fill[black] (2,2)     circle (0.1);
\fill[black] (2.2,2.2) circle (0.1);
\node[block] at (3,3) {};
\draw[thin] (current bounding box.south west) rectangle (current bounding box.north east);
\end{tikzpicture}
}%
\newcommand{\DDblock}[1]{%
\begin{tikzpicture}[x=#1,y=-#1,baseline=(current bounding box.center)]
\useasboundingbox (-0.5,-0.5) rectangle (0.5,0.5);
\node at (-0.2,-0.2) {.};
\node at (0,0) {.};
\node at (0.2,0.2) {.};
\end{tikzpicture}
}%
\begin{center}
$
 \vM = 
 (\vI_N \otimes \cF_{n_1,n_2}^{-1})
 \underbrace{
 \left[\!\!
 \mathord{\vcenter{\hbox{
 \begin{tikzpicture}[x=2.5ex,y=-2.5ex]
 \useasboundingbox (-0.5,-0.8) rectangle (3.5,3.8);
 \node at (0,0) {\Dblock{2.5ex}};
 \node at (0,1) {\Dblock{2.5ex}};
 \node at (1,0) {\Dblock{2.5ex}};
 \node at (1,1) {\Dblock{2.5ex}};
 \node at (2,1) {\DDblock{2.5ex}};
 \node at (1,2) {\DDblock{2.5ex}};
 \node at (2,2) {\DDblock{2.5ex}};
 \node at (2,3) {\Dblock{2.5ex}};
 \node at (3,2) {\Dblock{2.5ex}};
 \node at (3,3) {\Dblock{2.5ex}};
 \end{tikzpicture}}}}
 \right]}_{\hat\vM}
(\vI_N \otimes \cF_{n_1,n_2})
=
(\vI_N \otimes \vQ)
\underbrace{
\left[\!\!
 \mathord{\vcenter{\hbox{
 \begin{tikzpicture}[x=2.5ex,y=-2.5ex]
 \useasboundingbox (-0.5,-0.8) rectangle (3.5,3.8);
 \node at (0,0) {\BDblock{2.5ex}};
 \node at (0,1) {\BDblock{2.5ex}};
 \node at (1,0) {\BDblock{2.5ex}};
 \node at (1,1) {\BDblock{2.5ex}};
 \node at (2,1) {\DDblock{2.5ex}};
 \node at (1,2) {\DDblock{2.5ex}};
 \node at (2,2) {\DDblock{2.5ex}};
 \node at (2,3) {\BDblock{2.5ex}};
 \node at (3,2) {\BDblock{2.5ex}};
 \node at (3,3) {\BDblock{2.5ex}};
 \end{tikzpicture}}}}
\right]}_{\check \vM}
(\vI_N \otimes \vQ^*)
$
\end{center}
\caption{A schematic of the two factorizations of $\vM$ that we use.
The matrix $\vM$ is block tridiagonal with
$N\times N$ blocks each of size $|n| \times |n|$. All 
blocks are regular convolutions on signals of size $n_1 \times n_2$, and
can therefore be diagonalized by $\cF_{n_1,n_2}$, as seen in the first
factorization and \cref{sec:M}. In the second factorization, the regular
convolutions are viewed as periodic circular convolutions and can be
unitarily transformed with $\vQ$ from \cref{lem:blkdiag} into block diagonal matrices. Each
block diagonal matrix has $|q|$ diagonal blocks each of size $|p| \times
|p|$, see \cref{sec:inner}.}
\label{fig:Mfacto}
\end{figure}

\subsubsection{Solving systems with $\vM$}
\label{sec:M}
By construction $\vM$ is a $N\times N$ block matrix where each block of
size $|n| \times |n|$ can be interpreted as a circular convolution on
tensors of size $n_1 \times n_2$.  We
recall that the (circular) convolution of two $2-$nd order tensors
$\vF,\vG$ of same size $n_1 \times n_2$ is given  by
\begin{equation}
 (\vF * \vG)[i_1,i_2] = \sum_{j_1=0}^{n_1-1} \sum_{j_2=0}^{n_2-1}
 \vF[i_1-j_1,i_2-j_2] \vG[j_1,j_2].
 \label{eq:conv}
\end{equation}
By the convolution theorem, convolution is multiplication in Fourier space, i.e.
\begin{equation}
 \cF_{n_1,n_2} (\vF * \vG) = (\cF_{n_1,n_2} \vF) \hprod (\cF_{n_1,n_2}
 \vG),
 \label{eq:conv:thm}
\end{equation}
where ``$\hprod$'' denotes the componentwise or Hadamard product and 
$\cF_{n_1,n_2}$ is the discrete Fourier transform of a $2-$nd order
tensor $\vF$ of size $n_1 \times n_2$, that is
\begin{equation}
 (\cF_{n_1,n_2}\vF)[i_1,i_2] = \sum_{j_1=0}^{n_1-1} \sum_{j_2=0}^{n_2-1}
 \omega_{n_1}^{i_1j_1} \omega_{n_2}^{i_2j_2} \vF[j_1,j_2],
 ~i_1=0,\ldots,n_1,~i_2=0,\ldots,n_2,
 \label{eq:dft}
\end{equation}
where $\omega_n = \exp[-\imath 2\pi /n ]$ and $\imath = \sqrt{-1}$, see
e.g. \cite{Kincaid:1996:NA}. Therefore we may factorize $\vM$ as 
\begin{equation}
 \vM = (\vI_N \otimes \cF_{n_1,n_2}^{-1}) \hat \vM (\vI_N \otimes
 \cF_{n_1,n_2}),
 \label{eq:M:facto}
\end{equation}
where $\hat \vM$ is a block tridiagonal matrix with $N\times N$ blocks,
each block being a diagonal matrix of size $|n| \times |n|$. This
factorization of $\vM$ is visualized in \cref{fig:Mfacto}.
In a slight abuse of notation, we also use $\cF_{n_1,n_2}$ to denote
the matrix representation for the linear operation in \cref{eq:dft}.
Upon appropriate permutations, $\hat \vM$ can be viewed as block
diagonal matrix with $|n|$ blocks, each block Hermitian tridiagonal of
size $N \times N$. Each tridiagonal system takes about $8N$ 
flops to solve (see \cite[Algorithm 4.3.6]{Golub:2013:MC}),
so the cost of solving a linear system with $\hat \vM$ is about $8N|n|$
flops. Asssuming applying $\cF_{n_1,n_2}$ or its inverse costs $|n|
\log |n|$ flops, it takes about $2N |n| \log |n| + 8N
|n|$ flops to solve a linear system with $\vM$.

\subsubsection{Solving systems with the SMW inner matrix}
\label{sec:inner}
Our goal here is to solve systems with the SMW inner matrix $\vE =
\vI_{|n|} + \vA \vM^{-1} \vA^*$. Since $\vA$ is composed of
$(p_1,p_2)-$periodic convolutions, we may use \cref{lem:blkdiag} to
factorize $\vA$ as
\begin{equation}
  \vA = \vQ \check \vA  (\vI_N \otimes \vQ^*),
  ~\text{where}~
  \check \vA = [ \check \vA_1,\ldots,\check \vA_N ],
  \label{eq:A:facto}
\end{equation}
and each $\check \vA_j$ is a block diagonal matrix with $|q|$
blocks of size $|p| \times |p|$.  It is also convenient to
interpret the circular convolutions in $\vM$ as $(p_1,p_2)-$periodic
convolutions, to factorize $\vM$ as follows
\begin{equation}
 \vM =  (\vI_N \otimes \vQ) \check \vM
  (\vI_N \otimes \vQ^*),
\end{equation}
where $\check \vM$ is a $N\times N$ block tridiagonal Hermitian matrix,
and each block is a block diagonal matrix with $|q|$ blocks of size
$|p| \times |p|$. This factorization of $\vM$ is depicted in
\cref{fig:Mfacto}. The factor $\check \vL$ in the Cholesky
factorization of $\check \vM = \check \vL \check \vL^*$ is also a $N
\times N$ block matrix which is lower triangular and with only one block diagonal
below the main block diagonal. Moreover, the blocks in $\vL$ have the same
structure as the blocks of $\check \vM$. Putting it all together we may
write 
\begin{equation}
 \vE = \vQ \check \vE \vQ^*,
 ~\text{where}~ 
 \check \vE = \vI_{|n|} + \check \vA \check \vL^{-*} \check \vL^{-1}
 \check \vA^*.
 \label{eq:E:facto}
\end{equation}
Since all the matrices
involved in the definition \eqref{eq:E:facto} of $\check\vE$ are
composed of block diagonal matrices with $|q|$ blocks
of size $|p| \times |p|$, we conclude that $\check\vE$ has the
same structure and thus $\vE$ is a $(p_1,p_2)-$periodic circular
convolution.

Let us now estimate the cost of building $\check \vE$. The dominant cost in
computing the Cholesky factorization of $\check \vM$ is computing the
Cholesky factorization of the $N$ diagonal blocks, with each costing about
$|q||p|^3/3 = |p|^2|n|/3$ flops \cite[Algorithm 4.2.1]{Golub:2013:MC}. Therefore computing the Cholesky
factor $\check \vL$ costs about $N |p|^2 |n|/3$ flops. The forward and
backward substitution with $\check \vL$ take about $N |p|^2 |n|$
flops when applied to $\check \vA^*$, see \cite[\S3.1]{Golub:2013:MC}. Finally left multiplication by
$\check \vA$ takes about the same number of flops. Therefore forming
$\check \vE$ costs about $10N |p|^2 |n|/3$ flops.

Once $\check \vE$ has been computed, its Cholesky factorization can be
found in $|p|^2|n|/3$ flops. With the Cholesky factorization of $\check \vE$ available, 
a linear system solve with $\vE$ takes only about $2|n|\log|n|
+ 2 |p| |n|$ flops.

\subsection{Cost estimate}
\label{sec:estimate}
To summarize, solving one linear system with $\vB^* \vB$ using the
Sherman-Morrison-Woodbury formula \cref{eq:smw} involves the following:
\begin{itemize}
 \item[(a)] three linear solves with $\vM$ at a cost of about $2N |n| \log |n| +
 8 N |n|$ flops each,
 \item[(b)] an application of $\vA$ and $\vA^*$, each costing
 about $2N|n|\log|n| + N|p||n|$ flops,
 \item[(c)] a linear solve with $\vE$ which costs about $2|n|\log|n|
+ 2|p| |n|$ flops.
\end{itemize}
Therefore, an approximate total of $(10N + 2)|n| \log |n| + N |n| (24 + 2|p|) +
2|p| |n|$ flops are needed to apply $(\vB^* \vB)^{-1}$.  An application of $\vB^*$ involves the following:
\begin{itemize}
 \item[(a)] one application of $\vA^*$, costing about $2N |n| \log
 |n| + N |p| |n|$ flops,
 \item[(b)] one application of $\vR_1^*$, costing about $2 N|n|$ flops,
 \item[(c)] one application of $\vR_2^*$, costing about $(N-1)|n|$
 flops,
 \item[(d)] one application of $\vR_3^*=\vI_{N|n|}$, costing no flops.
\end{itemize}
Hence, applying $\vB^*$ takes about $2N|n| \log|n| + N|p||n| + (6N-1)|n|$ flops.

Finally, a matrix vector product with $\vB^\dagger = (\vB^*\vB)^{-1} \vB^*$ takes about
$(12N+2)|n|\log|n| + N|n|(30N+3|p|) + (2|p|-1)|n|$ flops. Hence for $|p|$ fixed and large
$|n|$, $N$, applying $\vB^\dagger$ takes approximately $12N|n|\log |n|$ flops.

For comparison, the cost involved in a LSQR iteration to solve a least
squares problem with $\vA$ is roughly the cost of applying $\vA$ and
$\vA^*$, each once per iteration, that is $4N |n| \log |n| + 2 N |n|$
flops.  Therefore the SMW formula enables to compute the solution to the
normal equations for the augmented systems with roughly the same cost as
three LSQR iterations, for $|p|$ fixed and large $|n|$, $N$.  This
estimate is pessimistic, as the number of DFTs can be reduced when
applying $\vB^\dagger$.  Our Julia \cite{Bezanson:2017:JFA}
implementation in \texttt{DiffractogramImaging.jl}
\cite{DiffractogramImaging} uses fewer DFTs, with savings of about $N$
DFTs on signals of size $|n|$.

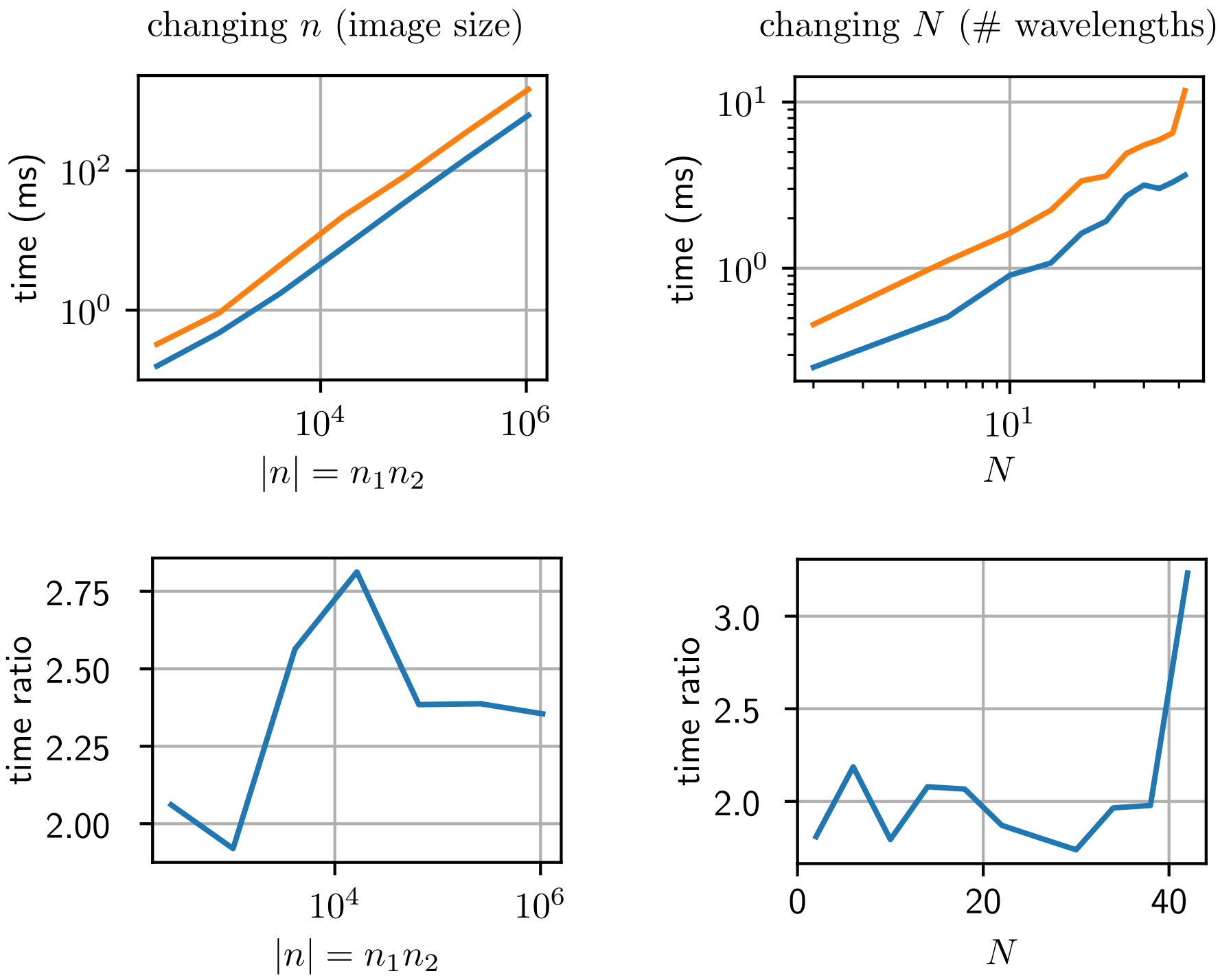
\begin{figure}
 \begin{center}
  \begin{tabular}{cc}
  \hspace{3em}changing $n$ (image size) & \hspace{3em}changing $N$ (\# wavelengths)\\
  \input{bench_vary_n_times.pgf} & \input{bench_vary_nw_times.pgf} \\
  \input{bench_vary_n_ratio.pgf} & \input{bench_vary_nw_ratio.pgf}
  \end{tabular}
 \end{center}
 \caption{A comparison between solving the augmented least squares
 problem with the direct method in \cref{sec:direct} (a matrix vector
 product with $\vB^\dagger$) and one step of the
 LSQR iteration to solve a least square problem with $\vA$. The left
 column shows a numerical experiment where we vary the image size $n$ while keeping $p$
 and $N$ fixed. The right column shows the number of wavelengths $N$ changing
 while keeping $n$ and $p$ fixed. The top row reports wall-clock times, while
 the bottom row has the ratio of times. As expected, a matrix vector product
 with $\vB^\dagger$ is about three times more expensive than a LSQR iteration.}
 \label{fig:bench}
\end{figure}

We checked this rough estimate against numerical experiments comparing
the computational cost of one LSQR step to solve a least squares problem
with $\vA$ (approximated by a matrix vector product by $\vA$ and
$\vA^*$) to one augmented least squares system solve (e.g. multiplying
by the pseudoinverse $\vB^\dagger$). The numerical experiments were
performed on a MacBook Pro M2 with an Apple M2 Max processor, 32GB of
RAM, running Mac OS X 26.6.2, Julia 1.13.2 and FFTW 3.3.12. The results
are summarized in \cref{fig:bench}. 

The first experiment (first column in \cref{fig:bench}) consists of
fixing $N=5$, $p=(8,8)$ and varying $n=(n_1,n_1)$ with $n_1 \in
\{2^4,2^5,\ldots,2^{10}\}$. The kernels used to build the matrix $\vA$
were randomly generated and the matrix vector products were against
randomly generated vectors of appropriate size. The top figure
corresponds to wall-clock times measured using
\texttt{BenchmarkTools.jl} with 20 samples, and does not include the
significant time it took to construct $\vA$ and $\vB^\dagger$. The blue
curve is the cost of one LSQR iteration and the orange curve is the cost
of one matrix vector multiplication with $\vB^\dagger$. As expected,
both operations become more expensive as we increase $|n|$, with
applying $\vB^\dagger$ being more expensive than one LSQR iteration. The
bottom figure shows the ratio of times between the two methods, which is
between 2 and 2.75, consistent with the fewer DFTs that are used in our
implementation.  The second experiment (second column in
\ref{fig:bench}) consists of fixing $n=(2^5,2^5)$, $p=(8,8)$ and picking
$N \in \{2,6,\ldots,42\}$.  Similarly, the computational cost increases
with $N$ and the time ratio is between 2 and 3.2, consistent with the
estimated time ration.

\section{RGB visualization of a hyperspectral image}
\label{sec:rgb}
We follow \cite{Foster:2019:HIC} to approximate what a hyperspectral
scene $u(x,\lambda)$ looks to the human eye. The main idea is to reduce
the hyperspectral data to red, green and blue channels by integrating
times color matching functions $R(\lambda)$, $G(\lambda)$ and
$B(\lambda)$, i.e.
\begin{equation}
 u_C(x) = \int C(\lambda) u(x,\lambda) \, d\lambda,~\text{for}~C \in \{R,G,B\}.
\end{equation}
The $R,G,B$ functions we use are a linear combination of the CIE
(Commission Internationale de l'Éclairage) XYZ tristimulus color
matching functions in \cite{xyz}. The XYZ values are a standard model
for the colors that the human eye can see. The specific $3 \times 3$
matrix we use to transform from XYZ values to RGB, assumes a sRGB color
space and a D50 illuminant and is taken from \cite{lindbloom}. The
resulting channels are rescaled to have maximal value one across
channels, projected to non-negative values and treated with a gamma
correction to enhance luminance.  To be more precise, the channels that
we display are $\max(u_C(x)/M,0))^{0.4}$ for $C \in \{R,G,B\}$ and $M =
\max \{ u_C(x) ~|~  x \in \real^2, C \in \{R,G,B\}\}$. 

\section{A synthetic example with an acquired PSF}
\label{sec:synth}
\begin{figure}
 \begin{center}
 \setlength{\imagesize}{0.2\textwidth}
 \begin{tabular}{c@{\hspace{0.4em}}c@{\hspace{0.4em}}c}
  true RGB & diffractogram & rec. RGB\\
  \includegraphics[width=\imagesize]{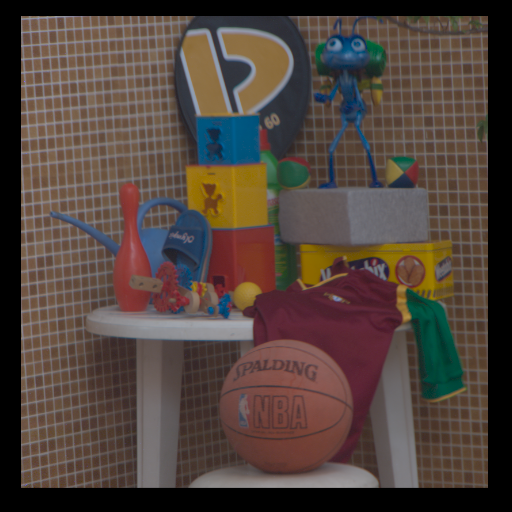} &
  \includegraphics[width=\imagesize]{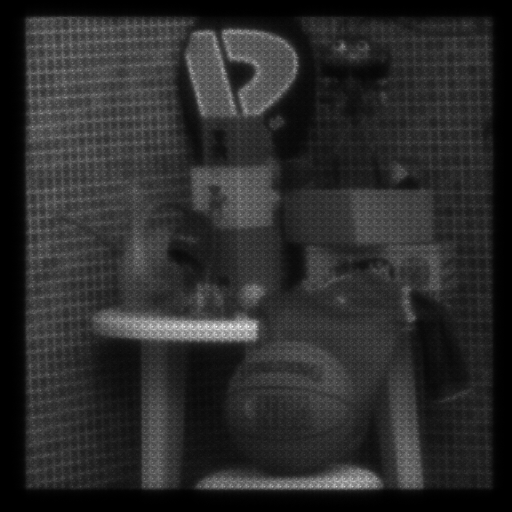} &
  \begin{tikzpicture}[x=\dimexpr\imagesize/512\relax,
                    y=\dimexpr\imagesize/512\relax,
		    every node/.style={font=\tiny}
		    ]
  \node[anchor=north west, inner sep=0] at (0,0)
    {\includegraphics[width=\imagesize]{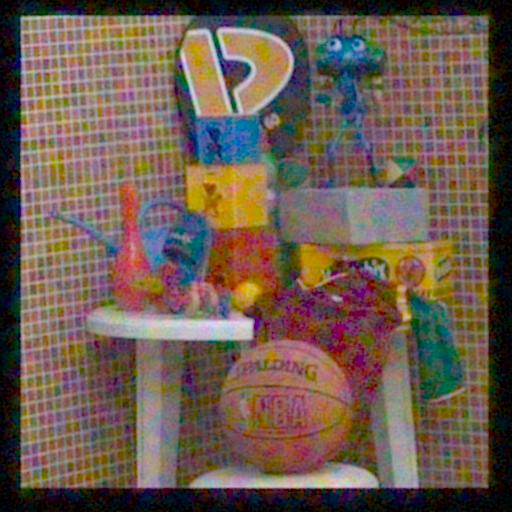}};
  \node at (242,-132) {1};
  \node at (242,-195) {2};
  \node at (245,-256) {3};
  \node at (341,-332) {4};
  \node at (380,-218) {5};
  \end{tikzpicture} \\
  \includegraphics[width=\imagesize]{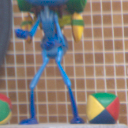} &
  \includegraphics[width=\imagesize]{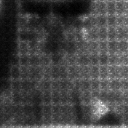} &
  \includegraphics[width=\imagesize]{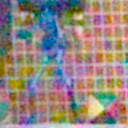}\\
  \includegraphics[width=\imagesize]{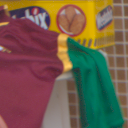} &
  \includegraphics[width=\imagesize]{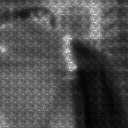} &
  \includegraphics[width=\imagesize]{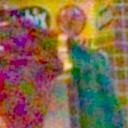}\\
 \end{tabular}
 \end{center}
 \caption{Synthetic reconstruction using a hyperspectral image from
 \cite{Foster:2006:FMN} and an acquired PSF from
 \cite{Majumder:2025:HDS}. The first column shows a RGB representation
 of the hyperspectral image and zoomed in areas. The second column
 shows the simulated diffractogram and zoomed in areas, with shot noise
 added. The third column is a RGB representation of the reconstruction.
 A comparison of the true and reconstructed spectra in the vicinity
 of locations labeled $1$--$5$ is included in \cref{fig:synth:conv}.
 The images in the first row are all $512 \times 512$ pixels while the
 zoomed in versions in the second and third rows are $128 \times 128$
 pixels. }
 \label{fig:synth}
\end{figure}
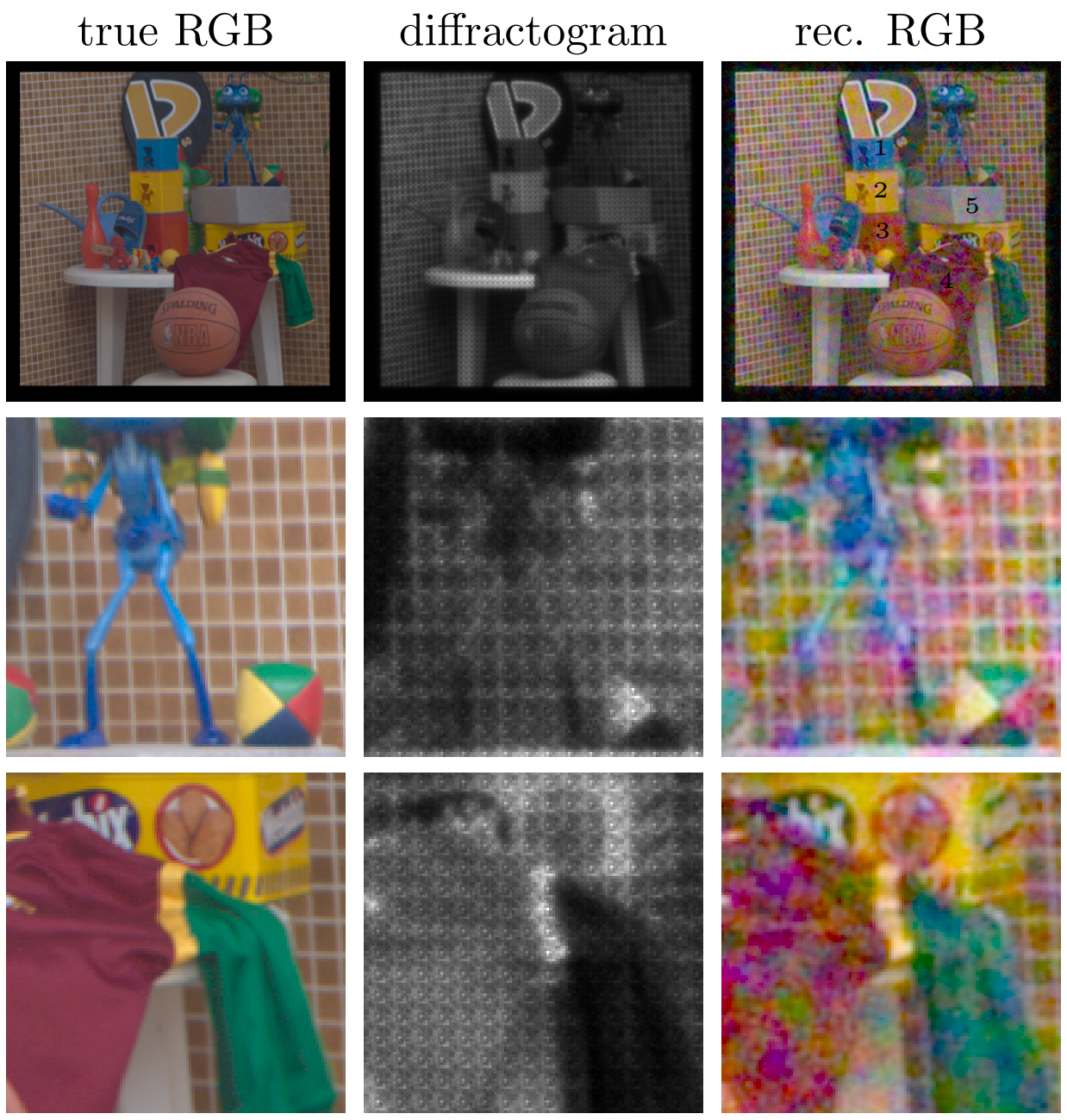

We start by testing our reconstruction method on a diffractogram
generated from a reference spectral image from \cite{Foster:2006:FMN}
and an experimentally acquired PSF from \cite{Majumder:2025:HDS}.  The
reference hyperspectral image from \cite{Foster:2006:FMN} is represented
by a $820\times 820 \times 31$ tensor, for wavelengths $\lambda \in
\{410,420,\ldots,710\}$nm. The data set was linearly interpolated to
reduce the number of pixels to $472 \times 472$ and to the $N=10$
wavelengths $\lambda \in \{440,470,\ldots,710\}$nm. As discussed in
\cref{sec:padding}, the image needs to be zero padded, which is done
here by adding $40$ pixels to each dimension, so that each $\vu_i$ has
dimensions $512 \times 512$. A RGB visualization of the reference or
``true'' hyperspectral data is included in the first column of
\cref{fig:synth}. 

The PSF we consider is $(8,8)-$periodic and is represented by a $48
\times 48 \times 8 \times 8 \times 26$ tensor. Each $48 \times 48$ image
corresponds to an experiment where a point-like object is placed at one
of the $8 \times 8$ positions within a period, for one of the particular
wavelengths $\lambda \in \{425,440,\ldots,800\}$nm. The PSF is
subsampled in wavelength to the $N=10$ wavelengths
$\{440,470,\ldots,710\}$nm, to match those of the reference
hyperspectral image.

The diffractogram (\cref{fig:synth} second column) is obtained by
calculating \eqref{eq:discrete} with $w_i=1$, $i=1,\ldots,N$ and then
adding shot noise corresponding to a $12-$bit sensor, like the sensor
used to acquire the PSF. To be more precise, the
values of $\vv = \vA \vu$ are rescaled to have maximum equal to
$2^{12} = 4096$, and the value of each pixel in the noisy image is 
drawn from a Poisson distribution with expectation equal to the pixel's
reference value. We recall that a Poisson distributed
random variable $X$ is characterized by its mean $\lambda$. The
probability that $X$ has the integer value $k$ is
\begin{equation}
 P(X=k) = \frac{1}{k!} \lambda^k e^{-\lambda},
\end{equation}
with expectation $\lambda$ equal to the reference value of the pixel.
The result is finally rescaled to have maximum one. Hence, the standard
deviation of the noise we add is roughly the square root of the signal
strength. Shot noise models noise in optical sensors such as the CMOS
sensor used to acquire the PSF, but ignores other noise sources
such as dark-current noise and noise added by the electronics, see
e.g.  \cite{EMVA1288}.

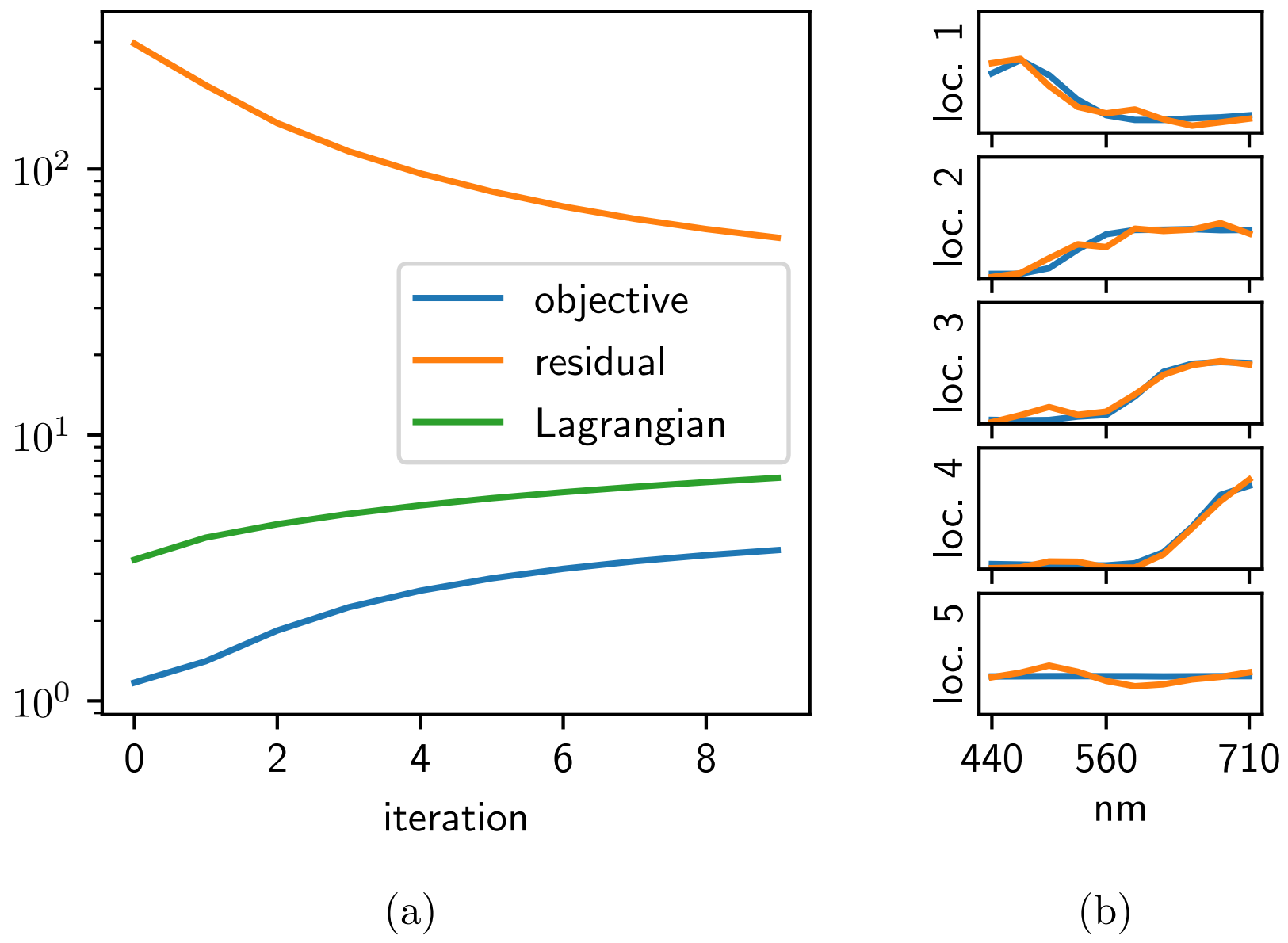
\begin{figure}
 \begin{tabular}{cc}
 \input{synth_convergence.pgf} &
 \input{synth_spec.pgf}
 \\
 (a) & (b)
 \end{tabular}
 \caption{(a) Synthetic reconstruction convergence history. (b) Spectra at
 marked locations in \cref{fig:synth}. Blue line is the true spectrum and orange
 the reconstructed one. The $y$ axes in (b) range from $0$ to $1$.}
 \label{fig:synth:conv}
\end{figure}

A RGB visualization of the reconstructions using the reconstruction
algorithm with $\rho=0.005$, $\alpha_1=0.5$, $\alpha_2=10^{-5}$, $\alpha_3
= 10^{-10}$, and stopping after 10 iterations. The reconstructions take
approximatively $17$ seconds on a MacBook Pro with an Apple M2 Max
processor. We report
in \cref{fig:synth:conv}(a) the convergence history of the objective
function, the residual $\sum_{j=1}^3 \| R_j \vu - \vz_j \|^2$ and the
augmented Lagrangian $L_\rho$ in \cref{eq:auglag}. We see that the
residual diminishes while the Lagrangian and objective function slightly
increase. Note that we cannot expect monotonic convergence from the ADMM
method.

In the zoomed-in versions in \cref{fig:synth}, we see that different
patterns appear in the diffractogram, depending on the spectrum of the
area. These fluctuations about a smooth background are the features that
contain information about the spectrum.  The reconstructions are far
from perfect, as could be expected from a severely ill-posed problem.
Nevertheless, we display in \cref{fig:synth:conv}(b) the reconstructed
(orange) vs true (blue) spectra in the vicinity of the five sample
locations indicated in the top right image in \cref{fig:synth}. The
spectra we report are averaged on a $8\times 8$ pixel area, and the
curves are rescaled by their $2-$norm. The average relative error in the
$2-$norm across all five locations is of about 14.5\%.

\section{Reconstructions from acquired diffractograms}
\label{sec:realdata}
We tested our approach with diffractograms acquired experimentally in
\cite{Majumder:2025:HDS}. Here the PSF is a
$48 \times 48 \times 8 \times 8 \times 25$ experimentally acquired 
$5-$th order tensor for wavelengths $\lambda \in \{ 440, 455,
\ldots , 800\}$ nm.

In \cref{fig:macbeth}, we reconstructed diffractograms obtained by
illuminating 18 different squares in the Macbeth/X-Rite ColorChecker
chart (a color calibration target used in photography), with a
collimated beam of white light restricted by a circular aperture. The
reconstructions were performed with regularization parameters
$\rho=0.03$, $\alpha_1=1$, $\alpha_2=10^{-8}$, $\alpha_3=10^{-9}$ and 10
ADMM iterations. For each of the diffractograms we report a RGB
visualization of the reconstruction and we sample the reconstructed
spectrum near the center of the image, averaged over $8\times 8$ pixels.
The reconstructed spectrum is compared to a spectrum measured
experimentally with a spectrometer, also from \cite{Majumder:2025:HDS}.
The comparison is done by renormalizing both spectra so that they have
the same $2-$norm and we additionally report the relative $2-$norm error
as a percentage. As can be seen from \cref{fig:macbeth}, the spectra we
obtain are qualitatively similar to the reference ones, but the average
$2-$norm error is about $32\%$, higher than the one in \cref{sec:synth},
and possibly due to the following factors.
First, the DFA fabrication and positioning over the sensor are not
perfect and thus PSF periodicity is only an approximation to the real optical
system. Second, the reconstructions clearly
overestimate the 440nm band. During the experiments to acquire the PSF,
the intensity of the source for this particular band was significantly
lower than for the other bands, leading to more noise.

\begin{figure}
\begin{center}
 \begin{tabular}{cc}
  \includegraphics[width=0.3\textwidth]{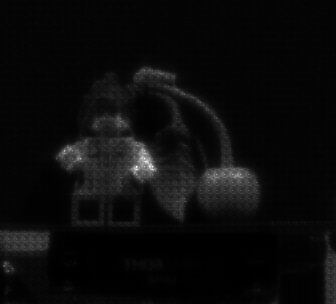} &
  \includegraphics[width=0.3\textwidth]{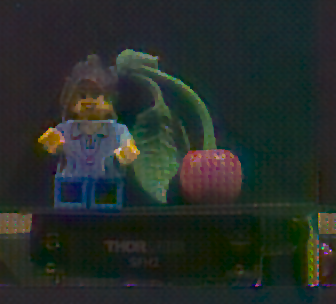}\\
  (a) & (b)
 \end{tabular}
\end{center}
\caption{(a) An experimentally acquired diffractogram of plastic toys on a
black background. (b) A RGB visualization of the reconstruction. Both
images are $304 \times 336$ pixels.}
\label{fig:lego}
\end{figure}

\Cref{fig:lego}(a) shows an acquired diffractogram of a scene
composed of toys on a black background. The RGB of the reconstruction is
shown in \cref{fig:lego}(b). The regularization parameters are
$\rho=0.01$, $\alpha_1=2$, $\alpha_2=10^{-3}$, $\alpha_3=10^{-6}$ and 20
ADMM iterations. No reference spectra are available for this experiment.
This image is included to illustrate two artifacts of our method. The
first one is due to using the PCC without appropriately zero padding the
image. Indeed, the PCC operation assumes the diffractogram is periodic
and this appears as bands on the right and bottom of the image, which
correspond to areas that are to the left and top of the of the actual
scene. The second is that because of parallax and other experimental errors,
the periodicity assumption of the kernel is not perfectly satisfied and
thus reconstructions have artifacts for pixels far from the calibration
region. In this particular example, there is a red (resp. green) tint in
the left (resp. right) side of the background.

\newcommand{\mb}[2]{%
\shortstack{%
\vspace{0.3em}\\
\includegraphics[trim=0 24 24 0,clip=true,width=0.08\textwidth]{mb_#1_rnrm}\\
\includegraphics[width=0.11\textwidth]{mb_#1_rnrm_spec}\\
\tiny{#2}}}

\begin{figure}
\begin{center}
 \begin{tabular}{|c|c|c|c|c|c|}
 \hline
 \mb{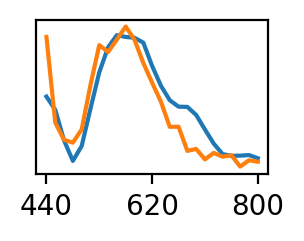} {31.8\%} & 
 \mb{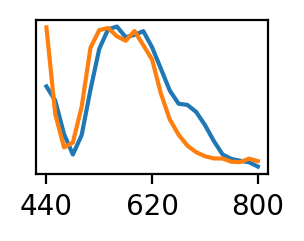} {30.8\%} & 
 \mb{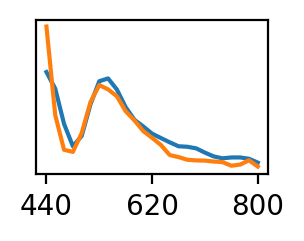} {30.7\%} & 
 \mb{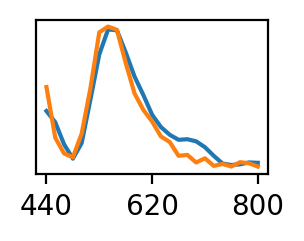} {20.2\%} & 
 \mb{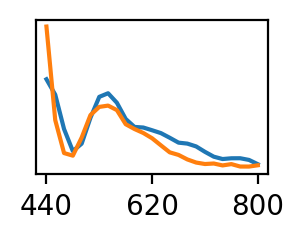} {38.1\%} & 
 \mb{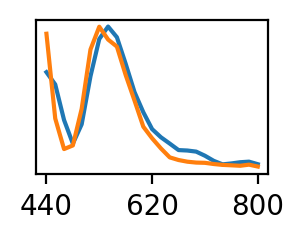} {24.5\%} \\\hline
 \mb{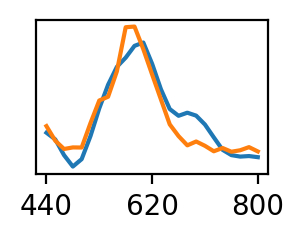} {28.8\%} & 
 \mb{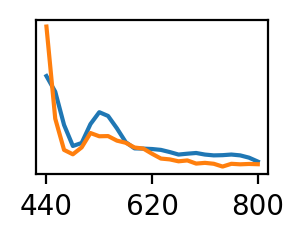} {47.4\%} & 
 \mb{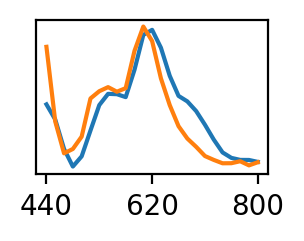} {35.5\%} & 
 \mb{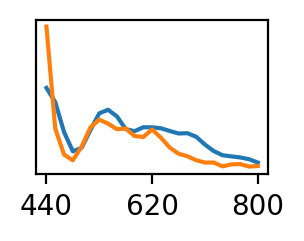}{48.3\%} & 
 \mb{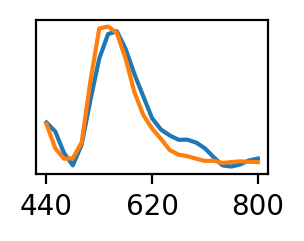}{21.8\%} & 
 \mb{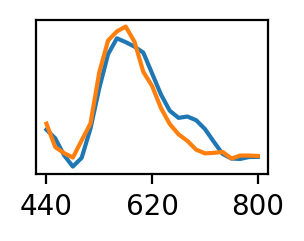}{23.5\%} \\\hline
 \mb{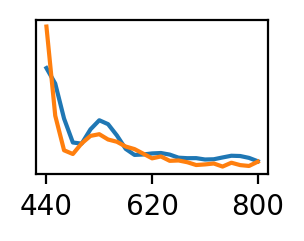}{43.3\%} & 
 \mb{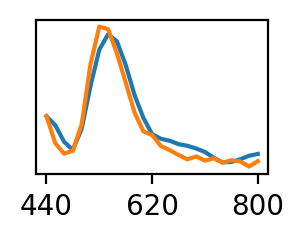}{21.8\%} & 
 \mb{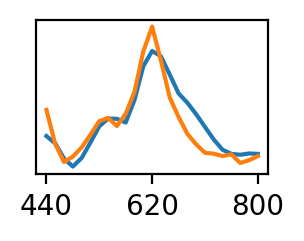}{31.1\%} & 
 \mb{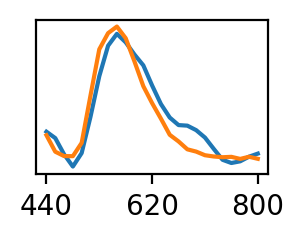}{23.9\%} & 
 \mb{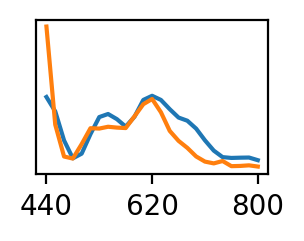}{45.2\%} & 
 \mb{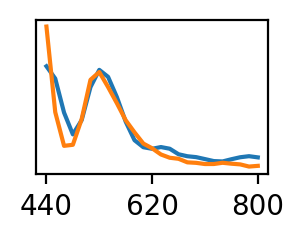}{29.2\%} \\\hline
 \end{tabular}
\end{center}
\caption{Reconstructions from experimental data taken in
\cite{Majumder:2025:HDS} of squares in a Macbeth/X-Rite ColorChecker
chart.  We compare the spectra at the center of each image with spectra
measured with a spectrometer from \cite{Majumder:2025:HDS}.  The two
signals are rescaled so that they have the same $2-$norm. The horizontal
axis is the wavelength in nanometers and the vertical axis is the
irradiance in arbitrary units.  The percentage is the $2-$norm relative
error between the spectra. The mean $2-$norm relative error is about
32\%.}
\label{fig:macbeth}
\end{figure}

\section{Summary and perspectives} 
\label{sec:summary} 
We presented a method to reconstruct the spatial and spectral
distribution of a scene from a single monochrome image acquired through
an optical system with a periodic wavelength-dependent PSF. While
recovering $N$ spectral channels from a single diffractogram is severely
underdetermined, the problem can be regularized using spatial and
spectral total variation and solved efficiently by exploiting the
periodicity of the PSF. Specifically, using periodic circular
convolutions allows the augmented linear systems within ADMM to be
solved with an exact direct method that costs roughly as much as three
LSQR iterations. Reconstructions from synthetic and experimental data
confirm that our method can reliably recover spatial and spectral
information in the visible spectrum.

The method remains sensitive to noise, both during PSF calibration
\cite{Majumder:2025:HDS} and in acquiring the diffractogram signal.
Because spectral information is encoded as spatial variations about a
smooth background, experimental noise can overwhelm these
variations. This may explain the higher reconstruction errors observed
in experimental data compared to synthetic benchmarks.

A key assumption of our approach is that the PSF is periodic.  As shown
in \cite{Majumder:2025:HDS} and illustrated in \cref{fig:lego}, the
physical PSF slowly drifts across the field of view due to e.g. parallax, DFA
fabrication tolerances or DFA positioning on the sensor. Extending the
direct method to account for such slow, long-range spatial
variations (with e.g. coordinate mappings or piecewise-periodic
approximations), is a direction of future study.

Finally, the periodic DFA used here was designed using a heuristic
parameter search \cite{Majumder:2025:HDS}. A deeper understanding of
what makes a DFA geometry good for spectral discrimination would be
welcome. Studies of wave propagation through random rough interfaces
such as \cite{Gomez:2026:IRI}, could help answer this question.

\bibliographystyle{siamplain}
\bibliography{hid}

\end{document}

%% file: bench_vary_n_times.pgf
\begingroup%
\makeatletter%
\begin{pgfpicture}%
\pgfpathrectangle{\pgfpointorigin}{\pgfqpoint{2.260496in}{1.781079in}}%
\pgfusepath{use as bounding box, clip}%
\begin{pgfscope}%
\pgfsetbuttcap%
\pgfsetmiterjoin%
\definecolor{currentfill}{rgb}{1.000000,1.000000,1.000000}%
\pgfsetfillcolor{currentfill}%
\pgfsetlinewidth{0.000000pt}%
\definecolor{currentstroke}{rgb}{1.000000,1.000000,1.000000}%
\pgfsetstrokecolor{currentstroke}%
\pgfsetdash{}{0pt}%
\pgfpathmoveto{\pgfqpoint{0.000000in}{0.000000in}}%
\pgfpathlineto{\pgfqpoint{2.260496in}{0.000000in}}%
\pgfpathlineto{\pgfqpoint{2.260496in}{1.781079in}}%
\pgfpathlineto{\pgfqpoint{0.000000in}{1.781079in}}%
\pgfpathlineto{\pgfqpoint{0.000000in}{0.000000in}}%
\pgfpathclose%
\pgfusepath{fill}%
\end{pgfscope}%
\begin{pgfscope}%
\pgfsetbuttcap%
\pgfsetmiterjoin%
\definecolor{currentfill}{rgb}{1.000000,1.000000,1.000000}%
\pgfsetfillcolor{currentfill}%
\pgfsetlinewidth{0.000000pt}%
\definecolor{currentstroke}{rgb}{0.000000,0.000000,0.000000}%
\pgfsetstrokecolor{currentstroke}%
\pgfsetstrokeopacity{0.000000}%
\pgfsetdash{}{0pt}%
\pgfpathmoveto{\pgfqpoint{0.588387in}{0.526079in}}%
\pgfpathlineto{\pgfqpoint{2.138387in}{0.526079in}}%
\pgfpathlineto{\pgfqpoint{2.138387in}{1.681079in}}%
\pgfpathlineto{\pgfqpoint{0.588387in}{1.681079in}}%
\pgfpathlineto{\pgfqpoint{0.588387in}{0.526079in}}%
\pgfpathclose%
\pgfusepath{fill}%
\end{pgfscope}%
\begin{pgfscope}%
\pgfpathrectangle{\pgfqpoint{0.588387in}{0.526079in}}{\pgfqpoint{1.550000in}{1.155000in}}%
\pgfusepath{clip}%
\pgfsetrectcap%
\pgfsetroundjoin%
\pgfsetlinewidth{0.803000pt}%
\definecolor{currentstroke}{rgb}{0.690196,0.690196,0.690196}%
\pgfsetstrokecolor{currentstroke}%
\pgfsetdash{}{0pt}%
\pgfpathmoveto{\pgfqpoint{1.279747in}{0.526079in}}%
\pgfpathlineto{\pgfqpoint{1.279747in}{1.681079in}}%
\pgfusepath{stroke}%
\end{pgfscope}%
\begin{pgfscope}%
\pgfsetbuttcap%
\pgfsetroundjoin%
\definecolor{currentfill}{rgb}{0.000000,0.000000,0.000000}%
\pgfsetfillcolor{currentfill}%
\pgfsetlinewidth{0.803000pt}%
\definecolor{currentstroke}{rgb}{0.000000,0.000000,0.000000}%
\pgfsetstrokecolor{currentstroke}%
\pgfsetdash{}{0pt}%
\pgfsys@defobject{currentmarker}{\pgfqpoint{0.000000in}{-0.048611in}}{\pgfqpoint{0.000000in}{0.000000in}}{%
\pgfpathmoveto{\pgfqpoint{0.000000in}{0.000000in}}%
\pgfpathlineto{\pgfqpoint{0.000000in}{-0.048611in}}%
\pgfusepath{stroke,fill}%
}%
\begin{pgfscope}%
\pgfsys@transformshift{1.279747in}{0.526079in}%
\pgfsys@useobject{currentmarker}{}%
\end{pgfscope}%
\end{pgfscope}%
\begin{pgfscope}%
\definecolor{textcolor}{rgb}{0.000000,0.000000,0.000000}%
\pgfsetstrokecolor{textcolor}%
\pgfsetfillcolor{textcolor}%
\pgftext[x=1.279747in,y=0.428857in,,top]{\color{textcolor}{\sffamily\fontsize{10.000000}{12.000000}\selectfont\catcode`\^=\active\def^{\ifmmode\sp\else\^{}\fi}\catcode`\%=\active\def
\end{pgfscope}%
\begin{pgfscope}%
\pgfpathrectangle{\pgfqpoint{0.588387in}{0.526079in}}{\pgfqpoint{1.550000in}{1.155000in}}%
\pgfusepath{clip}%
\pgfsetrectcap%
\pgfsetroundjoin%
\pgfsetlinewidth{0.803000pt}%
\definecolor{currentstroke}{rgb}{0.690196,0.690196,0.690196}%
\pgfsetstrokecolor{currentstroke}%
\pgfsetdash{}{0pt}%
\pgfpathmoveto{\pgfqpoint{2.059897in}{0.526079in}}%
\pgfpathlineto{\pgfqpoint{2.059897in}{1.681079in}}%
\pgfusepath{stroke}%
\end{pgfscope}%
\begin{pgfscope}%
\pgfsetbuttcap%
\pgfsetroundjoin%
\definecolor{currentfill}{rgb}{0.000000,0.000000,0.000000}%
\pgfsetfillcolor{currentfill}%
\pgfsetlinewidth{0.803000pt}%
\definecolor{currentstroke}{rgb}{0.000000,0.000000,0.000000}%
\pgfsetstrokecolor{currentstroke}%
\pgfsetdash{}{0pt}%
\pgfsys@defobject{currentmarker}{\pgfqpoint{0.000000in}{-0.048611in}}{\pgfqpoint{0.000000in}{0.000000in}}{%
\pgfpathmoveto{\pgfqpoint{0.000000in}{0.000000in}}%
\pgfpathlineto{\pgfqpoint{0.000000in}{-0.048611in}}%
\pgfusepath{stroke,fill}%
}%
\begin{pgfscope}%
\pgfsys@transformshift{2.059897in}{0.526079in}%
\pgfsys@useobject{currentmarker}{}%
\end{pgfscope}%
\end{pgfscope}%
\begin{pgfscope}%
\definecolor{textcolor}{rgb}{0.000000,0.000000,0.000000}%
\pgfsetstrokecolor{textcolor}%
\pgfsetfillcolor{textcolor}%
\pgftext[x=2.059897in,y=0.428857in,,top]{\color{textcolor}{\sffamily\fontsize{10.000000}{12.000000}\selectfont\catcode`\^=\active\def^{\ifmmode\sp\else\^{}\fi}\catcode`\%=\active\def
\end{pgfscope}%
\begin{pgfscope}%
\definecolor{textcolor}{rgb}{0.000000,0.000000,0.000000}%
\pgfsetstrokecolor{textcolor}%
\pgfsetfillcolor{textcolor}%
\pgftext[x=1.363387in,y=0.238889in,,top]{\color{textcolor}{\sffamily\fontsize{10.000000}{12.000000}\selectfont\catcode`\^=\active\def^{\ifmmode\sp\else\^{}\fi}\catcode`\%=\active\def
\end{pgfscope}%
\begin{pgfscope}%
\pgfpathrectangle{\pgfqpoint{0.588387in}{0.526079in}}{\pgfqpoint{1.550000in}{1.155000in}}%
\pgfusepath{clip}%
\pgfsetrectcap%
\pgfsetroundjoin%
\pgfsetlinewidth{0.803000pt}%
\definecolor{currentstroke}{rgb}{0.690196,0.690196,0.690196}%
\pgfsetstrokecolor{currentstroke}%
\pgfsetdash{}{0pt}%
\pgfpathmoveto{\pgfqpoint{0.588387in}{0.791034in}}%
\pgfpathlineto{\pgfqpoint{2.138387in}{0.791034in}}%
\pgfusepath{stroke}%
\end{pgfscope}%
\begin{pgfscope}%
\pgfsetbuttcap%
\pgfsetroundjoin%
\definecolor{currentfill}{rgb}{0.000000,0.000000,0.000000}%
\pgfsetfillcolor{currentfill}%
\pgfsetlinewidth{0.803000pt}%
\definecolor{currentstroke}{rgb}{0.000000,0.000000,0.000000}%
\pgfsetstrokecolor{currentstroke}%
\pgfsetdash{}{0pt}%
\pgfsys@defobject{currentmarker}{\pgfqpoint{-0.048611in}{0.000000in}}{\pgfqpoint{-0.000000in}{0.000000in}}{%
\pgfpathmoveto{\pgfqpoint{-0.000000in}{0.000000in}}%
\pgfpathlineto{\pgfqpoint{-0.048611in}{0.000000in}}%
\pgfusepath{stroke,fill}%
}%
\begin{pgfscope}%
\pgfsys@transformshift{0.588387in}{0.791034in}%
\pgfsys@useobject{currentmarker}{}%
\end{pgfscope}%
\end{pgfscope}%
\begin{pgfscope}%
\definecolor{textcolor}{rgb}{0.000000,0.000000,0.000000}%
\pgfsetstrokecolor{textcolor}%
\pgfsetfillcolor{textcolor}%
\pgftext[x=0.289968in, y=0.738272in, left, base]{\color{textcolor}{\sffamily\fontsize{10.000000}{12.000000}\selectfont\catcode`\^=\active\def^{\ifmmode\sp\else\^{}\fi}\catcode`\%=\active\def
\end{pgfscope}%
\begin{pgfscope}%
\pgfpathrectangle{\pgfqpoint{0.588387in}{0.526079in}}{\pgfqpoint{1.550000in}{1.155000in}}%
\pgfusepath{clip}%
\pgfsetrectcap%
\pgfsetroundjoin%
\pgfsetlinewidth{0.803000pt}%
\definecolor{currentstroke}{rgb}{0.690196,0.690196,0.690196}%
\pgfsetstrokecolor{currentstroke}%
\pgfsetdash{}{0pt}%
\pgfpathmoveto{\pgfqpoint{0.588387in}{1.320170in}}%
\pgfpathlineto{\pgfqpoint{2.138387in}{1.320170in}}%
\pgfusepath{stroke}%
\end{pgfscope}%
\begin{pgfscope}%
\pgfsetbuttcap%
\pgfsetroundjoin%
\definecolor{currentfill}{rgb}{0.000000,0.000000,0.000000}%
\pgfsetfillcolor{currentfill}%
\pgfsetlinewidth{0.803000pt}%
\definecolor{currentstroke}{rgb}{0.000000,0.000000,0.000000}%
\pgfsetstrokecolor{currentstroke}%
\pgfsetdash{}{0pt}%
\pgfsys@defobject{currentmarker}{\pgfqpoint{-0.048611in}{0.000000in}}{\pgfqpoint{-0.000000in}{0.000000in}}{%
\pgfpathmoveto{\pgfqpoint{-0.000000in}{0.000000in}}%
\pgfpathlineto{\pgfqpoint{-0.048611in}{0.000000in}}%
\pgfusepath{stroke,fill}%
}%
\begin{pgfscope}%
\pgfsys@transformshift{0.588387in}{1.320170in}%
\pgfsys@useobject{currentmarker}{}%
\end{pgfscope}%
\end{pgfscope}%
\begin{pgfscope}%
\definecolor{textcolor}{rgb}{0.000000,0.000000,0.000000}%
\pgfsetstrokecolor{textcolor}%
\pgfsetfillcolor{textcolor}%
\pgftext[x=0.289968in, y=1.267409in, left, base]{\color{textcolor}{\sffamily\fontsize{10.000000}{12.000000}\selectfont\catcode`\^=\active\def^{\ifmmode\sp\else\^{}\fi}\catcode`\%=\active\def
\end{pgfscope}%
\begin{pgfscope}%
\definecolor{textcolor}{rgb}{0.000000,0.000000,0.000000}%
\pgfsetstrokecolor{textcolor}%
\pgfsetfillcolor{textcolor}%
\pgftext[x=0.234413in,y=1.103579in,,bottom,rotate=90.000000]{\color{textcolor}{\sffamily\fontsize{10.000000}{12.000000}\selectfont\catcode`\^=\active\def^{\ifmmode\sp\else\^{}\fi}\catcode`\%=\active\def
\end{pgfscope}%
\begin{pgfscope}%
\pgfpathrectangle{\pgfqpoint{0.588387in}{0.526079in}}{\pgfqpoint{1.550000in}{1.155000in}}%
\pgfusepath{clip}%
\pgfsetrectcap%
\pgfsetroundjoin%
\pgfsetlinewidth{1.505625pt}%
\definecolor{currentstroke}{rgb}{0.121569,0.466667,0.705882}%
\pgfsetstrokecolor{currentstroke}%
\pgfsetdash{}{0pt}%
\pgfpathmoveto{\pgfqpoint{0.658842in}{0.578579in}}%
\pgfpathlineto{\pgfqpoint{0.893690in}{0.704593in}}%
\pgfpathlineto{\pgfqpoint{1.128539in}{0.855861in}}%
\pgfpathlineto{\pgfqpoint{1.363387in}{1.026899in}}%
\pgfpathlineto{\pgfqpoint{1.598236in}{1.199253in}}%
\pgfpathlineto{\pgfqpoint{1.833084in}{1.366801in}}%
\pgfpathlineto{\pgfqpoint{2.067933in}{1.530178in}}%
\pgfusepath{stroke}%
\end{pgfscope}%
\begin{pgfscope}%
\pgfpathrectangle{\pgfqpoint{0.588387in}{0.526079in}}{\pgfqpoint{1.550000in}{1.155000in}}%
\pgfusepath{clip}%
\pgfsetrectcap%
\pgfsetroundjoin%
\pgfsetlinewidth{1.505625pt}%
\definecolor{currentstroke}{rgb}{1.000000,0.498039,0.054902}%
\pgfsetstrokecolor{currentstroke}%
\pgfsetdash{}{0pt}%
\pgfpathmoveto{\pgfqpoint{0.658842in}{0.661676in}}%
\pgfpathlineto{\pgfqpoint{0.893690in}{0.779532in}}%
\pgfpathlineto{\pgfqpoint{1.128539in}{0.964034in}}%
\pgfpathlineto{\pgfqpoint{1.363387in}{1.145724in}}%
\pgfpathlineto{\pgfqpoint{1.598236in}{1.299097in}}%
\pgfpathlineto{\pgfqpoint{1.833084in}{1.466778in}}%
\pgfpathlineto{\pgfqpoint{2.067933in}{1.628579in}}%
\pgfusepath{stroke}%
\end{pgfscope}%
\begin{pgfscope}%
\pgfsetrectcap%
\pgfsetmiterjoin%
\pgfsetlinewidth{0.803000pt}%
\definecolor{currentstroke}{rgb}{0.000000,0.000000,0.000000}%
\pgfsetstrokecolor{currentstroke}%
\pgfsetdash{}{0pt}%
\pgfpathmoveto{\pgfqpoint{0.588387in}{0.526079in}}%
\pgfpathlineto{\pgfqpoint{0.588387in}{1.681079in}}%
\pgfusepath{stroke}%
\end{pgfscope}%
\begin{pgfscope}%
\pgfsetrectcap%
\pgfsetmiterjoin%
\pgfsetlinewidth{0.803000pt}%
\definecolor{currentstroke}{rgb}{0.000000,0.000000,0.000000}%
\pgfsetstrokecolor{currentstroke}%
\pgfsetdash{}{0pt}%
\pgfpathmoveto{\pgfqpoint{2.138387in}{0.526079in}}%
\pgfpathlineto{\pgfqpoint{2.138387in}{1.681079in}}%
\pgfusepath{stroke}%
\end{pgfscope}%
\begin{pgfscope}%
\pgfsetrectcap%
\pgfsetmiterjoin%
\pgfsetlinewidth{0.803000pt}%
\definecolor{currentstroke}{rgb}{0.000000,0.000000,0.000000}%
\pgfsetstrokecolor{currentstroke}%
\pgfsetdash{}{0pt}%
\pgfpathmoveto{\pgfqpoint{0.588387in}{0.526079in}}%
\pgfpathlineto{\pgfqpoint{2.138387in}{0.526079in}}%
\pgfusepath{stroke}%
\end{pgfscope}%
\begin{pgfscope}%
\pgfsetrectcap%
\pgfsetmiterjoin%
\pgfsetlinewidth{0.803000pt}%
\definecolor{currentstroke}{rgb}{0.000000,0.000000,0.000000}%
\pgfsetstrokecolor{currentstroke}%
\pgfsetdash{}{0pt}%
\pgfpathmoveto{\pgfqpoint{0.588387in}{1.681079in}}%
\pgfpathlineto{\pgfqpoint{2.138387in}{1.681079in}}%
\pgfusepath{stroke}%
\end{pgfscope}%
\end{pgfpicture}%
\makeatother%
\endgroup%

%% file: bench_vary_nw_times.pgf
\begingroup%
\makeatletter%
\begin{pgfpicture}%
\pgfpathrectangle{\pgfpointorigin}{\pgfqpoint{2.238387in}{1.776603in}}%
\pgfusepath{use as bounding box, clip}%
\begin{pgfscope}%
\pgfsetbuttcap%
\pgfsetmiterjoin%
\definecolor{currentfill}{rgb}{1.000000,1.000000,1.000000}%
\pgfsetfillcolor{currentfill}%
\pgfsetlinewidth{0.000000pt}%
\definecolor{currentstroke}{rgb}{1.000000,1.000000,1.000000}%
\pgfsetstrokecolor{currentstroke}%
\pgfsetdash{}{0pt}%
\pgfpathmoveto{\pgfqpoint{0.000000in}{0.000000in}}%
\pgfpathlineto{\pgfqpoint{2.238387in}{0.000000in}}%
\pgfpathlineto{\pgfqpoint{2.238387in}{1.776603in}}%
\pgfpathlineto{\pgfqpoint{0.000000in}{1.776603in}}%
\pgfpathlineto{\pgfqpoint{0.000000in}{0.000000in}}%
\pgfpathclose%
\pgfusepath{fill}%
\end{pgfscope}%
\begin{pgfscope}%
\pgfsetbuttcap%
\pgfsetmiterjoin%
\definecolor{currentfill}{rgb}{1.000000,1.000000,1.000000}%
\pgfsetfillcolor{currentfill}%
\pgfsetlinewidth{0.000000pt}%
\definecolor{currentstroke}{rgb}{0.000000,0.000000,0.000000}%
\pgfsetstrokecolor{currentstroke}%
\pgfsetstrokeopacity{0.000000}%
\pgfsetdash{}{0pt}%
\pgfpathmoveto{\pgfqpoint{0.588387in}{0.521603in}}%
\pgfpathlineto{\pgfqpoint{2.138387in}{0.521603in}}%
\pgfpathlineto{\pgfqpoint{2.138387in}{1.676603in}}%
\pgfpathlineto{\pgfqpoint{0.588387in}{1.676603in}}%
\pgfpathlineto{\pgfqpoint{0.588387in}{0.521603in}}%
\pgfpathclose%
\pgfusepath{fill}%
\end{pgfscope}%
\begin{pgfscope}%
\pgfpathrectangle{\pgfqpoint{0.588387in}{0.521603in}}{\pgfqpoint{1.550000in}{1.155000in}}%
\pgfusepath{clip}%
\pgfsetrectcap%
\pgfsetroundjoin%
\pgfsetlinewidth{0.803000pt}%
\definecolor{currentstroke}{rgb}{0.690196,0.690196,0.690196}%
\pgfsetstrokecolor{currentstroke}%
\pgfsetdash{}{0pt}%
\pgfpathmoveto{\pgfqpoint{1.403735in}{0.521603in}}%
\pgfpathlineto{\pgfqpoint{1.403735in}{1.676603in}}%
\pgfusepath{stroke}%
\end{pgfscope}%
\begin{pgfscope}%
\pgfsetbuttcap%
\pgfsetroundjoin%
\definecolor{currentfill}{rgb}{0.000000,0.000000,0.000000}%
\pgfsetfillcolor{currentfill}%
\pgfsetlinewidth{0.803000pt}%
\definecolor{currentstroke}{rgb}{0.000000,0.000000,0.000000}%
\pgfsetstrokecolor{currentstroke}%
\pgfsetdash{}{0pt}%
\pgfsys@defobject{currentmarker}{\pgfqpoint{0.000000in}{-0.048611in}}{\pgfqpoint{0.000000in}{0.000000in}}{%
\pgfpathmoveto{\pgfqpoint{0.000000in}{0.000000in}}%
\pgfpathlineto{\pgfqpoint{0.000000in}{-0.048611in}}%
\pgfusepath{stroke,fill}%
}%
\begin{pgfscope}%
\pgfsys@transformshift{1.403735in}{0.521603in}%
\pgfsys@useobject{currentmarker}{}%
\end{pgfscope}%
\end{pgfscope}%
\begin{pgfscope}%
\definecolor{textcolor}{rgb}{0.000000,0.000000,0.000000}%
\pgfsetstrokecolor{textcolor}%
\pgfsetfillcolor{textcolor}%
\pgftext[x=1.403735in,y=0.424381in,,top]{\color{textcolor}{\sffamily\fontsize{10.000000}{12.000000}\selectfont\catcode`\^=\active\def^{\ifmmode\sp\else\^{}\fi}\catcode`\%=\active\def
\end{pgfscope}%
\begin{pgfscope}%
\pgfsetbuttcap%
\pgfsetroundjoin%
\definecolor{currentfill}{rgb}{0.000000,0.000000,0.000000}%
\pgfsetfillcolor{currentfill}%
\pgfsetlinewidth{0.602250pt}%
\definecolor{currentstroke}{rgb}{0.000000,0.000000,0.000000}%
\pgfsetstrokecolor{currentstroke}%
\pgfsetdash{}{0pt}%
\pgfsys@defobject{currentmarker}{\pgfqpoint{0.000000in}{-0.027778in}}{\pgfqpoint{0.000000in}{0.000000in}}{%
\pgfpathmoveto{\pgfqpoint{0.000000in}{0.000000in}}%
\pgfpathlineto{\pgfqpoint{0.000000in}{-0.027778in}}%
\pgfusepath{stroke,fill}%
}%
\begin{pgfscope}%
\pgfsys@transformshift{0.658842in}{0.521603in}%
\pgfsys@useobject{currentmarker}{}%
\end{pgfscope}%
\end{pgfscope}%
\begin{pgfscope}%
\pgfsetbuttcap%
\pgfsetroundjoin%
\definecolor{currentfill}{rgb}{0.000000,0.000000,0.000000}%
\pgfsetfillcolor{currentfill}%
\pgfsetlinewidth{0.602250pt}%
\definecolor{currentstroke}{rgb}{0.000000,0.000000,0.000000}%
\pgfsetstrokecolor{currentstroke}%
\pgfsetdash{}{0pt}%
\pgfsys@defobject{currentmarker}{\pgfqpoint{0.000000in}{-0.027778in}}{\pgfqpoint{0.000000in}{0.000000in}}{%
\pgfpathmoveto{\pgfqpoint{0.000000in}{0.000000in}}%
\pgfpathlineto{\pgfqpoint{0.000000in}{-0.027778in}}%
\pgfusepath{stroke,fill}%
}%
\begin{pgfscope}%
\pgfsys@transformshift{0.846502in}{0.521603in}%
\pgfsys@useobject{currentmarker}{}%
\end{pgfscope}%
\end{pgfscope}%
\begin{pgfscope}%
\pgfsetbuttcap%
\pgfsetroundjoin%
\definecolor{currentfill}{rgb}{0.000000,0.000000,0.000000}%
\pgfsetfillcolor{currentfill}%
\pgfsetlinewidth{0.602250pt}%
\definecolor{currentstroke}{rgb}{0.000000,0.000000,0.000000}%
\pgfsetstrokecolor{currentstroke}%
\pgfsetdash{}{0pt}%
\pgfsys@defobject{currentmarker}{\pgfqpoint{0.000000in}{-0.027778in}}{\pgfqpoint{0.000000in}{0.000000in}}{%
\pgfpathmoveto{\pgfqpoint{0.000000in}{0.000000in}}%
\pgfpathlineto{\pgfqpoint{0.000000in}{-0.027778in}}%
\pgfusepath{stroke,fill}%
}%
\begin{pgfscope}%
\pgfsys@transformshift{0.979650in}{0.521603in}%
\pgfsys@useobject{currentmarker}{}%
\end{pgfscope}%
\end{pgfscope}%
\begin{pgfscope}%
\pgfsetbuttcap%
\pgfsetroundjoin%
\definecolor{currentfill}{rgb}{0.000000,0.000000,0.000000}%
\pgfsetfillcolor{currentfill}%
\pgfsetlinewidth{0.602250pt}%
\definecolor{currentstroke}{rgb}{0.000000,0.000000,0.000000}%
\pgfsetstrokecolor{currentstroke}%
\pgfsetdash{}{0pt}%
\pgfsys@defobject{currentmarker}{\pgfqpoint{0.000000in}{-0.027778in}}{\pgfqpoint{0.000000in}{0.000000in}}{%
\pgfpathmoveto{\pgfqpoint{0.000000in}{0.000000in}}%
\pgfpathlineto{\pgfqpoint{0.000000in}{-0.027778in}}%
\pgfusepath{stroke,fill}%
}%
\begin{pgfscope}%
\pgfsys@transformshift{1.082927in}{0.521603in}%
\pgfsys@useobject{currentmarker}{}%
\end{pgfscope}%
\end{pgfscope}%
\begin{pgfscope}%
\pgfsetbuttcap%
\pgfsetroundjoin%
\definecolor{currentfill}{rgb}{0.000000,0.000000,0.000000}%
\pgfsetfillcolor{currentfill}%
\pgfsetlinewidth{0.602250pt}%
\definecolor{currentstroke}{rgb}{0.000000,0.000000,0.000000}%
\pgfsetstrokecolor{currentstroke}%
\pgfsetdash{}{0pt}%
\pgfsys@defobject{currentmarker}{\pgfqpoint{0.000000in}{-0.027778in}}{\pgfqpoint{0.000000in}{0.000000in}}{%
\pgfpathmoveto{\pgfqpoint{0.000000in}{0.000000in}}%
\pgfpathlineto{\pgfqpoint{0.000000in}{-0.027778in}}%
\pgfusepath{stroke,fill}%
}%
\begin{pgfscope}%
\pgfsys@transformshift{1.167311in}{0.521603in}%
\pgfsys@useobject{currentmarker}{}%
\end{pgfscope}%
\end{pgfscope}%
\begin{pgfscope}%
\pgfsetbuttcap%
\pgfsetroundjoin%
\definecolor{currentfill}{rgb}{0.000000,0.000000,0.000000}%
\pgfsetfillcolor{currentfill}%
\pgfsetlinewidth{0.602250pt}%
\definecolor{currentstroke}{rgb}{0.000000,0.000000,0.000000}%
\pgfsetstrokecolor{currentstroke}%
\pgfsetdash{}{0pt}%
\pgfsys@defobject{currentmarker}{\pgfqpoint{0.000000in}{-0.027778in}}{\pgfqpoint{0.000000in}{0.000000in}}{%
\pgfpathmoveto{\pgfqpoint{0.000000in}{0.000000in}}%
\pgfpathlineto{\pgfqpoint{0.000000in}{-0.027778in}}%
\pgfusepath{stroke,fill}%
}%
\begin{pgfscope}%
\pgfsys@transformshift{1.238656in}{0.521603in}%
\pgfsys@useobject{currentmarker}{}%
\end{pgfscope}%
\end{pgfscope}%
\begin{pgfscope}%
\pgfsetbuttcap%
\pgfsetroundjoin%
\definecolor{currentfill}{rgb}{0.000000,0.000000,0.000000}%
\pgfsetfillcolor{currentfill}%
\pgfsetlinewidth{0.602250pt}%
\definecolor{currentstroke}{rgb}{0.000000,0.000000,0.000000}%
\pgfsetstrokecolor{currentstroke}%
\pgfsetdash{}{0pt}%
\pgfsys@defobject{currentmarker}{\pgfqpoint{0.000000in}{-0.027778in}}{\pgfqpoint{0.000000in}{0.000000in}}{%
\pgfpathmoveto{\pgfqpoint{0.000000in}{0.000000in}}%
\pgfpathlineto{\pgfqpoint{0.000000in}{-0.027778in}}%
\pgfusepath{stroke,fill}%
}%
\begin{pgfscope}%
\pgfsys@transformshift{1.300458in}{0.521603in}%
\pgfsys@useobject{currentmarker}{}%
\end{pgfscope}%
\end{pgfscope}%
\begin{pgfscope}%
\pgfsetbuttcap%
\pgfsetroundjoin%
\definecolor{currentfill}{rgb}{0.000000,0.000000,0.000000}%
\pgfsetfillcolor{currentfill}%
\pgfsetlinewidth{0.602250pt}%
\definecolor{currentstroke}{rgb}{0.000000,0.000000,0.000000}%
\pgfsetstrokecolor{currentstroke}%
\pgfsetdash{}{0pt}%
\pgfsys@defobject{currentmarker}{\pgfqpoint{0.000000in}{-0.027778in}}{\pgfqpoint{0.000000in}{0.000000in}}{%
\pgfpathmoveto{\pgfqpoint{0.000000in}{0.000000in}}%
\pgfpathlineto{\pgfqpoint{0.000000in}{-0.027778in}}%
\pgfusepath{stroke,fill}%
}%
\begin{pgfscope}%
\pgfsys@transformshift{1.354971in}{0.521603in}%
\pgfsys@useobject{currentmarker}{}%
\end{pgfscope}%
\end{pgfscope}%
\begin{pgfscope}%
\pgfsetbuttcap%
\pgfsetroundjoin%
\definecolor{currentfill}{rgb}{0.000000,0.000000,0.000000}%
\pgfsetfillcolor{currentfill}%
\pgfsetlinewidth{0.602250pt}%
\definecolor{currentstroke}{rgb}{0.000000,0.000000,0.000000}%
\pgfsetstrokecolor{currentstroke}%
\pgfsetdash{}{0pt}%
\pgfsys@defobject{currentmarker}{\pgfqpoint{0.000000in}{-0.027778in}}{\pgfqpoint{0.000000in}{0.000000in}}{%
\pgfpathmoveto{\pgfqpoint{0.000000in}{0.000000in}}%
\pgfpathlineto{\pgfqpoint{0.000000in}{-0.027778in}}%
\pgfusepath{stroke,fill}%
}%
\begin{pgfscope}%
\pgfsys@transformshift{1.724543in}{0.521603in}%
\pgfsys@useobject{currentmarker}{}%
\end{pgfscope}%
\end{pgfscope}%
\begin{pgfscope}%
\pgfsetbuttcap%
\pgfsetroundjoin%
\definecolor{currentfill}{rgb}{0.000000,0.000000,0.000000}%
\pgfsetfillcolor{currentfill}%
\pgfsetlinewidth{0.602250pt}%
\definecolor{currentstroke}{rgb}{0.000000,0.000000,0.000000}%
\pgfsetstrokecolor{currentstroke}%
\pgfsetdash{}{0pt}%
\pgfsys@defobject{currentmarker}{\pgfqpoint{0.000000in}{-0.027778in}}{\pgfqpoint{0.000000in}{0.000000in}}{%
\pgfpathmoveto{\pgfqpoint{0.000000in}{0.000000in}}%
\pgfpathlineto{\pgfqpoint{0.000000in}{-0.027778in}}%
\pgfusepath{stroke,fill}%
}%
\begin{pgfscope}%
\pgfsys@transformshift{1.912204in}{0.521603in}%
\pgfsys@useobject{currentmarker}{}%
\end{pgfscope}%
\end{pgfscope}%
\begin{pgfscope}%
\pgfsetbuttcap%
\pgfsetroundjoin%
\definecolor{currentfill}{rgb}{0.000000,0.000000,0.000000}%
\pgfsetfillcolor{currentfill}%
\pgfsetlinewidth{0.602250pt}%
\definecolor{currentstroke}{rgb}{0.000000,0.000000,0.000000}%
\pgfsetstrokecolor{currentstroke}%
\pgfsetdash{}{0pt}%
\pgfsys@defobject{currentmarker}{\pgfqpoint{0.000000in}{-0.027778in}}{\pgfqpoint{0.000000in}{0.000000in}}{%
\pgfpathmoveto{\pgfqpoint{0.000000in}{0.000000in}}%
\pgfpathlineto{\pgfqpoint{0.000000in}{-0.027778in}}%
\pgfusepath{stroke,fill}%
}%
\begin{pgfscope}%
\pgfsys@transformshift{2.045351in}{0.521603in}%
\pgfsys@useobject{currentmarker}{}%
\end{pgfscope}%
\end{pgfscope}%
\begin{pgfscope}%
\definecolor{textcolor}{rgb}{0.000000,0.000000,0.000000}%
\pgfsetstrokecolor{textcolor}%
\pgfsetfillcolor{textcolor}%
\pgftext[x=1.363387in,y=0.234413in,,top]{\color{textcolor}{\sffamily\fontsize{10.000000}{12.000000}\selectfont\catcode`\^=\active\def^{\ifmmode\sp\else\^{}\fi}\catcode`\%=\active\def
\end{pgfscope}%
\begin{pgfscope}%
\pgfpathrectangle{\pgfqpoint{0.588387in}{0.521603in}}{\pgfqpoint{1.550000in}{1.155000in}}%
\pgfusepath{clip}%
\pgfsetrectcap%
\pgfsetroundjoin%
\pgfsetlinewidth{0.803000pt}%
\definecolor{currentstroke}{rgb}{0.690196,0.690196,0.690196}%
\pgfsetstrokecolor{currentstroke}%
\pgfsetdash{}{0pt}%
\pgfpathmoveto{\pgfqpoint{0.588387in}{0.950063in}}%
\pgfpathlineto{\pgfqpoint{2.138387in}{0.950063in}}%
\pgfusepath{stroke}%
\end{pgfscope}%
\begin{pgfscope}%
\pgfsetbuttcap%
\pgfsetroundjoin%
\definecolor{currentfill}{rgb}{0.000000,0.000000,0.000000}%
\pgfsetfillcolor{currentfill}%
\pgfsetlinewidth{0.803000pt}%
\definecolor{currentstroke}{rgb}{0.000000,0.000000,0.000000}%
\pgfsetstrokecolor{currentstroke}%
\pgfsetdash{}{0pt}%
\pgfsys@defobject{currentmarker}{\pgfqpoint{-0.048611in}{0.000000in}}{\pgfqpoint{-0.000000in}{0.000000in}}{%
\pgfpathmoveto{\pgfqpoint{-0.000000in}{0.000000in}}%
\pgfpathlineto{\pgfqpoint{-0.048611in}{0.000000in}}%
\pgfusepath{stroke,fill}%
}%
\begin{pgfscope}%
\pgfsys@transformshift{0.588387in}{0.950063in}%
\pgfsys@useobject{currentmarker}{}%
\end{pgfscope}%
\end{pgfscope}%
\begin{pgfscope}%
\definecolor{textcolor}{rgb}{0.000000,0.000000,0.000000}%
\pgfsetstrokecolor{textcolor}%
\pgfsetfillcolor{textcolor}%
\pgftext[x=0.289968in, y=0.897301in, left, base]{\color{textcolor}{\sffamily\fontsize{10.000000}{12.000000}\selectfont\catcode`\^=\active\def^{\ifmmode\sp\else\^{}\fi}\catcode`\%=\active\def
\end{pgfscope}%
\begin{pgfscope}%
\pgfpathrectangle{\pgfqpoint{0.588387in}{0.521603in}}{\pgfqpoint{1.550000in}{1.155000in}}%
\pgfusepath{clip}%
\pgfsetrectcap%
\pgfsetroundjoin%
\pgfsetlinewidth{0.803000pt}%
\definecolor{currentstroke}{rgb}{0.690196,0.690196,0.690196}%
\pgfsetstrokecolor{currentstroke}%
\pgfsetdash{}{0pt}%
\pgfpathmoveto{\pgfqpoint{0.588387in}{1.580534in}}%
\pgfpathlineto{\pgfqpoint{2.138387in}{1.580534in}}%
\pgfusepath{stroke}%
\end{pgfscope}%
\begin{pgfscope}%
\pgfsetbuttcap%
\pgfsetroundjoin%
\definecolor{currentfill}{rgb}{0.000000,0.000000,0.000000}%
\pgfsetfillcolor{currentfill}%
\pgfsetlinewidth{0.803000pt}%
\definecolor{currentstroke}{rgb}{0.000000,0.000000,0.000000}%
\pgfsetstrokecolor{currentstroke}%
\pgfsetdash{}{0pt}%
\pgfsys@defobject{currentmarker}{\pgfqpoint{-0.048611in}{0.000000in}}{\pgfqpoint{-0.000000in}{0.000000in}}{%
\pgfpathmoveto{\pgfqpoint{-0.000000in}{0.000000in}}%
\pgfpathlineto{\pgfqpoint{-0.048611in}{0.000000in}}%
\pgfusepath{stroke,fill}%
}%
\begin{pgfscope}%
\pgfsys@transformshift{0.588387in}{1.580534in}%
\pgfsys@useobject{currentmarker}{}%
\end{pgfscope}%
\end{pgfscope}%
\begin{pgfscope}%
\definecolor{textcolor}{rgb}{0.000000,0.000000,0.000000}%
\pgfsetstrokecolor{textcolor}%
\pgfsetfillcolor{textcolor}%
\pgftext[x=0.289968in, y=1.527772in, left, base]{\color{textcolor}{\sffamily\fontsize{10.000000}{12.000000}\selectfont\catcode`\^=\active\def^{\ifmmode\sp\else\^{}\fi}\catcode`\%=\active\def
\end{pgfscope}%
\begin{pgfscope}%
\pgfsetbuttcap%
\pgfsetroundjoin%
\definecolor{currentfill}{rgb}{0.000000,0.000000,0.000000}%
\pgfsetfillcolor{currentfill}%
\pgfsetlinewidth{0.602250pt}%
\definecolor{currentstroke}{rgb}{0.000000,0.000000,0.000000}%
\pgfsetstrokecolor{currentstroke}%
\pgfsetdash{}{0pt}%
\pgfsys@defobject{currentmarker}{\pgfqpoint{-0.027778in}{0.000000in}}{\pgfqpoint{-0.000000in}{0.000000in}}{%
\pgfpathmoveto{\pgfqpoint{-0.000000in}{0.000000in}}%
\pgfpathlineto{\pgfqpoint{-0.027778in}{0.000000in}}%
\pgfusepath{stroke,fill}%
}%
\begin{pgfscope}%
\pgfsys@transformshift{0.588387in}{0.620403in}%
\pgfsys@useobject{currentmarker}{}%
\end{pgfscope}%
\end{pgfscope}%
\begin{pgfscope}%
\pgfsetbuttcap%
\pgfsetroundjoin%
\definecolor{currentfill}{rgb}{0.000000,0.000000,0.000000}%
\pgfsetfillcolor{currentfill}%
\pgfsetlinewidth{0.602250pt}%
\definecolor{currentstroke}{rgb}{0.000000,0.000000,0.000000}%
\pgfsetstrokecolor{currentstroke}%
\pgfsetdash{}{0pt}%
\pgfsys@defobject{currentmarker}{\pgfqpoint{-0.027778in}{0.000000in}}{\pgfqpoint{-0.000000in}{0.000000in}}{%
\pgfpathmoveto{\pgfqpoint{-0.000000in}{0.000000in}}%
\pgfpathlineto{\pgfqpoint{-0.027778in}{0.000000in}}%
\pgfusepath{stroke,fill}%
}%
\begin{pgfscope}%
\pgfsys@transformshift{0.588387in}{0.699173in}%
\pgfsys@useobject{currentmarker}{}%
\end{pgfscope}%
\end{pgfscope}%
\begin{pgfscope}%
\pgfsetbuttcap%
\pgfsetroundjoin%
\definecolor{currentfill}{rgb}{0.000000,0.000000,0.000000}%
\pgfsetfillcolor{currentfill}%
\pgfsetlinewidth{0.602250pt}%
\definecolor{currentstroke}{rgb}{0.000000,0.000000,0.000000}%
\pgfsetstrokecolor{currentstroke}%
\pgfsetdash{}{0pt}%
\pgfsys@defobject{currentmarker}{\pgfqpoint{-0.027778in}{0.000000in}}{\pgfqpoint{-0.000000in}{0.000000in}}{%
\pgfpathmoveto{\pgfqpoint{-0.000000in}{0.000000in}}%
\pgfpathlineto{\pgfqpoint{-0.027778in}{0.000000in}}%
\pgfusepath{stroke,fill}%
}%
\begin{pgfscope}%
\pgfsys@transformshift{0.588387in}{0.760272in}%
\pgfsys@useobject{currentmarker}{}%
\end{pgfscope}%
\end{pgfscope}%
\begin{pgfscope}%
\pgfsetbuttcap%
\pgfsetroundjoin%
\definecolor{currentfill}{rgb}{0.000000,0.000000,0.000000}%
\pgfsetfillcolor{currentfill}%
\pgfsetlinewidth{0.602250pt}%
\definecolor{currentstroke}{rgb}{0.000000,0.000000,0.000000}%
\pgfsetstrokecolor{currentstroke}%
\pgfsetdash{}{0pt}%
\pgfsys@defobject{currentmarker}{\pgfqpoint{-0.027778in}{0.000000in}}{\pgfqpoint{-0.000000in}{0.000000in}}{%
\pgfpathmoveto{\pgfqpoint{-0.000000in}{0.000000in}}%
\pgfpathlineto{\pgfqpoint{-0.027778in}{0.000000in}}%
\pgfusepath{stroke,fill}%
}%
\begin{pgfscope}%
\pgfsys@transformshift{0.588387in}{0.810193in}%
\pgfsys@useobject{currentmarker}{}%
\end{pgfscope}%
\end{pgfscope}%
\begin{pgfscope}%
\pgfsetbuttcap%
\pgfsetroundjoin%
\definecolor{currentfill}{rgb}{0.000000,0.000000,0.000000}%
\pgfsetfillcolor{currentfill}%
\pgfsetlinewidth{0.602250pt}%
\definecolor{currentstroke}{rgb}{0.000000,0.000000,0.000000}%
\pgfsetstrokecolor{currentstroke}%
\pgfsetdash{}{0pt}%
\pgfsys@defobject{currentmarker}{\pgfqpoint{-0.027778in}{0.000000in}}{\pgfqpoint{-0.000000in}{0.000000in}}{%
\pgfpathmoveto{\pgfqpoint{-0.000000in}{0.000000in}}%
\pgfpathlineto{\pgfqpoint{-0.027778in}{0.000000in}}%
\pgfusepath{stroke,fill}%
}%
\begin{pgfscope}%
\pgfsys@transformshift{0.588387in}{0.852401in}%
\pgfsys@useobject{currentmarker}{}%
\end{pgfscope}%
\end{pgfscope}%
\begin{pgfscope}%
\pgfsetbuttcap%
\pgfsetroundjoin%
\definecolor{currentfill}{rgb}{0.000000,0.000000,0.000000}%
\pgfsetfillcolor{currentfill}%
\pgfsetlinewidth{0.602250pt}%
\definecolor{currentstroke}{rgb}{0.000000,0.000000,0.000000}%
\pgfsetstrokecolor{currentstroke}%
\pgfsetdash{}{0pt}%
\pgfsys@defobject{currentmarker}{\pgfqpoint{-0.027778in}{0.000000in}}{\pgfqpoint{-0.000000in}{0.000000in}}{%
\pgfpathmoveto{\pgfqpoint{-0.000000in}{0.000000in}}%
\pgfpathlineto{\pgfqpoint{-0.027778in}{0.000000in}}%
\pgfusepath{stroke,fill}%
}%
\begin{pgfscope}%
\pgfsys@transformshift{0.588387in}{0.888964in}%
\pgfsys@useobject{currentmarker}{}%
\end{pgfscope}%
\end{pgfscope}%
\begin{pgfscope}%
\pgfsetbuttcap%
\pgfsetroundjoin%
\definecolor{currentfill}{rgb}{0.000000,0.000000,0.000000}%
\pgfsetfillcolor{currentfill}%
\pgfsetlinewidth{0.602250pt}%
\definecolor{currentstroke}{rgb}{0.000000,0.000000,0.000000}%
\pgfsetstrokecolor{currentstroke}%
\pgfsetdash{}{0pt}%
\pgfsys@defobject{currentmarker}{\pgfqpoint{-0.027778in}{0.000000in}}{\pgfqpoint{-0.000000in}{0.000000in}}{%
\pgfpathmoveto{\pgfqpoint{-0.000000in}{0.000000in}}%
\pgfpathlineto{\pgfqpoint{-0.027778in}{0.000000in}}%
\pgfusepath{stroke,fill}%
}%
\begin{pgfscope}%
\pgfsys@transformshift{0.588387in}{0.921214in}%
\pgfsys@useobject{currentmarker}{}%
\end{pgfscope}%
\end{pgfscope}%
\begin{pgfscope}%
\pgfsetbuttcap%
\pgfsetroundjoin%
\definecolor{currentfill}{rgb}{0.000000,0.000000,0.000000}%
\pgfsetfillcolor{currentfill}%
\pgfsetlinewidth{0.602250pt}%
\definecolor{currentstroke}{rgb}{0.000000,0.000000,0.000000}%
\pgfsetstrokecolor{currentstroke}%
\pgfsetdash{}{0pt}%
\pgfsys@defobject{currentmarker}{\pgfqpoint{-0.027778in}{0.000000in}}{\pgfqpoint{-0.000000in}{0.000000in}}{%
\pgfpathmoveto{\pgfqpoint{-0.000000in}{0.000000in}}%
\pgfpathlineto{\pgfqpoint{-0.027778in}{0.000000in}}%
\pgfusepath{stroke,fill}%
}%
\begin{pgfscope}%
\pgfsys@transformshift{0.588387in}{1.139853in}%
\pgfsys@useobject{currentmarker}{}%
\end{pgfscope}%
\end{pgfscope}%
\begin{pgfscope}%
\pgfsetbuttcap%
\pgfsetroundjoin%
\definecolor{currentfill}{rgb}{0.000000,0.000000,0.000000}%
\pgfsetfillcolor{currentfill}%
\pgfsetlinewidth{0.602250pt}%
\definecolor{currentstroke}{rgb}{0.000000,0.000000,0.000000}%
\pgfsetstrokecolor{currentstroke}%
\pgfsetdash{}{0pt}%
\pgfsys@defobject{currentmarker}{\pgfqpoint{-0.027778in}{0.000000in}}{\pgfqpoint{-0.000000in}{0.000000in}}{%
\pgfpathmoveto{\pgfqpoint{-0.000000in}{0.000000in}}%
\pgfpathlineto{\pgfqpoint{-0.027778in}{0.000000in}}%
\pgfusepath{stroke,fill}%
}%
\begin{pgfscope}%
\pgfsys@transformshift{0.588387in}{1.250874in}%
\pgfsys@useobject{currentmarker}{}%
\end{pgfscope}%
\end{pgfscope}%
\begin{pgfscope}%
\pgfsetbuttcap%
\pgfsetroundjoin%
\definecolor{currentfill}{rgb}{0.000000,0.000000,0.000000}%
\pgfsetfillcolor{currentfill}%
\pgfsetlinewidth{0.602250pt}%
\definecolor{currentstroke}{rgb}{0.000000,0.000000,0.000000}%
\pgfsetstrokecolor{currentstroke}%
\pgfsetdash{}{0pt}%
\pgfsys@defobject{currentmarker}{\pgfqpoint{-0.027778in}{0.000000in}}{\pgfqpoint{-0.000000in}{0.000000in}}{%
\pgfpathmoveto{\pgfqpoint{-0.000000in}{0.000000in}}%
\pgfpathlineto{\pgfqpoint{-0.027778in}{0.000000in}}%
\pgfusepath{stroke,fill}%
}%
\begin{pgfscope}%
\pgfsys@transformshift{0.588387in}{1.329644in}%
\pgfsys@useobject{currentmarker}{}%
\end{pgfscope}%
\end{pgfscope}%
\begin{pgfscope}%
\pgfsetbuttcap%
\pgfsetroundjoin%
\definecolor{currentfill}{rgb}{0.000000,0.000000,0.000000}%
\pgfsetfillcolor{currentfill}%
\pgfsetlinewidth{0.602250pt}%
\definecolor{currentstroke}{rgb}{0.000000,0.000000,0.000000}%
\pgfsetstrokecolor{currentstroke}%
\pgfsetdash{}{0pt}%
\pgfsys@defobject{currentmarker}{\pgfqpoint{-0.027778in}{0.000000in}}{\pgfqpoint{-0.000000in}{0.000000in}}{%
\pgfpathmoveto{\pgfqpoint{-0.000000in}{0.000000in}}%
\pgfpathlineto{\pgfqpoint{-0.027778in}{0.000000in}}%
\pgfusepath{stroke,fill}%
}%
\begin{pgfscope}%
\pgfsys@transformshift{0.588387in}{1.390743in}%
\pgfsys@useobject{currentmarker}{}%
\end{pgfscope}%
\end{pgfscope}%
\begin{pgfscope}%
\pgfsetbuttcap%
\pgfsetroundjoin%
\definecolor{currentfill}{rgb}{0.000000,0.000000,0.000000}%
\pgfsetfillcolor{currentfill}%
\pgfsetlinewidth{0.602250pt}%
\definecolor{currentstroke}{rgb}{0.000000,0.000000,0.000000}%
\pgfsetstrokecolor{currentstroke}%
\pgfsetdash{}{0pt}%
\pgfsys@defobject{currentmarker}{\pgfqpoint{-0.027778in}{0.000000in}}{\pgfqpoint{-0.000000in}{0.000000in}}{%
\pgfpathmoveto{\pgfqpoint{-0.000000in}{0.000000in}}%
\pgfpathlineto{\pgfqpoint{-0.027778in}{0.000000in}}%
\pgfusepath{stroke,fill}%
}%
\begin{pgfscope}%
\pgfsys@transformshift{0.588387in}{1.440664in}%
\pgfsys@useobject{currentmarker}{}%
\end{pgfscope}%
\end{pgfscope}%
\begin{pgfscope}%
\pgfsetbuttcap%
\pgfsetroundjoin%
\definecolor{currentfill}{rgb}{0.000000,0.000000,0.000000}%
\pgfsetfillcolor{currentfill}%
\pgfsetlinewidth{0.602250pt}%
\definecolor{currentstroke}{rgb}{0.000000,0.000000,0.000000}%
\pgfsetstrokecolor{currentstroke}%
\pgfsetdash{}{0pt}%
\pgfsys@defobject{currentmarker}{\pgfqpoint{-0.027778in}{0.000000in}}{\pgfqpoint{-0.000000in}{0.000000in}}{%
\pgfpathmoveto{\pgfqpoint{-0.000000in}{0.000000in}}%
\pgfpathlineto{\pgfqpoint{-0.027778in}{0.000000in}}%
\pgfusepath{stroke,fill}%
}%
\begin{pgfscope}%
\pgfsys@transformshift{0.588387in}{1.482872in}%
\pgfsys@useobject{currentmarker}{}%
\end{pgfscope}%
\end{pgfscope}%
\begin{pgfscope}%
\pgfsetbuttcap%
\pgfsetroundjoin%
\definecolor{currentfill}{rgb}{0.000000,0.000000,0.000000}%
\pgfsetfillcolor{currentfill}%
\pgfsetlinewidth{0.602250pt}%
\definecolor{currentstroke}{rgb}{0.000000,0.000000,0.000000}%
\pgfsetstrokecolor{currentstroke}%
\pgfsetdash{}{0pt}%
\pgfsys@defobject{currentmarker}{\pgfqpoint{-0.027778in}{0.000000in}}{\pgfqpoint{-0.000000in}{0.000000in}}{%
\pgfpathmoveto{\pgfqpoint{-0.000000in}{0.000000in}}%
\pgfpathlineto{\pgfqpoint{-0.027778in}{0.000000in}}%
\pgfusepath{stroke,fill}%
}%
\begin{pgfscope}%
\pgfsys@transformshift{0.588387in}{1.519435in}%
\pgfsys@useobject{currentmarker}{}%
\end{pgfscope}%
\end{pgfscope}%
\begin{pgfscope}%
\pgfsetbuttcap%
\pgfsetroundjoin%
\definecolor{currentfill}{rgb}{0.000000,0.000000,0.000000}%
\pgfsetfillcolor{currentfill}%
\pgfsetlinewidth{0.602250pt}%
\definecolor{currentstroke}{rgb}{0.000000,0.000000,0.000000}%
\pgfsetstrokecolor{currentstroke}%
\pgfsetdash{}{0pt}%
\pgfsys@defobject{currentmarker}{\pgfqpoint{-0.027778in}{0.000000in}}{\pgfqpoint{-0.000000in}{0.000000in}}{%
\pgfpathmoveto{\pgfqpoint{-0.000000in}{0.000000in}}%
\pgfpathlineto{\pgfqpoint{-0.027778in}{0.000000in}}%
\pgfusepath{stroke,fill}%
}%
\begin{pgfscope}%
\pgfsys@transformshift{0.588387in}{1.551685in}%
\pgfsys@useobject{currentmarker}{}%
\end{pgfscope}%
\end{pgfscope}%
\begin{pgfscope}%
\definecolor{textcolor}{rgb}{0.000000,0.000000,0.000000}%
\pgfsetstrokecolor{textcolor}%
\pgfsetfillcolor{textcolor}%
\pgftext[x=0.234413in,y=1.099103in,,bottom,rotate=90.000000]{\color{textcolor}{\sffamily\fontsize{10.000000}{12.000000}\selectfont\catcode`\^=\active\def^{\ifmmode\sp\else\^{}\fi}\catcode`\%=\active\def
\end{pgfscope}%
\begin{pgfscope}%
\pgfpathrectangle{\pgfqpoint{0.588387in}{0.521603in}}{\pgfqpoint{1.550000in}{1.155000in}}%
\pgfusepath{clip}%
\pgfsetrectcap%
\pgfsetroundjoin%
\pgfsetlinewidth{1.505625pt}%
\definecolor{currentstroke}{rgb}{0.121569,0.466667,0.705882}%
\pgfsetstrokecolor{currentstroke}%
\pgfsetdash{}{0pt}%
\pgfpathmoveto{\pgfqpoint{0.658842in}{0.574103in}}%
\pgfpathlineto{\pgfqpoint{1.167311in}{0.764655in}}%
\pgfpathlineto{\pgfqpoint{1.403735in}{0.923279in}}%
\pgfpathlineto{\pgfqpoint{1.559464in}{0.969908in}}%
\pgfpathlineto{\pgfqpoint{1.675779in}{1.082811in}}%
\pgfpathlineto{\pgfqpoint{1.768655in}{1.128042in}}%
\pgfpathlineto{\pgfqpoint{1.845973in}{1.224185in}}%
\pgfpathlineto{\pgfqpoint{1.912204in}{1.265538in}}%
\pgfpathlineto{\pgfqpoint{1.970133in}{1.252518in}}%
\pgfpathlineto{\pgfqpoint{2.021611in}{1.276414in}}%
\pgfpathlineto{\pgfqpoint{2.067933in}{1.302861in}}%
\pgfusepath{stroke}%
\end{pgfscope}%
\begin{pgfscope}%
\pgfpathrectangle{\pgfqpoint{0.588387in}{0.521603in}}{\pgfqpoint{1.550000in}{1.155000in}}%
\pgfusepath{clip}%
\pgfsetrectcap%
\pgfsetroundjoin%
\pgfsetlinewidth{1.505625pt}%
\definecolor{currentstroke}{rgb}{1.000000,0.498039,0.054902}%
\pgfsetstrokecolor{currentstroke}%
\pgfsetdash{}{0pt}%
\pgfpathmoveto{\pgfqpoint{0.658842in}{0.736554in}}%
\pgfpathlineto{\pgfqpoint{1.167311in}{0.978877in}}%
\pgfpathlineto{\pgfqpoint{1.403735in}{1.083363in}}%
\pgfpathlineto{\pgfqpoint{1.559464in}{1.170316in}}%
\pgfpathlineto{\pgfqpoint{1.675779in}{1.281594in}}%
\pgfpathlineto{\pgfqpoint{1.768655in}{1.299564in}}%
\pgfpathlineto{\pgfqpoint{1.845973in}{1.385796in}}%
\pgfpathlineto{\pgfqpoint{1.912204in}{1.417056in}}%
\pgfpathlineto{\pgfqpoint{1.970133in}{1.437518in}}%
\pgfpathlineto{\pgfqpoint{2.021611in}{1.463218in}}%
\pgfpathlineto{\pgfqpoint{2.067933in}{1.624103in}}%
\pgfusepath{stroke}%
\end{pgfscope}%
\begin{pgfscope}%
\pgfsetrectcap%
\pgfsetmiterjoin%
\pgfsetlinewidth{0.803000pt}%
\definecolor{currentstroke}{rgb}{0.000000,0.000000,0.000000}%
\pgfsetstrokecolor{currentstroke}%
\pgfsetdash{}{0pt}%
\pgfpathmoveto{\pgfqpoint{0.588387in}{0.521603in}}%
\pgfpathlineto{\pgfqpoint{0.588387in}{1.676603in}}%
\pgfusepath{stroke}%
\end{pgfscope}%
\begin{pgfscope}%
\pgfsetrectcap%
\pgfsetmiterjoin%
\pgfsetlinewidth{0.803000pt}%
\definecolor{currentstroke}{rgb}{0.000000,0.000000,0.000000}%
\pgfsetstrokecolor{currentstroke}%
\pgfsetdash{}{0pt}%
\pgfpathmoveto{\pgfqpoint{2.138387in}{0.521603in}}%
\pgfpathlineto{\pgfqpoint{2.138387in}{1.676603in}}%
\pgfusepath{stroke}%
\end{pgfscope}%
\begin{pgfscope}%
\pgfsetrectcap%
\pgfsetmiterjoin%
\pgfsetlinewidth{0.803000pt}%
\definecolor{currentstroke}{rgb}{0.000000,0.000000,0.000000}%
\pgfsetstrokecolor{currentstroke}%
\pgfsetdash{}{0pt}%
\pgfpathmoveto{\pgfqpoint{0.588387in}{0.521603in}}%
\pgfpathlineto{\pgfqpoint{2.138387in}{0.521603in}}%
\pgfusepath{stroke}%
\end{pgfscope}%
\begin{pgfscope}%
\pgfsetrectcap%
\pgfsetmiterjoin%
\pgfsetlinewidth{0.803000pt}%
\definecolor{currentstroke}{rgb}{0.000000,0.000000,0.000000}%
\pgfsetstrokecolor{currentstroke}%
\pgfsetdash{}{0pt}%
\pgfpathmoveto{\pgfqpoint{0.588387in}{1.676603in}}%
\pgfpathlineto{\pgfqpoint{2.138387in}{1.676603in}}%
\pgfusepath{stroke}%
\end{pgfscope}%
\end{pgfpicture}%
\makeatother%
\endgroup%

%% file: bench_vary_n_ratio.pgf
\begingroup%
\makeatletter%
\begin{pgfpicture}%
\pgfpathrectangle{\pgfpointorigin}{\pgfqpoint{2.368544in}{1.781079in}}%
\pgfusepath{use as bounding box, clip}%
\begin{pgfscope}%
\pgfsetbuttcap%
\pgfsetmiterjoin%
\definecolor{currentfill}{rgb}{1.000000,1.000000,1.000000}%
\pgfsetfillcolor{currentfill}%
\pgfsetlinewidth{0.000000pt}%
\definecolor{currentstroke}{rgb}{1.000000,1.000000,1.000000}%
\pgfsetstrokecolor{currentstroke}%
\pgfsetdash{}{0pt}%
\pgfpathmoveto{\pgfqpoint{0.000000in}{0.000000in}}%
\pgfpathlineto{\pgfqpoint{2.368544in}{0.000000in}}%
\pgfpathlineto{\pgfqpoint{2.368544in}{1.781079in}}%
\pgfpathlineto{\pgfqpoint{0.000000in}{1.781079in}}%
\pgfpathlineto{\pgfqpoint{0.000000in}{0.000000in}}%
\pgfpathclose%
\pgfusepath{fill}%
\end{pgfscope}%
\begin{pgfscope}%
\pgfsetbuttcap%
\pgfsetmiterjoin%
\definecolor{currentfill}{rgb}{1.000000,1.000000,1.000000}%
\pgfsetfillcolor{currentfill}%
\pgfsetlinewidth{0.000000pt}%
\definecolor{currentstroke}{rgb}{0.000000,0.000000,0.000000}%
\pgfsetstrokecolor{currentstroke}%
\pgfsetstrokeopacity{0.000000}%
\pgfsetdash{}{0pt}%
\pgfpathmoveto{\pgfqpoint{0.696435in}{0.526079in}}%
\pgfpathlineto{\pgfqpoint{2.246435in}{0.526079in}}%
\pgfpathlineto{\pgfqpoint{2.246435in}{1.681079in}}%
\pgfpathlineto{\pgfqpoint{0.696435in}{1.681079in}}%
\pgfpathlineto{\pgfqpoint{0.696435in}{0.526079in}}%
\pgfpathclose%
\pgfusepath{fill}%
\end{pgfscope}%
\begin{pgfscope}%
\pgfpathrectangle{\pgfqpoint{0.696435in}{0.526079in}}{\pgfqpoint{1.550000in}{1.155000in}}%
\pgfusepath{clip}%
\pgfsetrectcap%
\pgfsetroundjoin%
\pgfsetlinewidth{0.803000pt}%
\definecolor{currentstroke}{rgb}{0.690196,0.690196,0.690196}%
\pgfsetstrokecolor{currentstroke}%
\pgfsetdash{}{0pt}%
\pgfpathmoveto{\pgfqpoint{1.387796in}{0.526079in}}%
\pgfpathlineto{\pgfqpoint{1.387796in}{1.681079in}}%
\pgfusepath{stroke}%
\end{pgfscope}%
\begin{pgfscope}%
\pgfsetbuttcap%
\pgfsetroundjoin%
\definecolor{currentfill}{rgb}{0.000000,0.000000,0.000000}%
\pgfsetfillcolor{currentfill}%
\pgfsetlinewidth{0.803000pt}%
\definecolor{currentstroke}{rgb}{0.000000,0.000000,0.000000}%
\pgfsetstrokecolor{currentstroke}%
\pgfsetdash{}{0pt}%
\pgfsys@defobject{currentmarker}{\pgfqpoint{0.000000in}{-0.048611in}}{\pgfqpoint{0.000000in}{0.000000in}}{%
\pgfpathmoveto{\pgfqpoint{0.000000in}{0.000000in}}%
\pgfpathlineto{\pgfqpoint{0.000000in}{-0.048611in}}%
\pgfusepath{stroke,fill}%
}%
\begin{pgfscope}%
\pgfsys@transformshift{1.387796in}{0.526079in}%
\pgfsys@useobject{currentmarker}{}%
\end{pgfscope}%
\end{pgfscope}%
\begin{pgfscope}%
\definecolor{textcolor}{rgb}{0.000000,0.000000,0.000000}%
\pgfsetstrokecolor{textcolor}%
\pgfsetfillcolor{textcolor}%
\pgftext[x=1.387796in,y=0.428857in,,top]{\color{textcolor}{\sffamily\fontsize{10.000000}{12.000000}\selectfont\catcode`\^=\active\def^{\ifmmode\sp\else\^{}\fi}\catcode`\%=\active\def
\end{pgfscope}%
\begin{pgfscope}%
\pgfpathrectangle{\pgfqpoint{0.696435in}{0.526079in}}{\pgfqpoint{1.550000in}{1.155000in}}%
\pgfusepath{clip}%
\pgfsetrectcap%
\pgfsetroundjoin%
\pgfsetlinewidth{0.803000pt}%
\definecolor{currentstroke}{rgb}{0.690196,0.690196,0.690196}%
\pgfsetstrokecolor{currentstroke}%
\pgfsetdash{}{0pt}%
\pgfpathmoveto{\pgfqpoint{2.167945in}{0.526079in}}%
\pgfpathlineto{\pgfqpoint{2.167945in}{1.681079in}}%
\pgfusepath{stroke}%
\end{pgfscope}%
\begin{pgfscope}%
\pgfsetbuttcap%
\pgfsetroundjoin%
\definecolor{currentfill}{rgb}{0.000000,0.000000,0.000000}%
\pgfsetfillcolor{currentfill}%
\pgfsetlinewidth{0.803000pt}%
\definecolor{currentstroke}{rgb}{0.000000,0.000000,0.000000}%
\pgfsetstrokecolor{currentstroke}%
\pgfsetdash{}{0pt}%
\pgfsys@defobject{currentmarker}{\pgfqpoint{0.000000in}{-0.048611in}}{\pgfqpoint{0.000000in}{0.000000in}}{%
\pgfpathmoveto{\pgfqpoint{0.000000in}{0.000000in}}%
\pgfpathlineto{\pgfqpoint{0.000000in}{-0.048611in}}%
\pgfusepath{stroke,fill}%
}%
\begin{pgfscope}%
\pgfsys@transformshift{2.167945in}{0.526079in}%
\pgfsys@useobject{currentmarker}{}%
\end{pgfscope}%
\end{pgfscope}%
\begin{pgfscope}%
\definecolor{textcolor}{rgb}{0.000000,0.000000,0.000000}%
\pgfsetstrokecolor{textcolor}%
\pgfsetfillcolor{textcolor}%
\pgftext[x=2.167945in,y=0.428857in,,top]{\color{textcolor}{\sffamily\fontsize{10.000000}{12.000000}\selectfont\catcode`\^=\active\def^{\ifmmode\sp\else\^{}\fi}\catcode`\%=\active\def
\end{pgfscope}%
\begin{pgfscope}%
\definecolor{textcolor}{rgb}{0.000000,0.000000,0.000000}%
\pgfsetstrokecolor{textcolor}%
\pgfsetfillcolor{textcolor}%
\pgftext[x=1.471435in,y=0.238889in,,top]{\color{textcolor}{\sffamily\fontsize{10.000000}{12.000000}\selectfont\catcode`\^=\active\def^{\ifmmode\sp\else\^{}\fi}\catcode`\%=\active\def
\end{pgfscope}%
\begin{pgfscope}%
\pgfpathrectangle{\pgfqpoint{0.696435in}{0.526079in}}{\pgfqpoint{1.550000in}{1.155000in}}%
\pgfusepath{clip}%
\pgfsetrectcap%
\pgfsetroundjoin%
\pgfsetlinewidth{0.803000pt}%
\definecolor{currentstroke}{rgb}{0.690196,0.690196,0.690196}%
\pgfsetstrokecolor{currentstroke}%
\pgfsetdash{}{0pt}%
\pgfpathmoveto{\pgfqpoint{0.696435in}{0.672918in}}%
\pgfpathlineto{\pgfqpoint{2.246435in}{0.672918in}}%
\pgfusepath{stroke}%
\end{pgfscope}%
\begin{pgfscope}%
\pgfsetbuttcap%
\pgfsetroundjoin%
\definecolor{currentfill}{rgb}{0.000000,0.000000,0.000000}%
\pgfsetfillcolor{currentfill}%
\pgfsetlinewidth{0.803000pt}%
\definecolor{currentstroke}{rgb}{0.000000,0.000000,0.000000}%
\pgfsetstrokecolor{currentstroke}%
\pgfsetdash{}{0pt}%
\pgfsys@defobject{currentmarker}{\pgfqpoint{-0.048611in}{0.000000in}}{\pgfqpoint{-0.000000in}{0.000000in}}{%
\pgfpathmoveto{\pgfqpoint{-0.000000in}{0.000000in}}%
\pgfpathlineto{\pgfqpoint{-0.048611in}{0.000000in}}%
\pgfusepath{stroke,fill}%
}%
\begin{pgfscope}%
\pgfsys@transformshift{0.696435in}{0.672918in}%
\pgfsys@useobject{currentmarker}{}%
\end{pgfscope}%
\end{pgfscope}%
\begin{pgfscope}%
\definecolor{textcolor}{rgb}{0.000000,0.000000,0.000000}%
\pgfsetstrokecolor{textcolor}%
\pgfsetfillcolor{textcolor}%
\pgftext[x=0.289968in, y=0.620156in, left, base]{\color{textcolor}{\sffamily\fontsize{10.000000}{12.000000}\selectfont\catcode`\^=\active\def^{\ifmmode\sp\else\^{}\fi}\catcode`\%=\active\def
\end{pgfscope}%
\begin{pgfscope}%
\pgfpathrectangle{\pgfqpoint{0.696435in}{0.526079in}}{\pgfqpoint{1.550000in}{1.155000in}}%
\pgfusepath{clip}%
\pgfsetrectcap%
\pgfsetroundjoin%
\pgfsetlinewidth{0.803000pt}%
\definecolor{currentstroke}{rgb}{0.690196,0.690196,0.690196}%
\pgfsetstrokecolor{currentstroke}%
\pgfsetdash{}{0pt}%
\pgfpathmoveto{\pgfqpoint{0.696435in}{0.966882in}}%
\pgfpathlineto{\pgfqpoint{2.246435in}{0.966882in}}%
\pgfusepath{stroke}%
\end{pgfscope}%
\begin{pgfscope}%
\pgfsetbuttcap%
\pgfsetroundjoin%
\definecolor{currentfill}{rgb}{0.000000,0.000000,0.000000}%
\pgfsetfillcolor{currentfill}%
\pgfsetlinewidth{0.803000pt}%
\definecolor{currentstroke}{rgb}{0.000000,0.000000,0.000000}%
\pgfsetstrokecolor{currentstroke}%
\pgfsetdash{}{0pt}%
\pgfsys@defobject{currentmarker}{\pgfqpoint{-0.048611in}{0.000000in}}{\pgfqpoint{-0.000000in}{0.000000in}}{%
\pgfpathmoveto{\pgfqpoint{-0.000000in}{0.000000in}}%
\pgfpathlineto{\pgfqpoint{-0.048611in}{0.000000in}}%
\pgfusepath{stroke,fill}%
}%
\begin{pgfscope}%
\pgfsys@transformshift{0.696435in}{0.966882in}%
\pgfsys@useobject{currentmarker}{}%
\end{pgfscope}%
\end{pgfscope}%
\begin{pgfscope}%
\definecolor{textcolor}{rgb}{0.000000,0.000000,0.000000}%
\pgfsetstrokecolor{textcolor}%
\pgfsetfillcolor{textcolor}%
\pgftext[x=0.289968in, y=0.914121in, left, base]{\color{textcolor}{\sffamily\fontsize{10.000000}{12.000000}\selectfont\catcode`\^=\active\def^{\ifmmode\sp\else\^{}\fi}\catcode`\%=\active\def
\end{pgfscope}%
\begin{pgfscope}%
\pgfpathrectangle{\pgfqpoint{0.696435in}{0.526079in}}{\pgfqpoint{1.550000in}{1.155000in}}%
\pgfusepath{clip}%
\pgfsetrectcap%
\pgfsetroundjoin%
\pgfsetlinewidth{0.803000pt}%
\definecolor{currentstroke}{rgb}{0.690196,0.690196,0.690196}%
\pgfsetstrokecolor{currentstroke}%
\pgfsetdash{}{0pt}%
\pgfpathmoveto{\pgfqpoint{0.696435in}{1.260847in}}%
\pgfpathlineto{\pgfqpoint{2.246435in}{1.260847in}}%
\pgfusepath{stroke}%
\end{pgfscope}%
\begin{pgfscope}%
\pgfsetbuttcap%
\pgfsetroundjoin%
\definecolor{currentfill}{rgb}{0.000000,0.000000,0.000000}%
\pgfsetfillcolor{currentfill}%
\pgfsetlinewidth{0.803000pt}%
\definecolor{currentstroke}{rgb}{0.000000,0.000000,0.000000}%
\pgfsetstrokecolor{currentstroke}%
\pgfsetdash{}{0pt}%
\pgfsys@defobject{currentmarker}{\pgfqpoint{-0.048611in}{0.000000in}}{\pgfqpoint{-0.000000in}{0.000000in}}{%
\pgfpathmoveto{\pgfqpoint{-0.000000in}{0.000000in}}%
\pgfpathlineto{\pgfqpoint{-0.048611in}{0.000000in}}%
\pgfusepath{stroke,fill}%
}%
\begin{pgfscope}%
\pgfsys@transformshift{0.696435in}{1.260847in}%
\pgfsys@useobject{currentmarker}{}%
\end{pgfscope}%
\end{pgfscope}%
\begin{pgfscope}%
\definecolor{textcolor}{rgb}{0.000000,0.000000,0.000000}%
\pgfsetstrokecolor{textcolor}%
\pgfsetfillcolor{textcolor}%
\pgftext[x=0.289968in, y=1.208085in, left, base]{\color{textcolor}{\sffamily\fontsize{10.000000}{12.000000}\selectfont\catcode`\^=\active\def^{\ifmmode\sp\else\^{}\fi}\catcode`\%=\active\def
\end{pgfscope}%
\begin{pgfscope}%
\pgfpathrectangle{\pgfqpoint{0.696435in}{0.526079in}}{\pgfqpoint{1.550000in}{1.155000in}}%
\pgfusepath{clip}%
\pgfsetrectcap%
\pgfsetroundjoin%
\pgfsetlinewidth{0.803000pt}%
\definecolor{currentstroke}{rgb}{0.690196,0.690196,0.690196}%
\pgfsetstrokecolor{currentstroke}%
\pgfsetdash{}{0pt}%
\pgfpathmoveto{\pgfqpoint{0.696435in}{1.554811in}}%
\pgfpathlineto{\pgfqpoint{2.246435in}{1.554811in}}%
\pgfusepath{stroke}%
\end{pgfscope}%
\begin{pgfscope}%
\pgfsetbuttcap%
\pgfsetroundjoin%
\definecolor{currentfill}{rgb}{0.000000,0.000000,0.000000}%
\pgfsetfillcolor{currentfill}%
\pgfsetlinewidth{0.803000pt}%
\definecolor{currentstroke}{rgb}{0.000000,0.000000,0.000000}%
\pgfsetstrokecolor{currentstroke}%
\pgfsetdash{}{0pt}%
\pgfsys@defobject{currentmarker}{\pgfqpoint{-0.048611in}{0.000000in}}{\pgfqpoint{-0.000000in}{0.000000in}}{%
\pgfpathmoveto{\pgfqpoint{-0.000000in}{0.000000in}}%
\pgfpathlineto{\pgfqpoint{-0.048611in}{0.000000in}}%
\pgfusepath{stroke,fill}%
}%
\begin{pgfscope}%
\pgfsys@transformshift{0.696435in}{1.554811in}%
\pgfsys@useobject{currentmarker}{}%
\end{pgfscope}%
\end{pgfscope}%
\begin{pgfscope}%
\definecolor{textcolor}{rgb}{0.000000,0.000000,0.000000}%
\pgfsetstrokecolor{textcolor}%
\pgfsetfillcolor{textcolor}%
\pgftext[x=0.289968in, y=1.502049in, left, base]{\color{textcolor}{\sffamily\fontsize{10.000000}{12.000000}\selectfont\catcode`\^=\active\def^{\ifmmode\sp\else\^{}\fi}\catcode`\%=\active\def
\end{pgfscope}%
\begin{pgfscope}%
\definecolor{textcolor}{rgb}{0.000000,0.000000,0.000000}%
\pgfsetstrokecolor{textcolor}%
\pgfsetfillcolor{textcolor}%
\pgftext[x=0.234413in,y=1.103579in,,bottom,rotate=90.000000]{\color{textcolor}{\sffamily\fontsize{10.000000}{12.000000}\selectfont\catcode`\^=\active\def^{\ifmmode\sp\else\^{}\fi}\catcode`\%=\active\def
\end{pgfscope}%
\begin{pgfscope}%
\pgfpathrectangle{\pgfqpoint{0.696435in}{0.526079in}}{\pgfqpoint{1.550000in}{1.155000in}}%
\pgfusepath{clip}%
\pgfsetrectcap%
\pgfsetroundjoin%
\pgfsetlinewidth{1.505625pt}%
\definecolor{currentstroke}{rgb}{0.121569,0.466667,0.705882}%
\pgfsetstrokecolor{currentstroke}%
\pgfsetdash{}{0pt}%
\pgfpathmoveto{\pgfqpoint{0.766890in}{0.744672in}}%
\pgfpathlineto{\pgfqpoint{1.001738in}{0.578579in}}%
\pgfpathlineto{\pgfqpoint{1.236587in}{1.335751in}}%
\pgfpathlineto{\pgfqpoint{1.471435in}{1.628579in}}%
\pgfpathlineto{\pgfqpoint{1.706284in}{1.124956in}}%
\pgfpathlineto{\pgfqpoint{1.941132in}{1.128208in}}%
\pgfpathlineto{\pgfqpoint{2.175981in}{1.089964in}}%
\pgfusepath{stroke}%
\end{pgfscope}%
\begin{pgfscope}%
\pgfsetrectcap%
\pgfsetmiterjoin%
\pgfsetlinewidth{0.803000pt}%
\definecolor{currentstroke}{rgb}{0.000000,0.000000,0.000000}%
\pgfsetstrokecolor{currentstroke}%
\pgfsetdash{}{0pt}%
\pgfpathmoveto{\pgfqpoint{0.696435in}{0.526079in}}%
\pgfpathlineto{\pgfqpoint{0.696435in}{1.681079in}}%
\pgfusepath{stroke}%
\end{pgfscope}%
\begin{pgfscope}%
\pgfsetrectcap%
\pgfsetmiterjoin%
\pgfsetlinewidth{0.803000pt}%
\definecolor{currentstroke}{rgb}{0.000000,0.000000,0.000000}%
\pgfsetstrokecolor{currentstroke}%
\pgfsetdash{}{0pt}%
\pgfpathmoveto{\pgfqpoint{2.246435in}{0.526079in}}%
\pgfpathlineto{\pgfqpoint{2.246435in}{1.681079in}}%
\pgfusepath{stroke}%
\end{pgfscope}%
\begin{pgfscope}%
\pgfsetrectcap%
\pgfsetmiterjoin%
\pgfsetlinewidth{0.803000pt}%
\definecolor{currentstroke}{rgb}{0.000000,0.000000,0.000000}%
\pgfsetstrokecolor{currentstroke}%
\pgfsetdash{}{0pt}%
\pgfpathmoveto{\pgfqpoint{0.696435in}{0.526079in}}%
\pgfpathlineto{\pgfqpoint{2.246435in}{0.526079in}}%
\pgfusepath{stroke}%
\end{pgfscope}%
\begin{pgfscope}%
\pgfsetrectcap%
\pgfsetmiterjoin%
\pgfsetlinewidth{0.803000pt}%
\definecolor{currentstroke}{rgb}{0.000000,0.000000,0.000000}%
\pgfsetstrokecolor{currentstroke}%
\pgfsetdash{}{0pt}%
\pgfpathmoveto{\pgfqpoint{0.696435in}{1.681079in}}%
\pgfpathlineto{\pgfqpoint{2.246435in}{1.681079in}}%
\pgfusepath{stroke}%
\end{pgfscope}%
\end{pgfpicture}%
\makeatother%
\endgroup%

%% file: bench_vary_nw_ratio.pgf
\begingroup%
\makeatletter%
\begin{pgfpicture}%
\pgfpathrectangle{\pgfpointorigin}{\pgfqpoint{2.258070in}{1.776603in}}%
\pgfusepath{use as bounding box, clip}%
\begin{pgfscope}%
\pgfsetbuttcap%
\pgfsetmiterjoin%
\definecolor{currentfill}{rgb}{1.000000,1.000000,1.000000}%
\pgfsetfillcolor{currentfill}%
\pgfsetlinewidth{0.000000pt}%
\definecolor{currentstroke}{rgb}{1.000000,1.000000,1.000000}%
\pgfsetstrokecolor{currentstroke}%
\pgfsetdash{}{0pt}%
\pgfpathmoveto{\pgfqpoint{0.000000in}{0.000000in}}%
\pgfpathlineto{\pgfqpoint{2.258070in}{0.000000in}}%
\pgfpathlineto{\pgfqpoint{2.258070in}{1.776603in}}%
\pgfpathlineto{\pgfqpoint{0.000000in}{1.776603in}}%
\pgfpathlineto{\pgfqpoint{0.000000in}{0.000000in}}%
\pgfpathclose%
\pgfusepath{fill}%
\end{pgfscope}%
\begin{pgfscope}%
\pgfsetbuttcap%
\pgfsetmiterjoin%
\definecolor{currentfill}{rgb}{1.000000,1.000000,1.000000}%
\pgfsetfillcolor{currentfill}%
\pgfsetlinewidth{0.000000pt}%
\definecolor{currentstroke}{rgb}{0.000000,0.000000,0.000000}%
\pgfsetstrokecolor{currentstroke}%
\pgfsetstrokeopacity{0.000000}%
\pgfsetdash{}{0pt}%
\pgfpathmoveto{\pgfqpoint{0.608070in}{0.521603in}}%
\pgfpathlineto{\pgfqpoint{2.158070in}{0.521603in}}%
\pgfpathlineto{\pgfqpoint{2.158070in}{1.676603in}}%
\pgfpathlineto{\pgfqpoint{0.608070in}{1.676603in}}%
\pgfpathlineto{\pgfqpoint{0.608070in}{0.521603in}}%
\pgfpathclose%
\pgfusepath{fill}%
\end{pgfscope}%
\begin{pgfscope}%
\pgfpathrectangle{\pgfqpoint{0.608070in}{0.521603in}}{\pgfqpoint{1.550000in}{1.155000in}}%
\pgfusepath{clip}%
\pgfsetrectcap%
\pgfsetroundjoin%
\pgfsetlinewidth{0.803000pt}%
\definecolor{currentstroke}{rgb}{0.690196,0.690196,0.690196}%
\pgfsetstrokecolor{currentstroke}%
\pgfsetdash{}{0pt}%
\pgfpathmoveto{\pgfqpoint{0.608070in}{0.521603in}}%
\pgfpathlineto{\pgfqpoint{0.608070in}{1.676603in}}%
\pgfusepath{stroke}%
\end{pgfscope}%
\begin{pgfscope}%
\pgfsetbuttcap%
\pgfsetroundjoin%
\definecolor{currentfill}{rgb}{0.000000,0.000000,0.000000}%
\pgfsetfillcolor{currentfill}%
\pgfsetlinewidth{0.803000pt}%
\definecolor{currentstroke}{rgb}{0.000000,0.000000,0.000000}%
\pgfsetstrokecolor{currentstroke}%
\pgfsetdash{}{0pt}%
\pgfsys@defobject{currentmarker}{\pgfqpoint{0.000000in}{-0.048611in}}{\pgfqpoint{0.000000in}{0.000000in}}{%
\pgfpathmoveto{\pgfqpoint{0.000000in}{0.000000in}}%
\pgfpathlineto{\pgfqpoint{0.000000in}{-0.048611in}}%
\pgfusepath{stroke,fill}%
}%
\begin{pgfscope}%
\pgfsys@transformshift{0.608070in}{0.521603in}%
\pgfsys@useobject{currentmarker}{}%
\end{pgfscope}%
\end{pgfscope}%
\begin{pgfscope}%
\definecolor{textcolor}{rgb}{0.000000,0.000000,0.000000}%
\pgfsetstrokecolor{textcolor}%
\pgfsetfillcolor{textcolor}%
\pgftext[x=0.608070in,y=0.424381in,,top]{\color{textcolor}{\sffamily\fontsize{10.000000}{12.000000}\selectfont\catcode`\^=\active\def^{\ifmmode\sp\else\^{}\fi}\catcode`\%=\active\def
\end{pgfscope}%
\begin{pgfscope}%
\pgfpathrectangle{\pgfqpoint{0.608070in}{0.521603in}}{\pgfqpoint{1.550000in}{1.155000in}}%
\pgfusepath{clip}%
\pgfsetrectcap%
\pgfsetroundjoin%
\pgfsetlinewidth{0.803000pt}%
\definecolor{currentstroke}{rgb}{0.690196,0.690196,0.690196}%
\pgfsetstrokecolor{currentstroke}%
\pgfsetdash{}{0pt}%
\pgfpathmoveto{\pgfqpoint{1.312615in}{0.521603in}}%
\pgfpathlineto{\pgfqpoint{1.312615in}{1.676603in}}%
\pgfusepath{stroke}%
\end{pgfscope}%
\begin{pgfscope}%
\pgfsetbuttcap%
\pgfsetroundjoin%
\definecolor{currentfill}{rgb}{0.000000,0.000000,0.000000}%
\pgfsetfillcolor{currentfill}%
\pgfsetlinewidth{0.803000pt}%
\definecolor{currentstroke}{rgb}{0.000000,0.000000,0.000000}%
\pgfsetstrokecolor{currentstroke}%
\pgfsetdash{}{0pt}%
\pgfsys@defobject{currentmarker}{\pgfqpoint{0.000000in}{-0.048611in}}{\pgfqpoint{0.000000in}{0.000000in}}{%
\pgfpathmoveto{\pgfqpoint{0.000000in}{0.000000in}}%
\pgfpathlineto{\pgfqpoint{0.000000in}{-0.048611in}}%
\pgfusepath{stroke,fill}%
}%
\begin{pgfscope}%
\pgfsys@transformshift{1.312615in}{0.521603in}%
\pgfsys@useobject{currentmarker}{}%
\end{pgfscope}%
\end{pgfscope}%
\begin{pgfscope}%
\definecolor{textcolor}{rgb}{0.000000,0.000000,0.000000}%
\pgfsetstrokecolor{textcolor}%
\pgfsetfillcolor{textcolor}%
\pgftext[x=1.312615in,y=0.424381in,,top]{\color{textcolor}{\sffamily\fontsize{10.000000}{12.000000}\selectfont\catcode`\^=\active\def^{\ifmmode\sp\else\^{}\fi}\catcode`\%=\active\def
\end{pgfscope}%
\begin{pgfscope}%
\pgfpathrectangle{\pgfqpoint{0.608070in}{0.521603in}}{\pgfqpoint{1.550000in}{1.155000in}}%
\pgfusepath{clip}%
\pgfsetrectcap%
\pgfsetroundjoin%
\pgfsetlinewidth{0.803000pt}%
\definecolor{currentstroke}{rgb}{0.690196,0.690196,0.690196}%
\pgfsetstrokecolor{currentstroke}%
\pgfsetdash{}{0pt}%
\pgfpathmoveto{\pgfqpoint{2.017161in}{0.521603in}}%
\pgfpathlineto{\pgfqpoint{2.017161in}{1.676603in}}%
\pgfusepath{stroke}%
\end{pgfscope}%
\begin{pgfscope}%
\pgfsetbuttcap%
\pgfsetroundjoin%
\definecolor{currentfill}{rgb}{0.000000,0.000000,0.000000}%
\pgfsetfillcolor{currentfill}%
\pgfsetlinewidth{0.803000pt}%
\definecolor{currentstroke}{rgb}{0.000000,0.000000,0.000000}%
\pgfsetstrokecolor{currentstroke}%
\pgfsetdash{}{0pt}%
\pgfsys@defobject{currentmarker}{\pgfqpoint{0.000000in}{-0.048611in}}{\pgfqpoint{0.000000in}{0.000000in}}{%
\pgfpathmoveto{\pgfqpoint{0.000000in}{0.000000in}}%
\pgfpathlineto{\pgfqpoint{0.000000in}{-0.048611in}}%
\pgfusepath{stroke,fill}%
}%
\begin{pgfscope}%
\pgfsys@transformshift{2.017161in}{0.521603in}%
\pgfsys@useobject{currentmarker}{}%
\end{pgfscope}%
\end{pgfscope}%
\begin{pgfscope}%
\definecolor{textcolor}{rgb}{0.000000,0.000000,0.000000}%
\pgfsetstrokecolor{textcolor}%
\pgfsetfillcolor{textcolor}%
\pgftext[x=2.017161in,y=0.424381in,,top]{\color{textcolor}{\sffamily\fontsize{10.000000}{12.000000}\selectfont\catcode`\^=\active\def^{\ifmmode\sp\else\^{}\fi}\catcode`\%=\active\def
\end{pgfscope}%
\begin{pgfscope}%
\definecolor{textcolor}{rgb}{0.000000,0.000000,0.000000}%
\pgfsetstrokecolor{textcolor}%
\pgfsetfillcolor{textcolor}%
\pgftext[x=1.383070in,y=0.234413in,,top]{\color{textcolor}{\sffamily\fontsize{10.000000}{12.000000}\selectfont\catcode`\^=\active\def^{\ifmmode\sp\else\^{}\fi}\catcode`\%=\active\def
\end{pgfscope}%
\begin{pgfscope}%
\pgfpathrectangle{\pgfqpoint{0.608070in}{0.521603in}}{\pgfqpoint{1.550000in}{1.155000in}}%
\pgfusepath{clip}%
\pgfsetrectcap%
\pgfsetroundjoin%
\pgfsetlinewidth{0.803000pt}%
\definecolor{currentstroke}{rgb}{0.690196,0.690196,0.690196}%
\pgfsetstrokecolor{currentstroke}%
\pgfsetdash{}{0pt}%
\pgfpathmoveto{\pgfqpoint{0.608070in}{0.757549in}}%
\pgfpathlineto{\pgfqpoint{2.158070in}{0.757549in}}%
\pgfusepath{stroke}%
\end{pgfscope}%
\begin{pgfscope}%
\pgfsetbuttcap%
\pgfsetroundjoin%
\definecolor{currentfill}{rgb}{0.000000,0.000000,0.000000}%
\pgfsetfillcolor{currentfill}%
\pgfsetlinewidth{0.803000pt}%
\definecolor{currentstroke}{rgb}{0.000000,0.000000,0.000000}%
\pgfsetstrokecolor{currentstroke}%
\pgfsetdash{}{0pt}%
\pgfsys@defobject{currentmarker}{\pgfqpoint{-0.048611in}{0.000000in}}{\pgfqpoint{-0.000000in}{0.000000in}}{%
\pgfpathmoveto{\pgfqpoint{-0.000000in}{0.000000in}}%
\pgfpathlineto{\pgfqpoint{-0.048611in}{0.000000in}}%
\pgfusepath{stroke,fill}%
}%
\begin{pgfscope}%
\pgfsys@transformshift{0.608070in}{0.757549in}%
\pgfsys@useobject{currentmarker}{}%
\end{pgfscope}%
\end{pgfscope}%
\begin{pgfscope}%
\definecolor{textcolor}{rgb}{0.000000,0.000000,0.000000}%
\pgfsetstrokecolor{textcolor}%
\pgfsetfillcolor{textcolor}%
\pgftext[x=0.289968in, y=0.704788in, left, base]{\color{textcolor}{\sffamily\fontsize{10.000000}{12.000000}\selectfont\catcode`\^=\active\def^{\ifmmode\sp\else\^{}\fi}\catcode`\%=\active\def
\end{pgfscope}%
\begin{pgfscope}%
\pgfpathrectangle{\pgfqpoint{0.608070in}{0.521603in}}{\pgfqpoint{1.550000in}{1.155000in}}%
\pgfusepath{clip}%
\pgfsetrectcap%
\pgfsetroundjoin%
\pgfsetlinewidth{0.803000pt}%
\definecolor{currentstroke}{rgb}{0.690196,0.690196,0.690196}%
\pgfsetstrokecolor{currentstroke}%
\pgfsetdash{}{0pt}%
\pgfpathmoveto{\pgfqpoint{0.608070in}{1.109115in}}%
\pgfpathlineto{\pgfqpoint{2.158070in}{1.109115in}}%
\pgfusepath{stroke}%
\end{pgfscope}%
\begin{pgfscope}%
\pgfsetbuttcap%
\pgfsetroundjoin%
\definecolor{currentfill}{rgb}{0.000000,0.000000,0.000000}%
\pgfsetfillcolor{currentfill}%
\pgfsetlinewidth{0.803000pt}%
\definecolor{currentstroke}{rgb}{0.000000,0.000000,0.000000}%
\pgfsetstrokecolor{currentstroke}%
\pgfsetdash{}{0pt}%
\pgfsys@defobject{currentmarker}{\pgfqpoint{-0.048611in}{0.000000in}}{\pgfqpoint{-0.000000in}{0.000000in}}{%
\pgfpathmoveto{\pgfqpoint{-0.000000in}{0.000000in}}%
\pgfpathlineto{\pgfqpoint{-0.048611in}{0.000000in}}%
\pgfusepath{stroke,fill}%
}%
\begin{pgfscope}%
\pgfsys@transformshift{0.608070in}{1.109115in}%
\pgfsys@useobject{currentmarker}{}%
\end{pgfscope}%
\end{pgfscope}%
\begin{pgfscope}%
\definecolor{textcolor}{rgb}{0.000000,0.000000,0.000000}%
\pgfsetstrokecolor{textcolor}%
\pgfsetfillcolor{textcolor}%
\pgftext[x=0.289968in, y=1.056354in, left, base]{\color{textcolor}{\sffamily\fontsize{10.000000}{12.000000}\selectfont\catcode`\^=\active\def^{\ifmmode\sp\else\^{}\fi}\catcode`\%=\active\def
\end{pgfscope}%
\begin{pgfscope}%
\pgfpathrectangle{\pgfqpoint{0.608070in}{0.521603in}}{\pgfqpoint{1.550000in}{1.155000in}}%
\pgfusepath{clip}%
\pgfsetrectcap%
\pgfsetroundjoin%
\pgfsetlinewidth{0.803000pt}%
\definecolor{currentstroke}{rgb}{0.690196,0.690196,0.690196}%
\pgfsetstrokecolor{currentstroke}%
\pgfsetdash{}{0pt}%
\pgfpathmoveto{\pgfqpoint{0.608070in}{1.460681in}}%
\pgfpathlineto{\pgfqpoint{2.158070in}{1.460681in}}%
\pgfusepath{stroke}%
\end{pgfscope}%
\begin{pgfscope}%
\pgfsetbuttcap%
\pgfsetroundjoin%
\definecolor{currentfill}{rgb}{0.000000,0.000000,0.000000}%
\pgfsetfillcolor{currentfill}%
\pgfsetlinewidth{0.803000pt}%
\definecolor{currentstroke}{rgb}{0.000000,0.000000,0.000000}%
\pgfsetstrokecolor{currentstroke}%
\pgfsetdash{}{0pt}%
\pgfsys@defobject{currentmarker}{\pgfqpoint{-0.048611in}{0.000000in}}{\pgfqpoint{-0.000000in}{0.000000in}}{%
\pgfpathmoveto{\pgfqpoint{-0.000000in}{0.000000in}}%
\pgfpathlineto{\pgfqpoint{-0.048611in}{0.000000in}}%
\pgfusepath{stroke,fill}%
}%
\begin{pgfscope}%
\pgfsys@transformshift{0.608070in}{1.460681in}%
\pgfsys@useobject{currentmarker}{}%
\end{pgfscope}%
\end{pgfscope}%
\begin{pgfscope}%
\definecolor{textcolor}{rgb}{0.000000,0.000000,0.000000}%
\pgfsetstrokecolor{textcolor}%
\pgfsetfillcolor{textcolor}%
\pgftext[x=0.289968in, y=1.407920in, left, base]{\color{textcolor}{\sffamily\fontsize{10.000000}{12.000000}\selectfont\catcode`\^=\active\def^{\ifmmode\sp\else\^{}\fi}\catcode`\%=\active\def
\end{pgfscope}%
\begin{pgfscope}%
\definecolor{textcolor}{rgb}{0.000000,0.000000,0.000000}%
\pgfsetstrokecolor{textcolor}%
\pgfsetfillcolor{textcolor}%
\pgftext[x=0.234413in,y=1.099103in,,bottom,rotate=90.000000]{\color{textcolor}{\sffamily\fontsize{10.000000}{12.000000}\selectfont\catcode`\^=\active\def^{\ifmmode\sp\else\^{}\fi}\catcode`\%=\active\def
\end{pgfscope}%
\begin{pgfscope}%
\pgfpathrectangle{\pgfqpoint{0.608070in}{0.521603in}}{\pgfqpoint{1.550000in}{1.155000in}}%
\pgfusepath{clip}%
\pgfsetrectcap%
\pgfsetroundjoin%
\pgfsetlinewidth{1.505625pt}%
\definecolor{currentstroke}{rgb}{0.121569,0.466667,0.705882}%
\pgfsetstrokecolor{currentstroke}%
\pgfsetdash{}{0pt}%
\pgfpathmoveto{\pgfqpoint{0.678525in}{0.623917in}}%
\pgfpathlineto{\pgfqpoint{0.819434in}{0.888795in}}%
\pgfpathlineto{\pgfqpoint{0.960343in}{0.612964in}}%
\pgfpathlineto{\pgfqpoint{1.101252in}{0.813153in}}%
\pgfpathlineto{\pgfqpoint{1.242161in}{0.804496in}}%
\pgfpathlineto{\pgfqpoint{1.383070in}{0.666786in}}%
\pgfpathlineto{\pgfqpoint{1.523979in}{0.620021in}}%
\pgfpathlineto{\pgfqpoint{1.664888in}{0.574103in}}%
\pgfpathlineto{\pgfqpoint{1.805797in}{0.733157in}}%
\pgfpathlineto{\pgfqpoint{1.946706in}{0.742293in}}%
\pgfpathlineto{\pgfqpoint{2.087615in}{1.624103in}}%
\pgfusepath{stroke}%
\end{pgfscope}%
\begin{pgfscope}%
\pgfsetrectcap%
\pgfsetmiterjoin%
\pgfsetlinewidth{0.803000pt}%
\definecolor{currentstroke}{rgb}{0.000000,0.000000,0.000000}%
\pgfsetstrokecolor{currentstroke}%
\pgfsetdash{}{0pt}%
\pgfpathmoveto{\pgfqpoint{0.608070in}{0.521603in}}%
\pgfpathlineto{\pgfqpoint{0.608070in}{1.676603in}}%
\pgfusepath{stroke}%
\end{pgfscope}%
\begin{pgfscope}%
\pgfsetrectcap%
\pgfsetmiterjoin%
\pgfsetlinewidth{0.803000pt}%
\definecolor{currentstroke}{rgb}{0.000000,0.000000,0.000000}%
\pgfsetstrokecolor{currentstroke}%
\pgfsetdash{}{0pt}%
\pgfpathmoveto{\pgfqpoint{2.158070in}{0.521603in}}%
\pgfpathlineto{\pgfqpoint{2.158070in}{1.676603in}}%
\pgfusepath{stroke}%
\end{pgfscope}%
\begin{pgfscope}%
\pgfsetrectcap%
\pgfsetmiterjoin%
\pgfsetlinewidth{0.803000pt}%
\definecolor{currentstroke}{rgb}{0.000000,0.000000,0.000000}%
\pgfsetstrokecolor{currentstroke}%
\pgfsetdash{}{0pt}%
\pgfpathmoveto{\pgfqpoint{0.608070in}{0.521603in}}%
\pgfpathlineto{\pgfqpoint{2.158070in}{0.521603in}}%
\pgfusepath{stroke}%
\end{pgfscope}%
\begin{pgfscope}%
\pgfsetrectcap%
\pgfsetmiterjoin%
\pgfsetlinewidth{0.803000pt}%
\definecolor{currentstroke}{rgb}{0.000000,0.000000,0.000000}%
\pgfsetstrokecolor{currentstroke}%
\pgfsetdash{}{0pt}%
\pgfpathmoveto{\pgfqpoint{0.608070in}{1.676603in}}%
\pgfpathlineto{\pgfqpoint{2.158070in}{1.676603in}}%
\pgfusepath{stroke}%
\end{pgfscope}%
\end{pgfpicture}%
\makeatother%
\endgroup%

%% file: synth_convergence.pgf
\begingroup%
\makeatletter%
\begin{pgfpicture}%
\pgfpathrectangle{\pgfpointorigin}{\pgfqpoint{2.823419in}{2.931603in}}%
\pgfusepath{use as bounding box, clip}%
\begin{pgfscope}%
\pgfsetbuttcap%
\pgfsetmiterjoin%
\definecolor{currentfill}{rgb}{1.000000,1.000000,1.000000}%
\pgfsetfillcolor{currentfill}%
\pgfsetlinewidth{0.000000pt}%
\definecolor{currentstroke}{rgb}{1.000000,1.000000,1.000000}%
\pgfsetstrokecolor{currentstroke}%
\pgfsetdash{}{0pt}%
\pgfpathmoveto{\pgfqpoint{0.000000in}{0.000000in}}%
\pgfpathlineto{\pgfqpoint{2.823419in}{0.000000in}}%
\pgfpathlineto{\pgfqpoint{2.823419in}{2.931603in}}%
\pgfpathlineto{\pgfqpoint{0.000000in}{2.931603in}}%
\pgfpathlineto{\pgfqpoint{0.000000in}{0.000000in}}%
\pgfpathclose%
\pgfusepath{fill}%
\end{pgfscope}%
\begin{pgfscope}%
\pgfsetbuttcap%
\pgfsetmiterjoin%
\definecolor{currentfill}{rgb}{1.000000,1.000000,1.000000}%
\pgfsetfillcolor{currentfill}%
\pgfsetlinewidth{0.000000pt}%
\definecolor{currentstroke}{rgb}{0.000000,0.000000,0.000000}%
\pgfsetstrokecolor{currentstroke}%
\pgfsetstrokeopacity{0.000000}%
\pgfsetdash{}{0pt}%
\pgfpathmoveto{\pgfqpoint{0.398419in}{0.521603in}}%
\pgfpathlineto{\pgfqpoint{2.723419in}{0.521603in}}%
\pgfpathlineto{\pgfqpoint{2.723419in}{2.831603in}}%
\pgfpathlineto{\pgfqpoint{0.398419in}{2.831603in}}%
\pgfpathlineto{\pgfqpoint{0.398419in}{0.521603in}}%
\pgfpathclose%
\pgfusepath{fill}%
\end{pgfscope}%
\begin{pgfscope}%
\pgfsetbuttcap%
\pgfsetroundjoin%
\definecolor{currentfill}{rgb}{0.000000,0.000000,0.000000}%
\pgfsetfillcolor{currentfill}%
\pgfsetlinewidth{0.803000pt}%
\definecolor{currentstroke}{rgb}{0.000000,0.000000,0.000000}%
\pgfsetstrokecolor{currentstroke}%
\pgfsetdash{}{0pt}%
\pgfsys@defobject{currentmarker}{\pgfqpoint{0.000000in}{-0.048611in}}{\pgfqpoint{0.000000in}{0.000000in}}{%
\pgfpathmoveto{\pgfqpoint{0.000000in}{0.000000in}}%
\pgfpathlineto{\pgfqpoint{0.000000in}{-0.048611in}}%
\pgfusepath{stroke,fill}%
}%
\begin{pgfscope}%
\pgfsys@transformshift{0.504101in}{0.521603in}%
\pgfsys@useobject{currentmarker}{}%
\end{pgfscope}%
\end{pgfscope}%
\begin{pgfscope}%
\definecolor{textcolor}{rgb}{0.000000,0.000000,0.000000}%
\pgfsetstrokecolor{textcolor}%
\pgfsetfillcolor{textcolor}%
\pgftext[x=0.504101in,y=0.424381in,,top]{\color{textcolor}{\sffamily\fontsize{10.000000}{12.000000}\selectfont\catcode`\^=\active\def^{\ifmmode\sp\else\^{}\fi}\catcode`\%=\active\def
\end{pgfscope}%
\begin{pgfscope}%
\pgfsetbuttcap%
\pgfsetroundjoin%
\definecolor{currentfill}{rgb}{0.000000,0.000000,0.000000}%
\pgfsetfillcolor{currentfill}%
\pgfsetlinewidth{0.803000pt}%
\definecolor{currentstroke}{rgb}{0.000000,0.000000,0.000000}%
\pgfsetstrokecolor{currentstroke}%
\pgfsetdash{}{0pt}%
\pgfsys@defobject{currentmarker}{\pgfqpoint{0.000000in}{-0.048611in}}{\pgfqpoint{0.000000in}{0.000000in}}{%
\pgfpathmoveto{\pgfqpoint{0.000000in}{0.000000in}}%
\pgfpathlineto{\pgfqpoint{0.000000in}{-0.048611in}}%
\pgfusepath{stroke,fill}%
}%
\begin{pgfscope}%
\pgfsys@transformshift{0.973798in}{0.521603in}%
\pgfsys@useobject{currentmarker}{}%
\end{pgfscope}%
\end{pgfscope}%
\begin{pgfscope}%
\definecolor{textcolor}{rgb}{0.000000,0.000000,0.000000}%
\pgfsetstrokecolor{textcolor}%
\pgfsetfillcolor{textcolor}%
\pgftext[x=0.973798in,y=0.424381in,,top]{\color{textcolor}{\sffamily\fontsize{10.000000}{12.000000}\selectfont\catcode`\^=\active\def^{\ifmmode\sp\else\^{}\fi}\catcode`\%=\active\def
\end{pgfscope}%
\begin{pgfscope}%
\pgfsetbuttcap%
\pgfsetroundjoin%
\definecolor{currentfill}{rgb}{0.000000,0.000000,0.000000}%
\pgfsetfillcolor{currentfill}%
\pgfsetlinewidth{0.803000pt}%
\definecolor{currentstroke}{rgb}{0.000000,0.000000,0.000000}%
\pgfsetstrokecolor{currentstroke}%
\pgfsetdash{}{0pt}%
\pgfsys@defobject{currentmarker}{\pgfqpoint{0.000000in}{-0.048611in}}{\pgfqpoint{0.000000in}{0.000000in}}{%
\pgfpathmoveto{\pgfqpoint{0.000000in}{0.000000in}}%
\pgfpathlineto{\pgfqpoint{0.000000in}{-0.048611in}}%
\pgfusepath{stroke,fill}%
}%
\begin{pgfscope}%
\pgfsys@transformshift{1.443495in}{0.521603in}%
\pgfsys@useobject{currentmarker}{}%
\end{pgfscope}%
\end{pgfscope}%
\begin{pgfscope}%
\definecolor{textcolor}{rgb}{0.000000,0.000000,0.000000}%
\pgfsetstrokecolor{textcolor}%
\pgfsetfillcolor{textcolor}%
\pgftext[x=1.443495in,y=0.424381in,,top]{\color{textcolor}{\sffamily\fontsize{10.000000}{12.000000}\selectfont\catcode`\^=\active\def^{\ifmmode\sp\else\^{}\fi}\catcode`\%=\active\def
\end{pgfscope}%
\begin{pgfscope}%
\pgfsetbuttcap%
\pgfsetroundjoin%
\definecolor{currentfill}{rgb}{0.000000,0.000000,0.000000}%
\pgfsetfillcolor{currentfill}%
\pgfsetlinewidth{0.803000pt}%
\definecolor{currentstroke}{rgb}{0.000000,0.000000,0.000000}%
\pgfsetstrokecolor{currentstroke}%
\pgfsetdash{}{0pt}%
\pgfsys@defobject{currentmarker}{\pgfqpoint{0.000000in}{-0.048611in}}{\pgfqpoint{0.000000in}{0.000000in}}{%
\pgfpathmoveto{\pgfqpoint{0.000000in}{0.000000in}}%
\pgfpathlineto{\pgfqpoint{0.000000in}{-0.048611in}}%
\pgfusepath{stroke,fill}%
}%
\begin{pgfscope}%
\pgfsys@transformshift{1.913192in}{0.521603in}%
\pgfsys@useobject{currentmarker}{}%
\end{pgfscope}%
\end{pgfscope}%
\begin{pgfscope}%
\definecolor{textcolor}{rgb}{0.000000,0.000000,0.000000}%
\pgfsetstrokecolor{textcolor}%
\pgfsetfillcolor{textcolor}%
\pgftext[x=1.913192in,y=0.424381in,,top]{\color{textcolor}{\sffamily\fontsize{10.000000}{12.000000}\selectfont\catcode`\^=\active\def^{\ifmmode\sp\else\^{}\fi}\catcode`\%=\active\def
\end{pgfscope}%
\begin{pgfscope}%
\pgfsetbuttcap%
\pgfsetroundjoin%
\definecolor{currentfill}{rgb}{0.000000,0.000000,0.000000}%
\pgfsetfillcolor{currentfill}%
\pgfsetlinewidth{0.803000pt}%
\definecolor{currentstroke}{rgb}{0.000000,0.000000,0.000000}%
\pgfsetstrokecolor{currentstroke}%
\pgfsetdash{}{0pt}%
\pgfsys@defobject{currentmarker}{\pgfqpoint{0.000000in}{-0.048611in}}{\pgfqpoint{0.000000in}{0.000000in}}{%
\pgfpathmoveto{\pgfqpoint{0.000000in}{0.000000in}}%
\pgfpathlineto{\pgfqpoint{0.000000in}{-0.048611in}}%
\pgfusepath{stroke,fill}%
}%
\begin{pgfscope}%
\pgfsys@transformshift{2.382889in}{0.521603in}%
\pgfsys@useobject{currentmarker}{}%
\end{pgfscope}%
\end{pgfscope}%
\begin{pgfscope}%
\definecolor{textcolor}{rgb}{0.000000,0.000000,0.000000}%
\pgfsetstrokecolor{textcolor}%
\pgfsetfillcolor{textcolor}%
\pgftext[x=2.382889in,y=0.424381in,,top]{\color{textcolor}{\sffamily\fontsize{10.000000}{12.000000}\selectfont\catcode`\^=\active\def^{\ifmmode\sp\else\^{}\fi}\catcode`\%=\active\def
\end{pgfscope}%
\begin{pgfscope}%
\definecolor{textcolor}{rgb}{0.000000,0.000000,0.000000}%
\pgfsetstrokecolor{textcolor}%
\pgfsetfillcolor{textcolor}%
\pgftext[x=1.560919in,y=0.234413in,,top]{\color{textcolor}{\sffamily\fontsize{10.000000}{12.000000}\selectfont\catcode`\^=\active\def^{\ifmmode\sp\else\^{}\fi}\catcode`\%=\active\def
\end{pgfscope}%
\begin{pgfscope}%
\pgfsetbuttcap%
\pgfsetroundjoin%
\definecolor{currentfill}{rgb}{0.000000,0.000000,0.000000}%
\pgfsetfillcolor{currentfill}%
\pgfsetlinewidth{0.803000pt}%
\definecolor{currentstroke}{rgb}{0.000000,0.000000,0.000000}%
\pgfsetstrokecolor{currentstroke}%
\pgfsetdash{}{0pt}%
\pgfsys@defobject{currentmarker}{\pgfqpoint{-0.048611in}{0.000000in}}{\pgfqpoint{-0.000000in}{0.000000in}}{%
\pgfpathmoveto{\pgfqpoint{-0.000000in}{0.000000in}}%
\pgfpathlineto{\pgfqpoint{-0.048611in}{0.000000in}}%
\pgfusepath{stroke,fill}%
}%
\begin{pgfscope}%
\pgfsys@transformshift{0.398419in}{0.566944in}%
\pgfsys@useobject{currentmarker}{}%
\end{pgfscope}%
\end{pgfscope}%
\begin{pgfscope}%
\definecolor{textcolor}{rgb}{0.000000,0.000000,0.000000}%
\pgfsetstrokecolor{textcolor}%
\pgfsetfillcolor{textcolor}%
\pgftext[x=0.100000in, y=0.514183in, left, base]{\color{textcolor}{\sffamily\fontsize{10.000000}{12.000000}\selectfont\catcode`\^=\active\def^{\ifmmode\sp\else\^{}\fi}\catcode`\%=\active\def
\end{pgfscope}%
\begin{pgfscope}%
\pgfsetbuttcap%
\pgfsetroundjoin%
\definecolor{currentfill}{rgb}{0.000000,0.000000,0.000000}%
\pgfsetfillcolor{currentfill}%
\pgfsetlinewidth{0.803000pt}%
\definecolor{currentstroke}{rgb}{0.000000,0.000000,0.000000}%
\pgfsetstrokecolor{currentstroke}%
\pgfsetdash{}{0pt}%
\pgfsys@defobject{currentmarker}{\pgfqpoint{-0.048611in}{0.000000in}}{\pgfqpoint{-0.000000in}{0.000000in}}{%
\pgfpathmoveto{\pgfqpoint{-0.000000in}{0.000000in}}%
\pgfpathlineto{\pgfqpoint{-0.048611in}{0.000000in}}%
\pgfusepath{stroke,fill}%
}%
\begin{pgfscope}%
\pgfsys@transformshift{0.398419in}{1.440767in}%
\pgfsys@useobject{currentmarker}{}%
\end{pgfscope}%
\end{pgfscope}%
\begin{pgfscope}%
\definecolor{textcolor}{rgb}{0.000000,0.000000,0.000000}%
\pgfsetstrokecolor{textcolor}%
\pgfsetfillcolor{textcolor}%
\pgftext[x=0.100000in, y=1.388005in, left, base]{\color{textcolor}{\sffamily\fontsize{10.000000}{12.000000}\selectfont\catcode`\^=\active\def^{\ifmmode\sp\else\^{}\fi}\catcode`\%=\active\def
\end{pgfscope}%
\begin{pgfscope}%
\pgfsetbuttcap%
\pgfsetroundjoin%
\definecolor{currentfill}{rgb}{0.000000,0.000000,0.000000}%
\pgfsetfillcolor{currentfill}%
\pgfsetlinewidth{0.803000pt}%
\definecolor{currentstroke}{rgb}{0.000000,0.000000,0.000000}%
\pgfsetstrokecolor{currentstroke}%
\pgfsetdash{}{0pt}%
\pgfsys@defobject{currentmarker}{\pgfqpoint{-0.048611in}{0.000000in}}{\pgfqpoint{-0.000000in}{0.000000in}}{%
\pgfpathmoveto{\pgfqpoint{-0.000000in}{0.000000in}}%
\pgfpathlineto{\pgfqpoint{-0.048611in}{0.000000in}}%
\pgfusepath{stroke,fill}%
}%
\begin{pgfscope}%
\pgfsys@transformshift{0.398419in}{2.314589in}%
\pgfsys@useobject{currentmarker}{}%
\end{pgfscope}%
\end{pgfscope}%
\begin{pgfscope}%
\definecolor{textcolor}{rgb}{0.000000,0.000000,0.000000}%
\pgfsetstrokecolor{textcolor}%
\pgfsetfillcolor{textcolor}%
\pgftext[x=0.100000in, y=2.261827in, left, base]{\color{textcolor}{\sffamily\fontsize{10.000000}{12.000000}\selectfont\catcode`\^=\active\def^{\ifmmode\sp\else\^{}\fi}\catcode`\%=\active\def
\end{pgfscope}%
\begin{pgfscope}%
\pgfsetbuttcap%
\pgfsetroundjoin%
\definecolor{currentfill}{rgb}{0.000000,0.000000,0.000000}%
\pgfsetfillcolor{currentfill}%
\pgfsetlinewidth{0.602250pt}%
\definecolor{currentstroke}{rgb}{0.000000,0.000000,0.000000}%
\pgfsetstrokecolor{currentstroke}%
\pgfsetdash{}{0pt}%
\pgfsys@defobject{currentmarker}{\pgfqpoint{-0.027778in}{0.000000in}}{\pgfqpoint{-0.000000in}{0.000000in}}{%
\pgfpathmoveto{\pgfqpoint{-0.000000in}{0.000000in}}%
\pgfpathlineto{\pgfqpoint{-0.027778in}{0.000000in}}%
\pgfusepath{stroke,fill}%
}%
\begin{pgfscope}%
\pgfsys@transformshift{0.398419in}{0.526960in}%
\pgfsys@useobject{currentmarker}{}%
\end{pgfscope}%
\end{pgfscope}%
\begin{pgfscope}%
\pgfsetbuttcap%
\pgfsetroundjoin%
\definecolor{currentfill}{rgb}{0.000000,0.000000,0.000000}%
\pgfsetfillcolor{currentfill}%
\pgfsetlinewidth{0.602250pt}%
\definecolor{currentstroke}{rgb}{0.000000,0.000000,0.000000}%
\pgfsetstrokecolor{currentstroke}%
\pgfsetdash{}{0pt}%
\pgfsys@defobject{currentmarker}{\pgfqpoint{-0.027778in}{0.000000in}}{\pgfqpoint{-0.000000in}{0.000000in}}{%
\pgfpathmoveto{\pgfqpoint{-0.000000in}{0.000000in}}%
\pgfpathlineto{\pgfqpoint{-0.027778in}{0.000000in}}%
\pgfusepath{stroke,fill}%
}%
\begin{pgfscope}%
\pgfsys@transformshift{0.398419in}{0.829991in}%
\pgfsys@useobject{currentmarker}{}%
\end{pgfscope}%
\end{pgfscope}%
\begin{pgfscope}%
\pgfsetbuttcap%
\pgfsetroundjoin%
\definecolor{currentfill}{rgb}{0.000000,0.000000,0.000000}%
\pgfsetfillcolor{currentfill}%
\pgfsetlinewidth{0.602250pt}%
\definecolor{currentstroke}{rgb}{0.000000,0.000000,0.000000}%
\pgfsetstrokecolor{currentstroke}%
\pgfsetdash{}{0pt}%
\pgfsys@defobject{currentmarker}{\pgfqpoint{-0.027778in}{0.000000in}}{\pgfqpoint{-0.000000in}{0.000000in}}{%
\pgfpathmoveto{\pgfqpoint{-0.000000in}{0.000000in}}%
\pgfpathlineto{\pgfqpoint{-0.027778in}{0.000000in}}%
\pgfusepath{stroke,fill}%
}%
\begin{pgfscope}%
\pgfsys@transformshift{0.398419in}{0.983864in}%
\pgfsys@useobject{currentmarker}{}%
\end{pgfscope}%
\end{pgfscope}%
\begin{pgfscope}%
\pgfsetbuttcap%
\pgfsetroundjoin%
\definecolor{currentfill}{rgb}{0.000000,0.000000,0.000000}%
\pgfsetfillcolor{currentfill}%
\pgfsetlinewidth{0.602250pt}%
\definecolor{currentstroke}{rgb}{0.000000,0.000000,0.000000}%
\pgfsetstrokecolor{currentstroke}%
\pgfsetdash{}{0pt}%
\pgfsys@defobject{currentmarker}{\pgfqpoint{-0.027778in}{0.000000in}}{\pgfqpoint{-0.000000in}{0.000000in}}{%
\pgfpathmoveto{\pgfqpoint{-0.000000in}{0.000000in}}%
\pgfpathlineto{\pgfqpoint{-0.027778in}{0.000000in}}%
\pgfusepath{stroke,fill}%
}%
\begin{pgfscope}%
\pgfsys@transformshift{0.398419in}{1.093038in}%
\pgfsys@useobject{currentmarker}{}%
\end{pgfscope}%
\end{pgfscope}%
\begin{pgfscope}%
\pgfsetbuttcap%
\pgfsetroundjoin%
\definecolor{currentfill}{rgb}{0.000000,0.000000,0.000000}%
\pgfsetfillcolor{currentfill}%
\pgfsetlinewidth{0.602250pt}%
\definecolor{currentstroke}{rgb}{0.000000,0.000000,0.000000}%
\pgfsetstrokecolor{currentstroke}%
\pgfsetdash{}{0pt}%
\pgfsys@defobject{currentmarker}{\pgfqpoint{-0.027778in}{0.000000in}}{\pgfqpoint{-0.000000in}{0.000000in}}{%
\pgfpathmoveto{\pgfqpoint{-0.000000in}{0.000000in}}%
\pgfpathlineto{\pgfqpoint{-0.027778in}{0.000000in}}%
\pgfusepath{stroke,fill}%
}%
\begin{pgfscope}%
\pgfsys@transformshift{0.398419in}{1.177720in}%
\pgfsys@useobject{currentmarker}{}%
\end{pgfscope}%
\end{pgfscope}%
\begin{pgfscope}%
\pgfsetbuttcap%
\pgfsetroundjoin%
\definecolor{currentfill}{rgb}{0.000000,0.000000,0.000000}%
\pgfsetfillcolor{currentfill}%
\pgfsetlinewidth{0.602250pt}%
\definecolor{currentstroke}{rgb}{0.000000,0.000000,0.000000}%
\pgfsetstrokecolor{currentstroke}%
\pgfsetdash{}{0pt}%
\pgfsys@defobject{currentmarker}{\pgfqpoint{-0.027778in}{0.000000in}}{\pgfqpoint{-0.000000in}{0.000000in}}{%
\pgfpathmoveto{\pgfqpoint{-0.000000in}{0.000000in}}%
\pgfpathlineto{\pgfqpoint{-0.027778in}{0.000000in}}%
\pgfusepath{stroke,fill}%
}%
\begin{pgfscope}%
\pgfsys@transformshift{0.398419in}{1.246910in}%
\pgfsys@useobject{currentmarker}{}%
\end{pgfscope}%
\end{pgfscope}%
\begin{pgfscope}%
\pgfsetbuttcap%
\pgfsetroundjoin%
\definecolor{currentfill}{rgb}{0.000000,0.000000,0.000000}%
\pgfsetfillcolor{currentfill}%
\pgfsetlinewidth{0.602250pt}%
\definecolor{currentstroke}{rgb}{0.000000,0.000000,0.000000}%
\pgfsetstrokecolor{currentstroke}%
\pgfsetdash{}{0pt}%
\pgfsys@defobject{currentmarker}{\pgfqpoint{-0.027778in}{0.000000in}}{\pgfqpoint{-0.000000in}{0.000000in}}{%
\pgfpathmoveto{\pgfqpoint{-0.000000in}{0.000000in}}%
\pgfpathlineto{\pgfqpoint{-0.027778in}{0.000000in}}%
\pgfusepath{stroke,fill}%
}%
\begin{pgfscope}%
\pgfsys@transformshift{0.398419in}{1.305410in}%
\pgfsys@useobject{currentmarker}{}%
\end{pgfscope}%
\end{pgfscope}%
\begin{pgfscope}%
\pgfsetbuttcap%
\pgfsetroundjoin%
\definecolor{currentfill}{rgb}{0.000000,0.000000,0.000000}%
\pgfsetfillcolor{currentfill}%
\pgfsetlinewidth{0.602250pt}%
\definecolor{currentstroke}{rgb}{0.000000,0.000000,0.000000}%
\pgfsetstrokecolor{currentstroke}%
\pgfsetdash{}{0pt}%
\pgfsys@defobject{currentmarker}{\pgfqpoint{-0.027778in}{0.000000in}}{\pgfqpoint{-0.000000in}{0.000000in}}{%
\pgfpathmoveto{\pgfqpoint{-0.000000in}{0.000000in}}%
\pgfpathlineto{\pgfqpoint{-0.027778in}{0.000000in}}%
\pgfusepath{stroke,fill}%
}%
\begin{pgfscope}%
\pgfsys@transformshift{0.398419in}{1.356084in}%
\pgfsys@useobject{currentmarker}{}%
\end{pgfscope}%
\end{pgfscope}%
\begin{pgfscope}%
\pgfsetbuttcap%
\pgfsetroundjoin%
\definecolor{currentfill}{rgb}{0.000000,0.000000,0.000000}%
\pgfsetfillcolor{currentfill}%
\pgfsetlinewidth{0.602250pt}%
\definecolor{currentstroke}{rgb}{0.000000,0.000000,0.000000}%
\pgfsetstrokecolor{currentstroke}%
\pgfsetdash{}{0pt}%
\pgfsys@defobject{currentmarker}{\pgfqpoint{-0.027778in}{0.000000in}}{\pgfqpoint{-0.000000in}{0.000000in}}{%
\pgfpathmoveto{\pgfqpoint{-0.000000in}{0.000000in}}%
\pgfpathlineto{\pgfqpoint{-0.027778in}{0.000000in}}%
\pgfusepath{stroke,fill}%
}%
\begin{pgfscope}%
\pgfsys@transformshift{0.398419in}{1.400783in}%
\pgfsys@useobject{currentmarker}{}%
\end{pgfscope}%
\end{pgfscope}%
\begin{pgfscope}%
\pgfsetbuttcap%
\pgfsetroundjoin%
\definecolor{currentfill}{rgb}{0.000000,0.000000,0.000000}%
\pgfsetfillcolor{currentfill}%
\pgfsetlinewidth{0.602250pt}%
\definecolor{currentstroke}{rgb}{0.000000,0.000000,0.000000}%
\pgfsetstrokecolor{currentstroke}%
\pgfsetdash{}{0pt}%
\pgfsys@defobject{currentmarker}{\pgfqpoint{-0.027778in}{0.000000in}}{\pgfqpoint{-0.000000in}{0.000000in}}{%
\pgfpathmoveto{\pgfqpoint{-0.000000in}{0.000000in}}%
\pgfpathlineto{\pgfqpoint{-0.027778in}{0.000000in}}%
\pgfusepath{stroke,fill}%
}%
\begin{pgfscope}%
\pgfsys@transformshift{0.398419in}{1.703813in}%
\pgfsys@useobject{currentmarker}{}%
\end{pgfscope}%
\end{pgfscope}%
\begin{pgfscope}%
\pgfsetbuttcap%
\pgfsetroundjoin%
\definecolor{currentfill}{rgb}{0.000000,0.000000,0.000000}%
\pgfsetfillcolor{currentfill}%
\pgfsetlinewidth{0.602250pt}%
\definecolor{currentstroke}{rgb}{0.000000,0.000000,0.000000}%
\pgfsetstrokecolor{currentstroke}%
\pgfsetdash{}{0pt}%
\pgfsys@defobject{currentmarker}{\pgfqpoint{-0.027778in}{0.000000in}}{\pgfqpoint{-0.000000in}{0.000000in}}{%
\pgfpathmoveto{\pgfqpoint{-0.000000in}{0.000000in}}%
\pgfpathlineto{\pgfqpoint{-0.027778in}{0.000000in}}%
\pgfusepath{stroke,fill}%
}%
\begin{pgfscope}%
\pgfsys@transformshift{0.398419in}{1.857686in}%
\pgfsys@useobject{currentmarker}{}%
\end{pgfscope}%
\end{pgfscope}%
\begin{pgfscope}%
\pgfsetbuttcap%
\pgfsetroundjoin%
\definecolor{currentfill}{rgb}{0.000000,0.000000,0.000000}%
\pgfsetfillcolor{currentfill}%
\pgfsetlinewidth{0.602250pt}%
\definecolor{currentstroke}{rgb}{0.000000,0.000000,0.000000}%
\pgfsetstrokecolor{currentstroke}%
\pgfsetdash{}{0pt}%
\pgfsys@defobject{currentmarker}{\pgfqpoint{-0.027778in}{0.000000in}}{\pgfqpoint{-0.000000in}{0.000000in}}{%
\pgfpathmoveto{\pgfqpoint{-0.000000in}{0.000000in}}%
\pgfpathlineto{\pgfqpoint{-0.027778in}{0.000000in}}%
\pgfusepath{stroke,fill}%
}%
\begin{pgfscope}%
\pgfsys@transformshift{0.398419in}{1.966860in}%
\pgfsys@useobject{currentmarker}{}%
\end{pgfscope}%
\end{pgfscope}%
\begin{pgfscope}%
\pgfsetbuttcap%
\pgfsetroundjoin%
\definecolor{currentfill}{rgb}{0.000000,0.000000,0.000000}%
\pgfsetfillcolor{currentfill}%
\pgfsetlinewidth{0.602250pt}%
\definecolor{currentstroke}{rgb}{0.000000,0.000000,0.000000}%
\pgfsetstrokecolor{currentstroke}%
\pgfsetdash{}{0pt}%
\pgfsys@defobject{currentmarker}{\pgfqpoint{-0.027778in}{0.000000in}}{\pgfqpoint{-0.000000in}{0.000000in}}{%
\pgfpathmoveto{\pgfqpoint{-0.000000in}{0.000000in}}%
\pgfpathlineto{\pgfqpoint{-0.027778in}{0.000000in}}%
\pgfusepath{stroke,fill}%
}%
\begin{pgfscope}%
\pgfsys@transformshift{0.398419in}{2.051542in}%
\pgfsys@useobject{currentmarker}{}%
\end{pgfscope}%
\end{pgfscope}%
\begin{pgfscope}%
\pgfsetbuttcap%
\pgfsetroundjoin%
\definecolor{currentfill}{rgb}{0.000000,0.000000,0.000000}%
\pgfsetfillcolor{currentfill}%
\pgfsetlinewidth{0.602250pt}%
\definecolor{currentstroke}{rgb}{0.000000,0.000000,0.000000}%
\pgfsetstrokecolor{currentstroke}%
\pgfsetdash{}{0pt}%
\pgfsys@defobject{currentmarker}{\pgfqpoint{-0.027778in}{0.000000in}}{\pgfqpoint{-0.000000in}{0.000000in}}{%
\pgfpathmoveto{\pgfqpoint{-0.000000in}{0.000000in}}%
\pgfpathlineto{\pgfqpoint{-0.027778in}{0.000000in}}%
\pgfusepath{stroke,fill}%
}%
\begin{pgfscope}%
\pgfsys@transformshift{0.398419in}{2.120732in}%
\pgfsys@useobject{currentmarker}{}%
\end{pgfscope}%
\end{pgfscope}%
\begin{pgfscope}%
\pgfsetbuttcap%
\pgfsetroundjoin%
\definecolor{currentfill}{rgb}{0.000000,0.000000,0.000000}%
\pgfsetfillcolor{currentfill}%
\pgfsetlinewidth{0.602250pt}%
\definecolor{currentstroke}{rgb}{0.000000,0.000000,0.000000}%
\pgfsetstrokecolor{currentstroke}%
\pgfsetdash{}{0pt}%
\pgfsys@defobject{currentmarker}{\pgfqpoint{-0.027778in}{0.000000in}}{\pgfqpoint{-0.000000in}{0.000000in}}{%
\pgfpathmoveto{\pgfqpoint{-0.000000in}{0.000000in}}%
\pgfpathlineto{\pgfqpoint{-0.027778in}{0.000000in}}%
\pgfusepath{stroke,fill}%
}%
\begin{pgfscope}%
\pgfsys@transformshift{0.398419in}{2.179232in}%
\pgfsys@useobject{currentmarker}{}%
\end{pgfscope}%
\end{pgfscope}%
\begin{pgfscope}%
\pgfsetbuttcap%
\pgfsetroundjoin%
\definecolor{currentfill}{rgb}{0.000000,0.000000,0.000000}%
\pgfsetfillcolor{currentfill}%
\pgfsetlinewidth{0.602250pt}%
\definecolor{currentstroke}{rgb}{0.000000,0.000000,0.000000}%
\pgfsetstrokecolor{currentstroke}%
\pgfsetdash{}{0pt}%
\pgfsys@defobject{currentmarker}{\pgfqpoint{-0.027778in}{0.000000in}}{\pgfqpoint{-0.000000in}{0.000000in}}{%
\pgfpathmoveto{\pgfqpoint{-0.000000in}{0.000000in}}%
\pgfpathlineto{\pgfqpoint{-0.027778in}{0.000000in}}%
\pgfusepath{stroke,fill}%
}%
\begin{pgfscope}%
\pgfsys@transformshift{0.398419in}{2.229907in}%
\pgfsys@useobject{currentmarker}{}%
\end{pgfscope}%
\end{pgfscope}%
\begin{pgfscope}%
\pgfsetbuttcap%
\pgfsetroundjoin%
\definecolor{currentfill}{rgb}{0.000000,0.000000,0.000000}%
\pgfsetfillcolor{currentfill}%
\pgfsetlinewidth{0.602250pt}%
\definecolor{currentstroke}{rgb}{0.000000,0.000000,0.000000}%
\pgfsetstrokecolor{currentstroke}%
\pgfsetdash{}{0pt}%
\pgfsys@defobject{currentmarker}{\pgfqpoint{-0.027778in}{0.000000in}}{\pgfqpoint{-0.000000in}{0.000000in}}{%
\pgfpathmoveto{\pgfqpoint{-0.000000in}{0.000000in}}%
\pgfpathlineto{\pgfqpoint{-0.027778in}{0.000000in}}%
\pgfusepath{stroke,fill}%
}%
\begin{pgfscope}%
\pgfsys@transformshift{0.398419in}{2.274605in}%
\pgfsys@useobject{currentmarker}{}%
\end{pgfscope}%
\end{pgfscope}%
\begin{pgfscope}%
\pgfsetbuttcap%
\pgfsetroundjoin%
\definecolor{currentfill}{rgb}{0.000000,0.000000,0.000000}%
\pgfsetfillcolor{currentfill}%
\pgfsetlinewidth{0.602250pt}%
\definecolor{currentstroke}{rgb}{0.000000,0.000000,0.000000}%
\pgfsetstrokecolor{currentstroke}%
\pgfsetdash{}{0pt}%
\pgfsys@defobject{currentmarker}{\pgfqpoint{-0.027778in}{0.000000in}}{\pgfqpoint{-0.000000in}{0.000000in}}{%
\pgfpathmoveto{\pgfqpoint{-0.000000in}{0.000000in}}%
\pgfpathlineto{\pgfqpoint{-0.027778in}{0.000000in}}%
\pgfusepath{stroke,fill}%
}%
\begin{pgfscope}%
\pgfsys@transformshift{0.398419in}{2.577636in}%
\pgfsys@useobject{currentmarker}{}%
\end{pgfscope}%
\end{pgfscope}%
\begin{pgfscope}%
\pgfsetbuttcap%
\pgfsetroundjoin%
\definecolor{currentfill}{rgb}{0.000000,0.000000,0.000000}%
\pgfsetfillcolor{currentfill}%
\pgfsetlinewidth{0.602250pt}%
\definecolor{currentstroke}{rgb}{0.000000,0.000000,0.000000}%
\pgfsetstrokecolor{currentstroke}%
\pgfsetdash{}{0pt}%
\pgfsys@defobject{currentmarker}{\pgfqpoint{-0.027778in}{0.000000in}}{\pgfqpoint{-0.000000in}{0.000000in}}{%
\pgfpathmoveto{\pgfqpoint{-0.000000in}{0.000000in}}%
\pgfpathlineto{\pgfqpoint{-0.027778in}{0.000000in}}%
\pgfusepath{stroke,fill}%
}%
\begin{pgfscope}%
\pgfsys@transformshift{0.398419in}{2.731508in}%
\pgfsys@useobject{currentmarker}{}%
\end{pgfscope}%
\end{pgfscope}%
\begin{pgfscope}%
\pgfpathrectangle{\pgfqpoint{0.398419in}{0.521603in}}{\pgfqpoint{2.325000in}{2.310000in}}%
\pgfusepath{clip}%
\pgfsetrectcap%
\pgfsetroundjoin%
\pgfsetlinewidth{1.505625pt}%
\definecolor{currentstroke}{rgb}{0.121569,0.466667,0.705882}%
\pgfsetstrokecolor{currentstroke}%
\pgfsetdash{}{0pt}%
\pgfpathmoveto{\pgfqpoint{0.504101in}{0.626603in}}%
\pgfpathlineto{\pgfqpoint{0.738949in}{0.697055in}}%
\pgfpathlineto{\pgfqpoint{0.973798in}{0.797894in}}%
\pgfpathlineto{\pgfqpoint{1.208646in}{0.874339in}}%
\pgfpathlineto{\pgfqpoint{1.443495in}{0.928795in}}%
\pgfpathlineto{\pgfqpoint{1.678343in}{0.969355in}}%
\pgfpathlineto{\pgfqpoint{1.913192in}{1.000710in}}%
\pgfpathlineto{\pgfqpoint{2.148040in}{1.025527in}}%
\pgfpathlineto{\pgfqpoint{2.382889in}{1.045589in}}%
\pgfpathlineto{\pgfqpoint{2.617737in}{1.062147in}}%
\pgfusepath{stroke}%
\end{pgfscope}%
\begin{pgfscope}%
\pgfpathrectangle{\pgfqpoint{0.398419in}{0.521603in}}{\pgfqpoint{2.325000in}{2.310000in}}%
\pgfusepath{clip}%
\pgfsetrectcap%
\pgfsetroundjoin%
\pgfsetlinewidth{1.505625pt}%
\definecolor{currentstroke}{rgb}{1.000000,0.498039,0.054902}%
\pgfsetstrokecolor{currentstroke}%
\pgfsetdash{}{0pt}%
\pgfpathmoveto{\pgfqpoint{0.504101in}{2.726603in}}%
\pgfpathlineto{\pgfqpoint{0.738949in}{2.589728in}}%
\pgfpathlineto{\pgfqpoint{0.973798in}{2.465285in}}%
\pgfpathlineto{\pgfqpoint{1.208646in}{2.372624in}}%
\pgfpathlineto{\pgfqpoint{1.443495in}{2.300013in}}%
\pgfpathlineto{\pgfqpoint{1.678343in}{2.240700in}}%
\pgfpathlineto{\pgfqpoint{1.913192in}{2.191562in}}%
\pgfpathlineto{\pgfqpoint{2.148040in}{2.150920in}}%
\pgfpathlineto{\pgfqpoint{2.382889in}{2.117204in}}%
\pgfpathlineto{\pgfqpoint{2.617737in}{2.089050in}}%
\pgfusepath{stroke}%
\end{pgfscope}%
\begin{pgfscope}%
\pgfpathrectangle{\pgfqpoint{0.398419in}{0.521603in}}{\pgfqpoint{2.325000in}{2.310000in}}%
\pgfusepath{clip}%
\pgfsetrectcap%
\pgfsetroundjoin%
\pgfsetlinewidth{1.505625pt}%
\definecolor{currentstroke}{rgb}{0.172549,0.627451,0.172549}%
\pgfsetstrokecolor{currentstroke}%
\pgfsetdash{}{0pt}%
\pgfpathmoveto{\pgfqpoint{0.504101in}{1.030395in}}%
\pgfpathlineto{\pgfqpoint{0.738949in}{1.103222in}}%
\pgfpathlineto{\pgfqpoint{0.973798in}{1.147119in}}%
\pgfpathlineto{\pgfqpoint{1.208646in}{1.181335in}}%
\pgfpathlineto{\pgfqpoint{1.443495in}{1.209256in}}%
\pgfpathlineto{\pgfqpoint{1.678343in}{1.232582in}}%
\pgfpathlineto{\pgfqpoint{1.913192in}{1.252511in}}%
\pgfpathlineto{\pgfqpoint{2.148040in}{1.269931in}}%
\pgfpathlineto{\pgfqpoint{2.382889in}{1.285470in}}%
\pgfpathlineto{\pgfqpoint{2.617737in}{1.299557in}}%
\pgfusepath{stroke}%
\end{pgfscope}%
\begin{pgfscope}%
\pgfsetrectcap%
\pgfsetmiterjoin%
\pgfsetlinewidth{0.803000pt}%
\definecolor{currentstroke}{rgb}{0.000000,0.000000,0.000000}%
\pgfsetstrokecolor{currentstroke}%
\pgfsetdash{}{0pt}%
\pgfpathmoveto{\pgfqpoint{0.398419in}{0.521603in}}%
\pgfpathlineto{\pgfqpoint{0.398419in}{2.831603in}}%
\pgfusepath{stroke}%
\end{pgfscope}%
\begin{pgfscope}%
\pgfsetrectcap%
\pgfsetmiterjoin%
\pgfsetlinewidth{0.803000pt}%
\definecolor{currentstroke}{rgb}{0.000000,0.000000,0.000000}%
\pgfsetstrokecolor{currentstroke}%
\pgfsetdash{}{0pt}%
\pgfpathmoveto{\pgfqpoint{2.723419in}{0.521603in}}%
\pgfpathlineto{\pgfqpoint{2.723419in}{2.831603in}}%
\pgfusepath{stroke}%
\end{pgfscope}%
\begin{pgfscope}%
\pgfsetrectcap%
\pgfsetmiterjoin%
\pgfsetlinewidth{0.803000pt}%
\definecolor{currentstroke}{rgb}{0.000000,0.000000,0.000000}%
\pgfsetstrokecolor{currentstroke}%
\pgfsetdash{}{0pt}%
\pgfpathmoveto{\pgfqpoint{0.398419in}{0.521603in}}%
\pgfpathlineto{\pgfqpoint{2.723419in}{0.521603in}}%
\pgfusepath{stroke}%
\end{pgfscope}%
\begin{pgfscope}%
\pgfsetrectcap%
\pgfsetmiterjoin%
\pgfsetlinewidth{0.803000pt}%
\definecolor{currentstroke}{rgb}{0.000000,0.000000,0.000000}%
\pgfsetstrokecolor{currentstroke}%
\pgfsetdash{}{0pt}%
\pgfpathmoveto{\pgfqpoint{0.398419in}{2.831603in}}%
\pgfpathlineto{\pgfqpoint{2.723419in}{2.831603in}}%
\pgfusepath{stroke}%
\end{pgfscope}%
\begin{pgfscope}%
\pgfsetbuttcap%
\pgfsetmiterjoin%
\definecolor{currentfill}{rgb}{1.000000,1.000000,1.000000}%
\pgfsetfillcolor{currentfill}%
\pgfsetfillopacity{0.800000}%
\pgfsetlinewidth{1.003750pt}%
\definecolor{currentstroke}{rgb}{0.800000,0.800000,0.800000}%
\pgfsetstrokecolor{currentstroke}%
\pgfsetstrokeopacity{0.800000}%
\pgfsetdash{}{0pt}%
\pgfpathmoveto{\pgfqpoint{1.400977in}{1.349984in}}%
\pgfpathlineto{\pgfqpoint{2.626197in}{1.349984in}}%
\pgfpathquadraticcurveto{\pgfqpoint{2.653974in}{1.349984in}}{\pgfqpoint{2.653974in}{1.377762in}}%
\pgfpathlineto{\pgfqpoint{2.653974in}{1.975445in}}%
\pgfpathquadraticcurveto{\pgfqpoint{2.653974in}{2.003223in}}{\pgfqpoint{2.626197in}{2.003223in}}%
\pgfpathlineto{\pgfqpoint{1.400977in}{2.003223in}}%
\pgfpathquadraticcurveto{\pgfqpoint{1.373199in}{2.003223in}}{\pgfqpoint{1.373199in}{1.975445in}}%
\pgfpathlineto{\pgfqpoint{1.373199in}{1.377762in}}%
\pgfpathquadraticcurveto{\pgfqpoint{1.373199in}{1.349984in}}{\pgfqpoint{1.400977in}{1.349984in}}%
\pgfpathlineto{\pgfqpoint{1.400977in}{1.349984in}}%
\pgfpathclose%
\pgfusepath{stroke,fill}%
\end{pgfscope}%
\begin{pgfscope}%
\pgfsetrectcap%
\pgfsetroundjoin%
\pgfsetlinewidth{1.505625pt}%
\definecolor{currentstroke}{rgb}{0.121569,0.466667,0.705882}%
\pgfsetstrokecolor{currentstroke}%
\pgfsetdash{}{0pt}%
\pgfpathmoveto{\pgfqpoint{1.428755in}{1.890755in}}%
\pgfpathlineto{\pgfqpoint{1.567644in}{1.890755in}}%
\pgfpathlineto{\pgfqpoint{1.706532in}{1.890755in}}%
\pgfusepath{stroke}%
\end{pgfscope}%
\begin{pgfscope}%
\definecolor{textcolor}{rgb}{0.000000,0.000000,0.000000}%
\pgfsetstrokecolor{textcolor}%
\pgfsetfillcolor{textcolor}%
\pgftext[x=1.817644in,y=1.842144in,left,base]{\color{textcolor}{\sffamily\fontsize{10.000000}{12.000000}\selectfont\catcode`\^=\active\def^{\ifmmode\sp\else\^{}\fi}\catcode`\%=\active\def
\end{pgfscope}%
\begin{pgfscope}%
\pgfsetrectcap%
\pgfsetroundjoin%
\pgfsetlinewidth{1.505625pt}%
\definecolor{currentstroke}{rgb}{1.000000,0.498039,0.054902}%
\pgfsetstrokecolor{currentstroke}%
\pgfsetdash{}{0pt}%
\pgfpathmoveto{\pgfqpoint{1.428755in}{1.686898in}}%
\pgfpathlineto{\pgfqpoint{1.567644in}{1.686898in}}%
\pgfpathlineto{\pgfqpoint{1.706532in}{1.686898in}}%
\pgfusepath{stroke}%
\end{pgfscope}%
\begin{pgfscope}%
\definecolor{textcolor}{rgb}{0.000000,0.000000,0.000000}%
\pgfsetstrokecolor{textcolor}%
\pgfsetfillcolor{textcolor}%
\pgftext[x=1.817644in,y=1.638287in,left,base]{\color{textcolor}{\sffamily\fontsize{10.000000}{12.000000}\selectfont\catcode`\^=\active\def^{\ifmmode\sp\else\^{}\fi}\catcode`\%=\active\def
\end{pgfscope}%
\begin{pgfscope}%
\pgfsetrectcap%
\pgfsetroundjoin%
\pgfsetlinewidth{1.505625pt}%
\definecolor{currentstroke}{rgb}{0.172549,0.627451,0.172549}%
\pgfsetstrokecolor{currentstroke}%
\pgfsetdash{}{0pt}%
\pgfpathmoveto{\pgfqpoint{1.428755in}{1.483041in}}%
\pgfpathlineto{\pgfqpoint{1.567644in}{1.483041in}}%
\pgfpathlineto{\pgfqpoint{1.706532in}{1.483041in}}%
\pgfusepath{stroke}%
\end{pgfscope}%
\begin{pgfscope}%
\definecolor{textcolor}{rgb}{0.000000,0.000000,0.000000}%
\pgfsetstrokecolor{textcolor}%
\pgfsetfillcolor{textcolor}%
\pgftext[x=1.817644in,y=1.434429in,left,base]{\color{textcolor}{\sffamily\fontsize{10.000000}{12.000000}\selectfont\catcode`\^=\active\def^{\ifmmode\sp\else\^{}\fi}\catcode`\%=\active\def
\end{pgfscope}%
\end{pgfpicture}%
\makeatother%
\endgroup%

%% file: synth_spec.pgf
\begingroup%
\makeatletter%
\begin{pgfpicture}%
\pgfpathrectangle{\pgfpointorigin}{\pgfqpoint{1.410244in}{2.931603in}}%
\pgfusepath{use as bounding box, clip}%
\begin{pgfscope}%
\pgfsetbuttcap%
\pgfsetmiterjoin%
\definecolor{currentfill}{rgb}{1.000000,1.000000,1.000000}%
\pgfsetfillcolor{currentfill}%
\pgfsetlinewidth{0.000000pt}%
\definecolor{currentstroke}{rgb}{1.000000,1.000000,1.000000}%
\pgfsetstrokecolor{currentstroke}%
\pgfsetdash{}{0pt}%
\pgfpathmoveto{\pgfqpoint{0.000000in}{0.000000in}}%
\pgfpathlineto{\pgfqpoint{1.410244in}{0.000000in}}%
\pgfpathlineto{\pgfqpoint{1.410244in}{2.931603in}}%
\pgfpathlineto{\pgfqpoint{0.000000in}{2.931603in}}%
\pgfpathlineto{\pgfqpoint{0.000000in}{0.000000in}}%
\pgfpathclose%
\pgfusepath{fill}%
\end{pgfscope}%
\begin{pgfscope}%
\pgfsetbuttcap%
\pgfsetmiterjoin%
\definecolor{currentfill}{rgb}{1.000000,1.000000,1.000000}%
\pgfsetfillcolor{currentfill}%
\pgfsetlinewidth{0.000000pt}%
\definecolor{currentstroke}{rgb}{0.000000,0.000000,0.000000}%
\pgfsetstrokecolor{currentstroke}%
\pgfsetstrokeopacity{0.000000}%
\pgfsetdash{}{0pt}%
\pgfpathmoveto{\pgfqpoint{0.289968in}{2.433327in}}%
\pgfpathlineto{\pgfqpoint{1.219968in}{2.433327in}}%
\pgfpathlineto{\pgfqpoint{1.219968in}{2.831603in}}%
\pgfpathlineto{\pgfqpoint{0.289968in}{2.831603in}}%
\pgfpathlineto{\pgfqpoint{0.289968in}{2.433327in}}%
\pgfpathclose%
\pgfusepath{fill}%
\end{pgfscope}%
\begin{pgfscope}%
\pgfsetbuttcap%
\pgfsetroundjoin%
\definecolor{currentfill}{rgb}{0.000000,0.000000,0.000000}%
\pgfsetfillcolor{currentfill}%
\pgfsetlinewidth{0.803000pt}%
\definecolor{currentstroke}{rgb}{0.000000,0.000000,0.000000}%
\pgfsetstrokecolor{currentstroke}%
\pgfsetdash{}{0pt}%
\pgfsys@defobject{currentmarker}{\pgfqpoint{0.000000in}{-0.048611in}}{\pgfqpoint{0.000000in}{0.000000in}}{%
\pgfpathmoveto{\pgfqpoint{0.000000in}{0.000000in}}%
\pgfpathlineto{\pgfqpoint{0.000000in}{-0.048611in}}%
\pgfusepath{stroke,fill}%
}%
\begin{pgfscope}%
\pgfsys@transformshift{0.332241in}{2.433327in}%
\pgfsys@useobject{currentmarker}{}%
\end{pgfscope}%
\end{pgfscope}%
\begin{pgfscope}%
\pgfsetbuttcap%
\pgfsetroundjoin%
\definecolor{currentfill}{rgb}{0.000000,0.000000,0.000000}%
\pgfsetfillcolor{currentfill}%
\pgfsetlinewidth{0.803000pt}%
\definecolor{currentstroke}{rgb}{0.000000,0.000000,0.000000}%
\pgfsetstrokecolor{currentstroke}%
\pgfsetdash{}{0pt}%
\pgfsys@defobject{currentmarker}{\pgfqpoint{0.000000in}{-0.048611in}}{\pgfqpoint{0.000000in}{0.000000in}}{%
\pgfpathmoveto{\pgfqpoint{0.000000in}{0.000000in}}%
\pgfpathlineto{\pgfqpoint{0.000000in}{-0.048611in}}%
\pgfusepath{stroke,fill}%
}%
\begin{pgfscope}%
\pgfsys@transformshift{0.707999in}{2.433327in}%
\pgfsys@useobject{currentmarker}{}%
\end{pgfscope}%
\end{pgfscope}%
\begin{pgfscope}%
\pgfsetbuttcap%
\pgfsetroundjoin%
\definecolor{currentfill}{rgb}{0.000000,0.000000,0.000000}%
\pgfsetfillcolor{currentfill}%
\pgfsetlinewidth{0.803000pt}%
\definecolor{currentstroke}{rgb}{0.000000,0.000000,0.000000}%
\pgfsetstrokecolor{currentstroke}%
\pgfsetdash{}{0pt}%
\pgfsys@defobject{currentmarker}{\pgfqpoint{0.000000in}{-0.048611in}}{\pgfqpoint{0.000000in}{0.000000in}}{%
\pgfpathmoveto{\pgfqpoint{0.000000in}{0.000000in}}%
\pgfpathlineto{\pgfqpoint{0.000000in}{-0.048611in}}%
\pgfusepath{stroke,fill}%
}%
\begin{pgfscope}%
\pgfsys@transformshift{1.177696in}{2.433327in}%
\pgfsys@useobject{currentmarker}{}%
\end{pgfscope}%
\end{pgfscope}%
\begin{pgfscope}%
\definecolor{textcolor}{rgb}{0.000000,0.000000,0.000000}%
\pgfsetstrokecolor{textcolor}%
\pgfsetfillcolor{textcolor}%
\pgftext[x=0.234413in,y=2.632465in,,bottom,rotate=90.000000]{\color{textcolor}{\sffamily\fontsize{10.000000}{12.000000}\selectfont\catcode`\^=\active\def^{\ifmmode\sp\else\^{}\fi}\catcode`\%=\active\def
\end{pgfscope}%
\begin{pgfscope}%
\pgfpathrectangle{\pgfqpoint{0.289968in}{2.433327in}}{\pgfqpoint{0.930000in}{0.398276in}}%
\pgfusepath{clip}%
\pgfsetrectcap%
\pgfsetroundjoin%
\pgfsetlinewidth{1.505625pt}%
\definecolor{currentstroke}{rgb}{0.121569,0.466667,0.705882}%
\pgfsetstrokecolor{currentstroke}%
\pgfsetdash{}{0pt}%
\pgfpathmoveto{\pgfqpoint{0.332241in}{2.629236in}}%
\pgfpathlineto{\pgfqpoint{0.426180in}{2.672451in}}%
\pgfpathlineto{\pgfqpoint{0.520120in}{2.623103in}}%
\pgfpathlineto{\pgfqpoint{0.614059in}{2.542197in}}%
\pgfpathlineto{\pgfqpoint{0.707999in}{2.491152in}}%
\pgfpathlineto{\pgfqpoint{0.801938in}{2.475851in}}%
\pgfpathlineto{\pgfqpoint{0.895877in}{2.475720in}}%
\pgfpathlineto{\pgfqpoint{0.989817in}{2.481290in}}%
\pgfpathlineto{\pgfqpoint{1.083756in}{2.484423in}}%
\pgfpathlineto{\pgfqpoint{1.177696in}{2.491088in}}%
\pgfusepath{stroke}%
\end{pgfscope}%
\begin{pgfscope}%
\pgfpathrectangle{\pgfqpoint{0.289968in}{2.433327in}}{\pgfqpoint{0.930000in}{0.398276in}}%
\pgfusepath{clip}%
\pgfsetrectcap%
\pgfsetroundjoin%
\pgfsetlinewidth{1.505625pt}%
\definecolor{currentstroke}{rgb}{1.000000,0.498039,0.054902}%
\pgfsetstrokecolor{currentstroke}%
\pgfsetdash{}{0pt}%
\pgfpathmoveto{\pgfqpoint{0.332241in}{2.662422in}}%
\pgfpathlineto{\pgfqpoint{0.426180in}{2.676486in}}%
\pgfpathlineto{\pgfqpoint{0.520120in}{2.587562in}}%
\pgfpathlineto{\pgfqpoint{0.614059in}{2.518726in}}%
\pgfpathlineto{\pgfqpoint{0.707999in}{2.497447in}}%
\pgfpathlineto{\pgfqpoint{0.801938in}{2.510251in}}%
\pgfpathlineto{\pgfqpoint{0.895877in}{2.477728in}}%
\pgfpathlineto{\pgfqpoint{0.989817in}{2.456788in}}%
\pgfpathlineto{\pgfqpoint{1.083756in}{2.468149in}}%
\pgfpathlineto{\pgfqpoint{1.177696in}{2.479901in}}%
\pgfusepath{stroke}%
\end{pgfscope}%
\begin{pgfscope}%
\pgfsetrectcap%
\pgfsetmiterjoin%
\pgfsetlinewidth{0.803000pt}%
\definecolor{currentstroke}{rgb}{0.000000,0.000000,0.000000}%
\pgfsetstrokecolor{currentstroke}%
\pgfsetdash{}{0pt}%
\pgfpathmoveto{\pgfqpoint{0.289968in}{2.433327in}}%
\pgfpathlineto{\pgfqpoint{0.289968in}{2.831603in}}%
\pgfusepath{stroke}%
\end{pgfscope}%
\begin{pgfscope}%
\pgfsetrectcap%
\pgfsetmiterjoin%
\pgfsetlinewidth{0.803000pt}%
\definecolor{currentstroke}{rgb}{0.000000,0.000000,0.000000}%
\pgfsetstrokecolor{currentstroke}%
\pgfsetdash{}{0pt}%
\pgfpathmoveto{\pgfqpoint{1.219968in}{2.433327in}}%
\pgfpathlineto{\pgfqpoint{1.219968in}{2.831603in}}%
\pgfusepath{stroke}%
\end{pgfscope}%
\begin{pgfscope}%
\pgfsetrectcap%
\pgfsetmiterjoin%
\pgfsetlinewidth{0.803000pt}%
\definecolor{currentstroke}{rgb}{0.000000,0.000000,0.000000}%
\pgfsetstrokecolor{currentstroke}%
\pgfsetdash{}{0pt}%
\pgfpathmoveto{\pgfqpoint{0.289968in}{2.433327in}}%
\pgfpathlineto{\pgfqpoint{1.219968in}{2.433327in}}%
\pgfusepath{stroke}%
\end{pgfscope}%
\begin{pgfscope}%
\pgfsetrectcap%
\pgfsetmiterjoin%
\pgfsetlinewidth{0.803000pt}%
\definecolor{currentstroke}{rgb}{0.000000,0.000000,0.000000}%
\pgfsetstrokecolor{currentstroke}%
\pgfsetdash{}{0pt}%
\pgfpathmoveto{\pgfqpoint{0.289968in}{2.831603in}}%
\pgfpathlineto{\pgfqpoint{1.219968in}{2.831603in}}%
\pgfusepath{stroke}%
\end{pgfscope}%
\begin{pgfscope}%
\pgfsetbuttcap%
\pgfsetmiterjoin%
\definecolor{currentfill}{rgb}{1.000000,1.000000,1.000000}%
\pgfsetfillcolor{currentfill}%
\pgfsetlinewidth{0.000000pt}%
\definecolor{currentstroke}{rgb}{0.000000,0.000000,0.000000}%
\pgfsetstrokecolor{currentstroke}%
\pgfsetstrokeopacity{0.000000}%
\pgfsetdash{}{0pt}%
\pgfpathmoveto{\pgfqpoint{0.289968in}{1.955396in}}%
\pgfpathlineto{\pgfqpoint{1.219968in}{1.955396in}}%
\pgfpathlineto{\pgfqpoint{1.219968in}{2.353672in}}%
\pgfpathlineto{\pgfqpoint{0.289968in}{2.353672in}}%
\pgfpathlineto{\pgfqpoint{0.289968in}{1.955396in}}%
\pgfpathclose%
\pgfusepath{fill}%
\end{pgfscope}%
\begin{pgfscope}%
\pgfsetbuttcap%
\pgfsetroundjoin%
\definecolor{currentfill}{rgb}{0.000000,0.000000,0.000000}%
\pgfsetfillcolor{currentfill}%
\pgfsetlinewidth{0.803000pt}%
\definecolor{currentstroke}{rgb}{0.000000,0.000000,0.000000}%
\pgfsetstrokecolor{currentstroke}%
\pgfsetdash{}{0pt}%
\pgfsys@defobject{currentmarker}{\pgfqpoint{0.000000in}{-0.048611in}}{\pgfqpoint{0.000000in}{0.000000in}}{%
\pgfpathmoveto{\pgfqpoint{0.000000in}{0.000000in}}%
\pgfpathlineto{\pgfqpoint{0.000000in}{-0.048611in}}%
\pgfusepath{stroke,fill}%
}%
\begin{pgfscope}%
\pgfsys@transformshift{0.332241in}{1.955396in}%
\pgfsys@useobject{currentmarker}{}%
\end{pgfscope}%
\end{pgfscope}%
\begin{pgfscope}%
\pgfsetbuttcap%
\pgfsetroundjoin%
\definecolor{currentfill}{rgb}{0.000000,0.000000,0.000000}%
\pgfsetfillcolor{currentfill}%
\pgfsetlinewidth{0.803000pt}%
\definecolor{currentstroke}{rgb}{0.000000,0.000000,0.000000}%
\pgfsetstrokecolor{currentstroke}%
\pgfsetdash{}{0pt}%
\pgfsys@defobject{currentmarker}{\pgfqpoint{0.000000in}{-0.048611in}}{\pgfqpoint{0.000000in}{0.000000in}}{%
\pgfpathmoveto{\pgfqpoint{0.000000in}{0.000000in}}%
\pgfpathlineto{\pgfqpoint{0.000000in}{-0.048611in}}%
\pgfusepath{stroke,fill}%
}%
\begin{pgfscope}%
\pgfsys@transformshift{0.707999in}{1.955396in}%
\pgfsys@useobject{currentmarker}{}%
\end{pgfscope}%
\end{pgfscope}%
\begin{pgfscope}%
\pgfsetbuttcap%
\pgfsetroundjoin%
\definecolor{currentfill}{rgb}{0.000000,0.000000,0.000000}%
\pgfsetfillcolor{currentfill}%
\pgfsetlinewidth{0.803000pt}%
\definecolor{currentstroke}{rgb}{0.000000,0.000000,0.000000}%
\pgfsetstrokecolor{currentstroke}%
\pgfsetdash{}{0pt}%
\pgfsys@defobject{currentmarker}{\pgfqpoint{0.000000in}{-0.048611in}}{\pgfqpoint{0.000000in}{0.000000in}}{%
\pgfpathmoveto{\pgfqpoint{0.000000in}{0.000000in}}%
\pgfpathlineto{\pgfqpoint{0.000000in}{-0.048611in}}%
\pgfusepath{stroke,fill}%
}%
\begin{pgfscope}%
\pgfsys@transformshift{1.177696in}{1.955396in}%
\pgfsys@useobject{currentmarker}{}%
\end{pgfscope}%
\end{pgfscope}%
\begin{pgfscope}%
\definecolor{textcolor}{rgb}{0.000000,0.000000,0.000000}%
\pgfsetstrokecolor{textcolor}%
\pgfsetfillcolor{textcolor}%
\pgftext[x=0.234413in,y=2.154534in,,bottom,rotate=90.000000]{\color{textcolor}{\sffamily\fontsize{10.000000}{12.000000}\selectfont\catcode`\^=\active\def^{\ifmmode\sp\else\^{}\fi}\catcode`\%=\active\def
\end{pgfscope}%
\begin{pgfscope}%
\pgfpathrectangle{\pgfqpoint{0.289968in}{1.955396in}}{\pgfqpoint{0.930000in}{0.398276in}}%
\pgfusepath{clip}%
\pgfsetrectcap%
\pgfsetroundjoin%
\pgfsetlinewidth{1.505625pt}%
\definecolor{currentstroke}{rgb}{0.121569,0.466667,0.705882}%
\pgfsetstrokecolor{currentstroke}%
\pgfsetdash{}{0pt}%
\pgfpathmoveto{\pgfqpoint{0.332241in}{1.969992in}}%
\pgfpathlineto{\pgfqpoint{0.426180in}{1.970169in}}%
\pgfpathlineto{\pgfqpoint{0.520120in}{1.988600in}}%
\pgfpathlineto{\pgfqpoint{0.614059in}{2.050079in}}%
\pgfpathlineto{\pgfqpoint{0.707999in}{2.100134in}}%
\pgfpathlineto{\pgfqpoint{0.801938in}{2.114248in}}%
\pgfpathlineto{\pgfqpoint{0.895877in}{2.115129in}}%
\pgfpathlineto{\pgfqpoint{0.989817in}{2.116960in}}%
\pgfpathlineto{\pgfqpoint{1.083756in}{2.113550in}}%
\pgfpathlineto{\pgfqpoint{1.177696in}{2.114495in}}%
\pgfusepath{stroke}%
\end{pgfscope}%
\begin{pgfscope}%
\pgfpathrectangle{\pgfqpoint{0.289968in}{1.955396in}}{\pgfqpoint{0.930000in}{0.398276in}}%
\pgfusepath{clip}%
\pgfsetrectcap%
\pgfsetroundjoin%
\pgfsetlinewidth{1.505625pt}%
\definecolor{currentstroke}{rgb}{1.000000,0.498039,0.054902}%
\pgfsetstrokecolor{currentstroke}%
\pgfsetdash{}{0pt}%
\pgfpathmoveto{\pgfqpoint{0.332241in}{1.959499in}}%
\pgfpathlineto{\pgfqpoint{0.426180in}{1.972308in}}%
\pgfpathlineto{\pgfqpoint{0.520120in}{2.021401in}}%
\pgfpathlineto{\pgfqpoint{0.614059in}{2.067275in}}%
\pgfpathlineto{\pgfqpoint{0.707999in}{2.058037in}}%
\pgfpathlineto{\pgfqpoint{0.801938in}{2.118615in}}%
\pgfpathlineto{\pgfqpoint{0.895877in}{2.110843in}}%
\pgfpathlineto{\pgfqpoint{0.989817in}{2.115504in}}%
\pgfpathlineto{\pgfqpoint{1.083756in}{2.136764in}}%
\pgfpathlineto{\pgfqpoint{1.177696in}{2.102297in}}%
\pgfusepath{stroke}%
\end{pgfscope}%
\begin{pgfscope}%
\pgfsetrectcap%
\pgfsetmiterjoin%
\pgfsetlinewidth{0.803000pt}%
\definecolor{currentstroke}{rgb}{0.000000,0.000000,0.000000}%
\pgfsetstrokecolor{currentstroke}%
\pgfsetdash{}{0pt}%
\pgfpathmoveto{\pgfqpoint{0.289968in}{1.955396in}}%
\pgfpathlineto{\pgfqpoint{0.289968in}{2.353672in}}%
\pgfusepath{stroke}%
\end{pgfscope}%
\begin{pgfscope}%
\pgfsetrectcap%
\pgfsetmiterjoin%
\pgfsetlinewidth{0.803000pt}%
\definecolor{currentstroke}{rgb}{0.000000,0.000000,0.000000}%
\pgfsetstrokecolor{currentstroke}%
\pgfsetdash{}{0pt}%
\pgfpathmoveto{\pgfqpoint{1.219968in}{1.955396in}}%
\pgfpathlineto{\pgfqpoint{1.219968in}{2.353672in}}%
\pgfusepath{stroke}%
\end{pgfscope}%
\begin{pgfscope}%
\pgfsetrectcap%
\pgfsetmiterjoin%
\pgfsetlinewidth{0.803000pt}%
\definecolor{currentstroke}{rgb}{0.000000,0.000000,0.000000}%
\pgfsetstrokecolor{currentstroke}%
\pgfsetdash{}{0pt}%
\pgfpathmoveto{\pgfqpoint{0.289968in}{1.955396in}}%
\pgfpathlineto{\pgfqpoint{1.219968in}{1.955396in}}%
\pgfusepath{stroke}%
\end{pgfscope}%
\begin{pgfscope}%
\pgfsetrectcap%
\pgfsetmiterjoin%
\pgfsetlinewidth{0.803000pt}%
\definecolor{currentstroke}{rgb}{0.000000,0.000000,0.000000}%
\pgfsetstrokecolor{currentstroke}%
\pgfsetdash{}{0pt}%
\pgfpathmoveto{\pgfqpoint{0.289968in}{2.353672in}}%
\pgfpathlineto{\pgfqpoint{1.219968in}{2.353672in}}%
\pgfusepath{stroke}%
\end{pgfscope}%
\begin{pgfscope}%
\pgfsetbuttcap%
\pgfsetmiterjoin%
\definecolor{currentfill}{rgb}{1.000000,1.000000,1.000000}%
\pgfsetfillcolor{currentfill}%
\pgfsetlinewidth{0.000000pt}%
\definecolor{currentstroke}{rgb}{0.000000,0.000000,0.000000}%
\pgfsetstrokecolor{currentstroke}%
\pgfsetstrokeopacity{0.000000}%
\pgfsetdash{}{0pt}%
\pgfpathmoveto{\pgfqpoint{0.289968in}{1.477465in}}%
\pgfpathlineto{\pgfqpoint{1.219968in}{1.477465in}}%
\pgfpathlineto{\pgfqpoint{1.219968in}{1.875741in}}%
\pgfpathlineto{\pgfqpoint{0.289968in}{1.875741in}}%
\pgfpathlineto{\pgfqpoint{0.289968in}{1.477465in}}%
\pgfpathclose%
\pgfusepath{fill}%
\end{pgfscope}%
\begin{pgfscope}%
\pgfsetbuttcap%
\pgfsetroundjoin%
\definecolor{currentfill}{rgb}{0.000000,0.000000,0.000000}%
\pgfsetfillcolor{currentfill}%
\pgfsetlinewidth{0.803000pt}%
\definecolor{currentstroke}{rgb}{0.000000,0.000000,0.000000}%
\pgfsetstrokecolor{currentstroke}%
\pgfsetdash{}{0pt}%
\pgfsys@defobject{currentmarker}{\pgfqpoint{0.000000in}{-0.048611in}}{\pgfqpoint{0.000000in}{0.000000in}}{%
\pgfpathmoveto{\pgfqpoint{0.000000in}{0.000000in}}%
\pgfpathlineto{\pgfqpoint{0.000000in}{-0.048611in}}%
\pgfusepath{stroke,fill}%
}%
\begin{pgfscope}%
\pgfsys@transformshift{0.332241in}{1.477465in}%
\pgfsys@useobject{currentmarker}{}%
\end{pgfscope}%
\end{pgfscope}%
\begin{pgfscope}%
\pgfsetbuttcap%
\pgfsetroundjoin%
\definecolor{currentfill}{rgb}{0.000000,0.000000,0.000000}%
\pgfsetfillcolor{currentfill}%
\pgfsetlinewidth{0.803000pt}%
\definecolor{currentstroke}{rgb}{0.000000,0.000000,0.000000}%
\pgfsetstrokecolor{currentstroke}%
\pgfsetdash{}{0pt}%
\pgfsys@defobject{currentmarker}{\pgfqpoint{0.000000in}{-0.048611in}}{\pgfqpoint{0.000000in}{0.000000in}}{%
\pgfpathmoveto{\pgfqpoint{0.000000in}{0.000000in}}%
\pgfpathlineto{\pgfqpoint{0.000000in}{-0.048611in}}%
\pgfusepath{stroke,fill}%
}%
\begin{pgfscope}%
\pgfsys@transformshift{0.707999in}{1.477465in}%
\pgfsys@useobject{currentmarker}{}%
\end{pgfscope}%
\end{pgfscope}%
\begin{pgfscope}%
\pgfsetbuttcap%
\pgfsetroundjoin%
\definecolor{currentfill}{rgb}{0.000000,0.000000,0.000000}%
\pgfsetfillcolor{currentfill}%
\pgfsetlinewidth{0.803000pt}%
\definecolor{currentstroke}{rgb}{0.000000,0.000000,0.000000}%
\pgfsetstrokecolor{currentstroke}%
\pgfsetdash{}{0pt}%
\pgfsys@defobject{currentmarker}{\pgfqpoint{0.000000in}{-0.048611in}}{\pgfqpoint{0.000000in}{0.000000in}}{%
\pgfpathmoveto{\pgfqpoint{0.000000in}{0.000000in}}%
\pgfpathlineto{\pgfqpoint{0.000000in}{-0.048611in}}%
\pgfusepath{stroke,fill}%
}%
\begin{pgfscope}%
\pgfsys@transformshift{1.177696in}{1.477465in}%
\pgfsys@useobject{currentmarker}{}%
\end{pgfscope}%
\end{pgfscope}%
\begin{pgfscope}%
\definecolor{textcolor}{rgb}{0.000000,0.000000,0.000000}%
\pgfsetstrokecolor{textcolor}%
\pgfsetfillcolor{textcolor}%
\pgftext[x=0.234413in,y=1.676603in,,bottom,rotate=90.000000]{\color{textcolor}{\sffamily\fontsize{10.000000}{12.000000}\selectfont\catcode`\^=\active\def^{\ifmmode\sp\else\^{}\fi}\catcode`\%=\active\def
\end{pgfscope}%
\begin{pgfscope}%
\pgfpathrectangle{\pgfqpoint{0.289968in}{1.477465in}}{\pgfqpoint{0.930000in}{0.398276in}}%
\pgfusepath{clip}%
\pgfsetrectcap%
\pgfsetroundjoin%
\pgfsetlinewidth{1.505625pt}%
\definecolor{currentstroke}{rgb}{0.121569,0.466667,0.705882}%
\pgfsetstrokecolor{currentstroke}%
\pgfsetdash{}{0pt}%
\pgfpathmoveto{\pgfqpoint{0.332241in}{1.490204in}}%
\pgfpathlineto{\pgfqpoint{0.426180in}{1.488709in}}%
\pgfpathlineto{\pgfqpoint{0.520120in}{1.489731in}}%
\pgfpathlineto{\pgfqpoint{0.614059in}{1.500970in}}%
\pgfpathlineto{\pgfqpoint{0.707999in}{1.506417in}}%
\pgfpathlineto{\pgfqpoint{0.801938in}{1.567035in}}%
\pgfpathlineto{\pgfqpoint{0.895877in}{1.648128in}}%
\pgfpathlineto{\pgfqpoint{0.989817in}{1.674560in}}%
\pgfpathlineto{\pgfqpoint{1.083756in}{1.680056in}}%
\pgfpathlineto{\pgfqpoint{1.177696in}{1.676855in}}%
\pgfusepath{stroke}%
\end{pgfscope}%
\begin{pgfscope}%
\pgfpathrectangle{\pgfqpoint{0.289968in}{1.477465in}}{\pgfqpoint{0.930000in}{0.398276in}}%
\pgfusepath{clip}%
\pgfsetrectcap%
\pgfsetroundjoin%
\pgfsetlinewidth{1.505625pt}%
\definecolor{currentstroke}{rgb}{1.000000,0.498039,0.054902}%
\pgfsetstrokecolor{currentstroke}%
\pgfsetdash{}{0pt}%
\pgfpathmoveto{\pgfqpoint{0.332241in}{1.482641in}}%
\pgfpathlineto{\pgfqpoint{0.426180in}{1.505778in}}%
\pgfpathlineto{\pgfqpoint{0.520120in}{1.532308in}}%
\pgfpathlineto{\pgfqpoint{0.614059in}{1.506852in}}%
\pgfpathlineto{\pgfqpoint{0.707999in}{1.517177in}}%
\pgfpathlineto{\pgfqpoint{0.801938in}{1.573747in}}%
\pgfpathlineto{\pgfqpoint{0.895877in}{1.638194in}}%
\pgfpathlineto{\pgfqpoint{0.989817in}{1.669980in}}%
\pgfpathlineto{\pgfqpoint{1.083756in}{1.683299in}}%
\pgfpathlineto{\pgfqpoint{1.177696in}{1.671923in}}%
\pgfusepath{stroke}%
\end{pgfscope}%
\begin{pgfscope}%
\pgfsetrectcap%
\pgfsetmiterjoin%
\pgfsetlinewidth{0.803000pt}%
\definecolor{currentstroke}{rgb}{0.000000,0.000000,0.000000}%
\pgfsetstrokecolor{currentstroke}%
\pgfsetdash{}{0pt}%
\pgfpathmoveto{\pgfqpoint{0.289968in}{1.477465in}}%
\pgfpathlineto{\pgfqpoint{0.289968in}{1.875741in}}%
\pgfusepath{stroke}%
\end{pgfscope}%
\begin{pgfscope}%
\pgfsetrectcap%
\pgfsetmiterjoin%
\pgfsetlinewidth{0.803000pt}%
\definecolor{currentstroke}{rgb}{0.000000,0.000000,0.000000}%
\pgfsetstrokecolor{currentstroke}%
\pgfsetdash{}{0pt}%
\pgfpathmoveto{\pgfqpoint{1.219968in}{1.477465in}}%
\pgfpathlineto{\pgfqpoint{1.219968in}{1.875741in}}%
\pgfusepath{stroke}%
\end{pgfscope}%
\begin{pgfscope}%
\pgfsetrectcap%
\pgfsetmiterjoin%
\pgfsetlinewidth{0.803000pt}%
\definecolor{currentstroke}{rgb}{0.000000,0.000000,0.000000}%
\pgfsetstrokecolor{currentstroke}%
\pgfsetdash{}{0pt}%
\pgfpathmoveto{\pgfqpoint{0.289968in}{1.477465in}}%
\pgfpathlineto{\pgfqpoint{1.219968in}{1.477465in}}%
\pgfusepath{stroke}%
\end{pgfscope}%
\begin{pgfscope}%
\pgfsetrectcap%
\pgfsetmiterjoin%
\pgfsetlinewidth{0.803000pt}%
\definecolor{currentstroke}{rgb}{0.000000,0.000000,0.000000}%
\pgfsetstrokecolor{currentstroke}%
\pgfsetdash{}{0pt}%
\pgfpathmoveto{\pgfqpoint{0.289968in}{1.875741in}}%
\pgfpathlineto{\pgfqpoint{1.219968in}{1.875741in}}%
\pgfusepath{stroke}%
\end{pgfscope}%
\begin{pgfscope}%
\pgfsetbuttcap%
\pgfsetmiterjoin%
\definecolor{currentfill}{rgb}{1.000000,1.000000,1.000000}%
\pgfsetfillcolor{currentfill}%
\pgfsetlinewidth{0.000000pt}%
\definecolor{currentstroke}{rgb}{0.000000,0.000000,0.000000}%
\pgfsetstrokecolor{currentstroke}%
\pgfsetstrokeopacity{0.000000}%
\pgfsetdash{}{0pt}%
\pgfpathmoveto{\pgfqpoint{0.289968in}{0.999534in}}%
\pgfpathlineto{\pgfqpoint{1.219968in}{0.999534in}}%
\pgfpathlineto{\pgfqpoint{1.219968in}{1.397810in}}%
\pgfpathlineto{\pgfqpoint{0.289968in}{1.397810in}}%
\pgfpathlineto{\pgfqpoint{0.289968in}{0.999534in}}%
\pgfpathclose%
\pgfusepath{fill}%
\end{pgfscope}%
\begin{pgfscope}%
\pgfsetbuttcap%
\pgfsetroundjoin%
\definecolor{currentfill}{rgb}{0.000000,0.000000,0.000000}%
\pgfsetfillcolor{currentfill}%
\pgfsetlinewidth{0.803000pt}%
\definecolor{currentstroke}{rgb}{0.000000,0.000000,0.000000}%
\pgfsetstrokecolor{currentstroke}%
\pgfsetdash{}{0pt}%
\pgfsys@defobject{currentmarker}{\pgfqpoint{0.000000in}{-0.048611in}}{\pgfqpoint{0.000000in}{0.000000in}}{%
\pgfpathmoveto{\pgfqpoint{0.000000in}{0.000000in}}%
\pgfpathlineto{\pgfqpoint{0.000000in}{-0.048611in}}%
\pgfusepath{stroke,fill}%
}%
\begin{pgfscope}%
\pgfsys@transformshift{0.332241in}{0.999534in}%
\pgfsys@useobject{currentmarker}{}%
\end{pgfscope}%
\end{pgfscope}%
\begin{pgfscope}%
\pgfsetbuttcap%
\pgfsetroundjoin%
\definecolor{currentfill}{rgb}{0.000000,0.000000,0.000000}%
\pgfsetfillcolor{currentfill}%
\pgfsetlinewidth{0.803000pt}%
\definecolor{currentstroke}{rgb}{0.000000,0.000000,0.000000}%
\pgfsetstrokecolor{currentstroke}%
\pgfsetdash{}{0pt}%
\pgfsys@defobject{currentmarker}{\pgfqpoint{0.000000in}{-0.048611in}}{\pgfqpoint{0.000000in}{0.000000in}}{%
\pgfpathmoveto{\pgfqpoint{0.000000in}{0.000000in}}%
\pgfpathlineto{\pgfqpoint{0.000000in}{-0.048611in}}%
\pgfusepath{stroke,fill}%
}%
\begin{pgfscope}%
\pgfsys@transformshift{0.707999in}{0.999534in}%
\pgfsys@useobject{currentmarker}{}%
\end{pgfscope}%
\end{pgfscope}%
\begin{pgfscope}%
\pgfsetbuttcap%
\pgfsetroundjoin%
\definecolor{currentfill}{rgb}{0.000000,0.000000,0.000000}%
\pgfsetfillcolor{currentfill}%
\pgfsetlinewidth{0.803000pt}%
\definecolor{currentstroke}{rgb}{0.000000,0.000000,0.000000}%
\pgfsetstrokecolor{currentstroke}%
\pgfsetdash{}{0pt}%
\pgfsys@defobject{currentmarker}{\pgfqpoint{0.000000in}{-0.048611in}}{\pgfqpoint{0.000000in}{0.000000in}}{%
\pgfpathmoveto{\pgfqpoint{0.000000in}{0.000000in}}%
\pgfpathlineto{\pgfqpoint{0.000000in}{-0.048611in}}%
\pgfusepath{stroke,fill}%
}%
\begin{pgfscope}%
\pgfsys@transformshift{1.177696in}{0.999534in}%
\pgfsys@useobject{currentmarker}{}%
\end{pgfscope}%
\end{pgfscope}%
\begin{pgfscope}%
\definecolor{textcolor}{rgb}{0.000000,0.000000,0.000000}%
\pgfsetstrokecolor{textcolor}%
\pgfsetfillcolor{textcolor}%
\pgftext[x=0.234413in,y=1.198672in,,bottom,rotate=90.000000]{\color{textcolor}{\sffamily\fontsize{10.000000}{12.000000}\selectfont\catcode`\^=\active\def^{\ifmmode\sp\else\^{}\fi}\catcode`\%=\active\def
\end{pgfscope}%
\begin{pgfscope}%
\pgfpathrectangle{\pgfqpoint{0.289968in}{0.999534in}}{\pgfqpoint{0.930000in}{0.398276in}}%
\pgfusepath{clip}%
\pgfsetrectcap%
\pgfsetroundjoin%
\pgfsetlinewidth{1.505625pt}%
\definecolor{currentstroke}{rgb}{0.121569,0.466667,0.705882}%
\pgfsetstrokecolor{currentstroke}%
\pgfsetdash{}{0pt}%
\pgfpathmoveto{\pgfqpoint{0.332241in}{1.016309in}}%
\pgfpathlineto{\pgfqpoint{0.426180in}{1.014767in}}%
\pgfpathlineto{\pgfqpoint{0.520120in}{1.012324in}}%
\pgfpathlineto{\pgfqpoint{0.614059in}{1.011596in}}%
\pgfpathlineto{\pgfqpoint{0.707999in}{1.011794in}}%
\pgfpathlineto{\pgfqpoint{0.801938in}{1.018810in}}%
\pgfpathlineto{\pgfqpoint{0.895877in}{1.055150in}}%
\pgfpathlineto{\pgfqpoint{0.989817in}{1.138726in}}%
\pgfpathlineto{\pgfqpoint{1.083756in}{1.243786in}}%
\pgfpathlineto{\pgfqpoint{1.177696in}{1.273675in}}%
\pgfusepath{stroke}%
\end{pgfscope}%
\begin{pgfscope}%
\pgfpathrectangle{\pgfqpoint{0.289968in}{0.999534in}}{\pgfqpoint{0.930000in}{0.398276in}}%
\pgfusepath{clip}%
\pgfsetrectcap%
\pgfsetroundjoin%
\pgfsetlinewidth{1.505625pt}%
\definecolor{currentstroke}{rgb}{1.000000,0.498039,0.054902}%
\pgfsetstrokecolor{currentstroke}%
\pgfsetdash{}{0pt}%
\pgfpathmoveto{\pgfqpoint{0.332241in}{1.002975in}}%
\pgfpathlineto{\pgfqpoint{0.426180in}{1.005155in}}%
\pgfpathlineto{\pgfqpoint{0.520120in}{1.024815in}}%
\pgfpathlineto{\pgfqpoint{0.614059in}{1.023930in}}%
\pgfpathlineto{\pgfqpoint{0.707999in}{1.005539in}}%
\pgfpathlineto{\pgfqpoint{0.801938in}{1.004640in}}%
\pgfpathlineto{\pgfqpoint{0.895877in}{1.047961in}}%
\pgfpathlineto{\pgfqpoint{0.989817in}{1.134380in}}%
\pgfpathlineto{\pgfqpoint{1.083756in}{1.223274in}}%
\pgfpathlineto{\pgfqpoint{1.177696in}{1.293976in}}%
\pgfusepath{stroke}%
\end{pgfscope}%
\begin{pgfscope}%
\pgfsetrectcap%
\pgfsetmiterjoin%
\pgfsetlinewidth{0.803000pt}%
\definecolor{currentstroke}{rgb}{0.000000,0.000000,0.000000}%
\pgfsetstrokecolor{currentstroke}%
\pgfsetdash{}{0pt}%
\pgfpathmoveto{\pgfqpoint{0.289968in}{0.999534in}}%
\pgfpathlineto{\pgfqpoint{0.289968in}{1.397810in}}%
\pgfusepath{stroke}%
\end{pgfscope}%
\begin{pgfscope}%
\pgfsetrectcap%
\pgfsetmiterjoin%
\pgfsetlinewidth{0.803000pt}%
\definecolor{currentstroke}{rgb}{0.000000,0.000000,0.000000}%
\pgfsetstrokecolor{currentstroke}%
\pgfsetdash{}{0pt}%
\pgfpathmoveto{\pgfqpoint{1.219968in}{0.999534in}}%
\pgfpathlineto{\pgfqpoint{1.219968in}{1.397810in}}%
\pgfusepath{stroke}%
\end{pgfscope}%
\begin{pgfscope}%
\pgfsetrectcap%
\pgfsetmiterjoin%
\pgfsetlinewidth{0.803000pt}%
\definecolor{currentstroke}{rgb}{0.000000,0.000000,0.000000}%
\pgfsetstrokecolor{currentstroke}%
\pgfsetdash{}{0pt}%
\pgfpathmoveto{\pgfqpoint{0.289968in}{0.999534in}}%
\pgfpathlineto{\pgfqpoint{1.219968in}{0.999534in}}%
\pgfusepath{stroke}%
\end{pgfscope}%
\begin{pgfscope}%
\pgfsetrectcap%
\pgfsetmiterjoin%
\pgfsetlinewidth{0.803000pt}%
\definecolor{currentstroke}{rgb}{0.000000,0.000000,0.000000}%
\pgfsetstrokecolor{currentstroke}%
\pgfsetdash{}{0pt}%
\pgfpathmoveto{\pgfqpoint{0.289968in}{1.397810in}}%
\pgfpathlineto{\pgfqpoint{1.219968in}{1.397810in}}%
\pgfusepath{stroke}%
\end{pgfscope}%
\begin{pgfscope}%
\pgfsetbuttcap%
\pgfsetmiterjoin%
\definecolor{currentfill}{rgb}{1.000000,1.000000,1.000000}%
\pgfsetfillcolor{currentfill}%
\pgfsetlinewidth{0.000000pt}%
\definecolor{currentstroke}{rgb}{0.000000,0.000000,0.000000}%
\pgfsetstrokecolor{currentstroke}%
\pgfsetstrokeopacity{0.000000}%
\pgfsetdash{}{0pt}%
\pgfpathmoveto{\pgfqpoint{0.289968in}{0.521603in}}%
\pgfpathlineto{\pgfqpoint{1.219968in}{0.521603in}}%
\pgfpathlineto{\pgfqpoint{1.219968in}{0.919879in}}%
\pgfpathlineto{\pgfqpoint{0.289968in}{0.919879in}}%
\pgfpathlineto{\pgfqpoint{0.289968in}{0.521603in}}%
\pgfpathclose%
\pgfusepath{fill}%
\end{pgfscope}%
\begin{pgfscope}%
\pgfsetbuttcap%
\pgfsetroundjoin%
\definecolor{currentfill}{rgb}{0.000000,0.000000,0.000000}%
\pgfsetfillcolor{currentfill}%
\pgfsetlinewidth{0.803000pt}%
\definecolor{currentstroke}{rgb}{0.000000,0.000000,0.000000}%
\pgfsetstrokecolor{currentstroke}%
\pgfsetdash{}{0pt}%
\pgfsys@defobject{currentmarker}{\pgfqpoint{0.000000in}{-0.048611in}}{\pgfqpoint{0.000000in}{0.000000in}}{%
\pgfpathmoveto{\pgfqpoint{0.000000in}{0.000000in}}%
\pgfpathlineto{\pgfqpoint{0.000000in}{-0.048611in}}%
\pgfusepath{stroke,fill}%
}%
\begin{pgfscope}%
\pgfsys@transformshift{0.332241in}{0.521603in}%
\pgfsys@useobject{currentmarker}{}%
\end{pgfscope}%
\end{pgfscope}%
\begin{pgfscope}%
\definecolor{textcolor}{rgb}{0.000000,0.000000,0.000000}%
\pgfsetstrokecolor{textcolor}%
\pgfsetfillcolor{textcolor}%
\pgftext[x=0.332241in,y=0.424381in,,top]{\color{textcolor}{\sffamily\fontsize{10.000000}{12.000000}\selectfont\catcode`\^=\active\def^{\ifmmode\sp\else\^{}\fi}\catcode`\%=\active\def
\end{pgfscope}%
\begin{pgfscope}%
\pgfsetbuttcap%
\pgfsetroundjoin%
\definecolor{currentfill}{rgb}{0.000000,0.000000,0.000000}%
\pgfsetfillcolor{currentfill}%
\pgfsetlinewidth{0.803000pt}%
\definecolor{currentstroke}{rgb}{0.000000,0.000000,0.000000}%
\pgfsetstrokecolor{currentstroke}%
\pgfsetdash{}{0pt}%
\pgfsys@defobject{currentmarker}{\pgfqpoint{0.000000in}{-0.048611in}}{\pgfqpoint{0.000000in}{0.000000in}}{%
\pgfpathmoveto{\pgfqpoint{0.000000in}{0.000000in}}%
\pgfpathlineto{\pgfqpoint{0.000000in}{-0.048611in}}%
\pgfusepath{stroke,fill}%
}%
\begin{pgfscope}%
\pgfsys@transformshift{0.707999in}{0.521603in}%
\pgfsys@useobject{currentmarker}{}%
\end{pgfscope}%
\end{pgfscope}%
\begin{pgfscope}%
\definecolor{textcolor}{rgb}{0.000000,0.000000,0.000000}%
\pgfsetstrokecolor{textcolor}%
\pgfsetfillcolor{textcolor}%
\pgftext[x=0.707999in,y=0.424381in,,top]{\color{textcolor}{\sffamily\fontsize{10.000000}{12.000000}\selectfont\catcode`\^=\active\def^{\ifmmode\sp\else\^{}\fi}\catcode`\%=\active\def
\end{pgfscope}%
\begin{pgfscope}%
\pgfsetbuttcap%
\pgfsetroundjoin%
\definecolor{currentfill}{rgb}{0.000000,0.000000,0.000000}%
\pgfsetfillcolor{currentfill}%
\pgfsetlinewidth{0.803000pt}%
\definecolor{currentstroke}{rgb}{0.000000,0.000000,0.000000}%
\pgfsetstrokecolor{currentstroke}%
\pgfsetdash{}{0pt}%
\pgfsys@defobject{currentmarker}{\pgfqpoint{0.000000in}{-0.048611in}}{\pgfqpoint{0.000000in}{0.000000in}}{%
\pgfpathmoveto{\pgfqpoint{0.000000in}{0.000000in}}%
\pgfpathlineto{\pgfqpoint{0.000000in}{-0.048611in}}%
\pgfusepath{stroke,fill}%
}%
\begin{pgfscope}%
\pgfsys@transformshift{1.177696in}{0.521603in}%
\pgfsys@useobject{currentmarker}{}%
\end{pgfscope}%
\end{pgfscope}%
\begin{pgfscope}%
\definecolor{textcolor}{rgb}{0.000000,0.000000,0.000000}%
\pgfsetstrokecolor{textcolor}%
\pgfsetfillcolor{textcolor}%
\pgftext[x=1.177696in,y=0.424381in,,top]{\color{textcolor}{\sffamily\fontsize{10.000000}{12.000000}\selectfont\catcode`\^=\active\def^{\ifmmode\sp\else\^{}\fi}\catcode`\%=\active\def
\end{pgfscope}%
\begin{pgfscope}%
\definecolor{textcolor}{rgb}{0.000000,0.000000,0.000000}%
\pgfsetstrokecolor{textcolor}%
\pgfsetfillcolor{textcolor}%
\pgftext[x=0.754968in,y=0.234413in,,top]{\color{textcolor}{\sffamily\fontsize{10.000000}{12.000000}\selectfont\catcode`\^=\active\def^{\ifmmode\sp\else\^{}\fi}\catcode`\%=\active\def
\end{pgfscope}%
\begin{pgfscope}%
\definecolor{textcolor}{rgb}{0.000000,0.000000,0.000000}%
\pgfsetstrokecolor{textcolor}%
\pgfsetfillcolor{textcolor}%
\pgftext[x=0.234413in,y=0.720741in,,bottom,rotate=90.000000]{\color{textcolor}{\sffamily\fontsize{10.000000}{12.000000}\selectfont\catcode`\^=\active\def^{\ifmmode\sp\else\^{}\fi}\catcode`\%=\active\def
\end{pgfscope}%
\begin{pgfscope}%
\pgfpathrectangle{\pgfqpoint{0.289968in}{0.521603in}}{\pgfqpoint{0.930000in}{0.398276in}}%
\pgfusepath{clip}%
\pgfsetrectcap%
\pgfsetroundjoin%
\pgfsetlinewidth{1.505625pt}%
\definecolor{currentstroke}{rgb}{0.121569,0.466667,0.705882}%
\pgfsetstrokecolor{currentstroke}%
\pgfsetdash{}{0pt}%
\pgfpathmoveto{\pgfqpoint{0.332241in}{0.646940in}}%
\pgfpathlineto{\pgfqpoint{0.426180in}{0.647395in}}%
\pgfpathlineto{\pgfqpoint{0.520120in}{0.647846in}}%
\pgfpathlineto{\pgfqpoint{0.614059in}{0.647880in}}%
\pgfpathlineto{\pgfqpoint{0.707999in}{0.647757in}}%
\pgfpathlineto{\pgfqpoint{0.801938in}{0.647679in}}%
\pgfpathlineto{\pgfqpoint{0.895877in}{0.646997in}}%
\pgfpathlineto{\pgfqpoint{0.989817in}{0.647478in}}%
\pgfpathlineto{\pgfqpoint{1.083756in}{0.647718in}}%
\pgfpathlineto{\pgfqpoint{1.177696in}{0.647798in}}%
\pgfusepath{stroke}%
\end{pgfscope}%
\begin{pgfscope}%
\pgfpathrectangle{\pgfqpoint{0.289968in}{0.521603in}}{\pgfqpoint{0.930000in}{0.398276in}}%
\pgfusepath{clip}%
\pgfsetrectcap%
\pgfsetroundjoin%
\pgfsetlinewidth{1.505625pt}%
\definecolor{currentstroke}{rgb}{1.000000,0.498039,0.054902}%
\pgfsetstrokecolor{currentstroke}%
\pgfsetdash{}{0pt}%
\pgfpathmoveto{\pgfqpoint{0.332241in}{0.644270in}}%
\pgfpathlineto{\pgfqpoint{0.426180in}{0.659571in}}%
\pgfpathlineto{\pgfqpoint{0.520120in}{0.682583in}}%
\pgfpathlineto{\pgfqpoint{0.614059in}{0.662535in}}%
\pgfpathlineto{\pgfqpoint{0.707999in}{0.632122in}}%
\pgfpathlineto{\pgfqpoint{0.801938in}{0.614825in}}%
\pgfpathlineto{\pgfqpoint{0.895877in}{0.621003in}}%
\pgfpathlineto{\pgfqpoint{0.989817in}{0.636943in}}%
\pgfpathlineto{\pgfqpoint{1.083756in}{0.645510in}}%
\pgfpathlineto{\pgfqpoint{1.177696in}{0.660616in}}%
\pgfusepath{stroke}%
\end{pgfscope}%
\begin{pgfscope}%
\pgfsetrectcap%
\pgfsetmiterjoin%
\pgfsetlinewidth{0.803000pt}%
\definecolor{currentstroke}{rgb}{0.000000,0.000000,0.000000}%
\pgfsetstrokecolor{currentstroke}%
\pgfsetdash{}{0pt}%
\pgfpathmoveto{\pgfqpoint{0.289968in}{0.521603in}}%
\pgfpathlineto{\pgfqpoint{0.289968in}{0.919879in}}%
\pgfusepath{stroke}%
\end{pgfscope}%
\begin{pgfscope}%
\pgfsetrectcap%
\pgfsetmiterjoin%
\pgfsetlinewidth{0.803000pt}%
\definecolor{currentstroke}{rgb}{0.000000,0.000000,0.000000}%
\pgfsetstrokecolor{currentstroke}%
\pgfsetdash{}{0pt}%
\pgfpathmoveto{\pgfqpoint{1.219968in}{0.521603in}}%
\pgfpathlineto{\pgfqpoint{1.219968in}{0.919879in}}%
\pgfusepath{stroke}%
\end{pgfscope}%
\begin{pgfscope}%
\pgfsetrectcap%
\pgfsetmiterjoin%
\pgfsetlinewidth{0.803000pt}%
\definecolor{currentstroke}{rgb}{0.000000,0.000000,0.000000}%
\pgfsetstrokecolor{currentstroke}%
\pgfsetdash{}{0pt}%
\pgfpathmoveto{\pgfqpoint{0.289968in}{0.521603in}}%
\pgfpathlineto{\pgfqpoint{1.219968in}{0.521603in}}%
\pgfusepath{stroke}%
\end{pgfscope}%
\begin{pgfscope}%
\pgfsetrectcap%
\pgfsetmiterjoin%
\pgfsetlinewidth{0.803000pt}%
\definecolor{currentstroke}{rgb}{0.000000,0.000000,0.000000}%
\pgfsetstrokecolor{currentstroke}%
\pgfsetdash{}{0pt}%
\pgfpathmoveto{\pgfqpoint{0.289968in}{0.919879in}}%
\pgfpathlineto{\pgfqpoint{1.219968in}{0.919879in}}%
\pgfusepath{stroke}%
\end{pgfscope}%
\end{pgfpicture}%
\makeatother%
\endgroup%